\documentclass{article}
\usepackage[utf8]{inputenc}
\usepackage{latexsym}
\usepackage[normalem]{ulem}

\usepackage{authblk}
\usepackage[top=0.8in, bottom=1in, left=0.6in, right=0.6in]{geometry}
 
\usepackage[outline]{contour}
\usepackage{comment}
\usepackage{graphicx}
\usepackage{color}
\usepackage{amsmath, amsthm, amssymb}
\usepackage{enumerate}
\usepackage{float}
\usepackage{url}
\usepackage{stackrel}
\usepackage{mathrsfs,dsfont}
\usepackage{caption}
\usepackage{wrapfig}
\usepackage{bbm}
\usepackage{bm}
\usepackage{epigraph}
\usepackage{tikz-qtree}
\usepackage{tikz}
\usepackage{indentfirst}
\usepackage{float}
\usepackage{mathtools}
\usepackage{subfigure}
\usepackage{tabularx}
\usepackage{booktabs}

\usepackage{lipsum}

\usepackage{abstract}

\newcommand\blfootnote[1]{%
  \begingroup
  \renewcommand\thefootnote{}\footnote{#1}%
  \addtocounter{footnote}{-1}%
  \endgroup
}

\newtheorem{theorem}{Theorem}[section]

\newtheorem{remark}[theorem]{Remark}

\newcommand{\RNum}[1]
{\uppercase\expandafter{\romannumeral #1\relax}}
\title{A Mathematical Model of Clonal Hematopoiesis Explaining Phase Transitions in Chronic Myeloid Leukemia}
\author[1,2,3,*]{Lorand Gabriel Parajdi}
\author[1]{Xue Bai}
\author[3,4,5]{Dávid Kegyes}
\author[3,4,5]{Ciprian Tomuleasa}
\affil[1]{Department of Mathematics, West Virginia University, 26506 Morgantown, WV, USA}
\affil[2]{Department of Mathematics, ``Babeş–Bolyai" University, 400084 Cluj-Napoca, Romania}
\affil[3]{Academy of Romanian Scientists, Ilfov 3, 050044 Bucharest, Romania}
\affil[4]{Department of Hematology, ``Ion Chiricuță" Clinical Cancer Center, 400124 Cluj-Napoca, Romania}
\affil[5]{Research Center for Advanced Medicine, ``Iuliu Hațieganu" University of Medicine and Pharmacy, 400349 Cluj-Napoca, Romania}
\date{}                  

\begin{document}
\maketitle

\blfootnote{*Corresponding author: lorand.parajdi@ubbcluj.ro}
\blfootnote{Email: lorand.parajdi@ubbcluj.ro (L.G. Parajdi),\  xb0002@mix.wvu.edu (X. Bai),\  david.vale.kegyes@elearn.umfcluj.ro (D. Kegyes),}
\blfootnote{ciprian.tomuleasa@umfcluj.ro (C. Tomuleasa). This version of the manuscript was published in Mathematical Medicine and Biology: A} \blfootnote{Journal of the IMA (2025), 42(3), 253–288. https://doi.org/10.1093/imammb/dqaf004}
\vspace{-0.8cm}

\textbf{Abstract:} This study presents a mathematical model describing cloned hematopoiesis in chronic myeloid leukemia (CML) through a nonlinear system of differential equations. The primary objective is to understand the progression from healthy hematopoiesis to the chronic and accelerated-acute phases in myeloid leukemia. The model incorporates intrinsic cellular division events in hematopoiesis and delineates the evolution of chronic myeloid leukemia into five compartments: cycling stem cells, quiescent stem cells, progenitor cells, differentiated cells and terminally differentiated cells. Our analysis reveals the existence of three distinct non-zero steady states within the dynamical system, representing healthy hematopoiesis, the chronic phase and the accelerated-acute stage of the disease. We investigate the local and global stability of these steady states and provide a characterization of the hematopoietic states based on this analysis. Additionally, numerical simulations are included to illustrate the theoretical results. \\


\section{Introduction}

\noindent Mathematics provides both qualitative and quantitative tools to enhance our understanding, prediction and management of biological processes. These tools and models have proven invaluable in studying blood cell production processes and hematological disorders. Early contributions in this direction include the works of Rubinow and Lebowitz \cite{rl1,rl2}, Mackey and Glass \cite{mg}, Mackey \cite{m1} and Djulbegovic and Svetina \cite{ds}. Recent advancements are evident in publications by Fokas \textit{et al.} \cite{f}, Neiman \cite{n}, Andersen and Mackey \cite{am}, Colijn and Mackey \cite{cm}, Adimy \textit{et al.} \cite{acr}, Dingli and Michor \cite{dm}, Kim \textit{et al.} \cite{kll}, Cucuianu and Precup \cite{cp}, Doumic-Jauffret \textit{et al.} \cite{dj}, Komarova \cite{ko}, Stiehl and Marciniak-Czochra \cite{s}, MacLean \textit{et al.} \cite{mls,mfs}, Radulescu \textit{et al.} \cite{rch}, Bianca \textit{et al.} \cite{bpmr,bppr}, Ragusa and Russo \cite{rr}, Badralexi \textit{et al.} \cite{bhm} and the references therein. Additionally, studies related to stem cell transplantation have been conducted by Vincent \textit{et al.} \cite{vi}, DeConde \textit{et al.} \cite{dkll}, Kim \textit{et al.} \cite{K}, Marciniak-Czochra and Stiehl \cite{ms}, Precup \textit{et al.} \cite{pacs,pst,pdtsp}, Stiehl \textit{et al.} \cite{sh1} and Parajdi \cite{plg}. Furthermore, valuable insights into mathematical models for cancer, with a focus on chronic myeloid leukemia, can be found in reviews by Afenya \cite{af}, Michor \cite{mi}, Foley and Mackey \cite{fm} and Clapp and Levy \cite{cld}.

Cycling and quiescent hematopoietic stem cells within the bone marrow play a pivotal role in the intricate process of cell production, possessing the unique capability for self-renewal and the potential to differentiate into various blood cell types. Disruptions in this intricate biological mechanism underlie the development of several hematological disorders, including chronic myeloid leukemia (CML). Chronic myeloid leukemia is a malignancy originating from mutated stem cells, affecting the myeloid cell lineage and progressing through three distinct phases: the chronic phase, the accelerated or transitional phase and the acute or blast phase. Due to the challenge in distinguishing between the last two phases, they will be collectively referred to as the accelerated-acute phase in this context.

In this paper, we propose a mathematical model of cloned hematopoiesis in chronic myeloid leukemia, inspired by the work of Molina-Peña \textit{et al.} \cite{mta}, which is based on pseudo-chemical reactions. We will construct this model starting from the healthy-leukemic stem cell system introduced by Parajdi \textit{et al.} \cite{ppbt}, which aids in understanding and describing the progression of chronic myeloid leukemia. The model enables us to differentiate among three hematopoietic conditions associated with chronic myeloid leukemia: the healthy hematopoietic state, the chronic leukemic state and the accelerated-acute state, also known as the blast phase. In our mathematical-biological study, we introduce a novel model represented by a ten-dimensional nonlinear differential system. This model exhibits three distinct nontrivial steady states and their local and global stability analysis reveals that the asymptotic stability of these states depends on a cumulative parameter $R$. This parameter encapsulates factors such as the intrinsic reaction rate of symmetric self-renewal or the growth rate, cell death rate, the intrinsic reaction rate of the direct and symmetric differentiation and the sensitivity rate of malignant cells, representing the equilibrium amount of the malignant cells. Values of $R$ below a specific threshold correspond to the healthy hematological state, while $R$ values within a defined interval characterize the chronic phase. Larger $R$ values indicate the accelerated-acute phase, as depicted in Figure \ref{fig:2} below. Mathematically, the transition between different hematological states is driven by variations in the fundamental parameters of leukemic cells aggregated within $R$. The transition is progressive as the disease advances and regressive with the aid of treatment. In the context of therapeutic advancements, our analysis may serve as a valuable guide for improving therapeutic agents and strategies.

The paper is structured as follows: In Section \ref{Section:2}, we propose a new model and proceed with its mathematical analysis. The model comprehensively describes cell divisions (symmetric division - self-renewing and asymmetric division), differentiations (direct differentiation and symmetric differentiation) and transitions (between active/cycling phase and resting/quiescent phase, applicable only to stem cells) on five levels: cycling stem cells, quiescent stem cells, progenitor cells, differentiated cells and terminally differentiated cells. Based on this model, the healthy hematopoietic state and the chronic and accelerated-acute phases of chronic myeloid leukemia (CML) are mathematically characterized in terms of parameter $R$. We provide comprehensive proofs establishing the local and global stability of the observed steady states. In Section \ref{Section:3}, we provide some numerical simulations of the model along with parameter estimation. Finally, in Section \ref{Section:4}, we include medical discussions and conclusions. Therefore, we conclude this introductory section by providing a medical background for readers seeking more in-depth biological and medical insights into hematopoiesis, malignant disorders and relevant literature.

\subsection{Medical Background}

\noindent \textbf{Hematopoiesis} is the complex process through which the body produces blood cells by differentiation of hematopoietic stem cells (HSCs) into various lineages, including red blood cells (erythrocytes), white blood cells (leukocytes) and platelets (thrombocytes). This is an active, lifetime process primarily occurring in the bone marrow, which plays a crucial role in regulating HSC lineage commitment (Pinho and Frenette \cite{pf}). Before HSCs commit to a specific lineage (such as erythroid, megakaryocytic or lymphoid), they undergo differentiation through a series of progenitor cell stages, gradually losing their potential to develop into multiple lineages. These stem and progenitor cells can be identified biologically either by specific cell surface markers using flow cytometry or by functional assays assessing key stem cell characteristics, such as self-renewal capacity, clonogenicity and the ability to differentiate into a particular lineage (Savino \textit{et al.} \cite{ssd}). The continuous regeneration of the hematopoietic system requires the precise generation of specific cell types at the right time and location. The process of clonal hematopoietic cell production, known as steady-state hematopoiesis, occurs in the bone marrow of healthy adults, where mature cells fully develop from blast cells. All stages of this transformation (from the colony-forming unit (CFU) to the differentiated state) take place within the hematopoietic niche, a specialized microenvironment composed of a central macrophage surrounded by developing differentiating cells. To maintain normal hematopoiesis, HSCs must consistently make accurate fate decisions, balancing factors such as quiescence versus proliferation, self-renewal versus differentiation and choices related to survival versus programmed cell death (apoptosis). A complex network of physiological communication pathways regulates the generation, distribution and turnover of differentiated blood cells in healthy individuals, ensuring the maintenance of normal clonal hematopoiesis throughout lifetime. Disruptions in these healthy cell fate decisions can lead to hematological disorders, with specific molecular pathways governing both physiological and pathological processes (Jagannathan-Bogdan and Zon \cite{jbz}).

As \textbf{hematopoietic stem cells} (HSCs) are enduring cells, the accumulation of somatic mutations over time is inevitable with aging. This accumulation can lead to the emergence of mutations that confer a survival advantage to the mutated clone, a process referred to as clonal hematopoiesis. Clonal hematopoiesis is implicated in various hematological malignancies, including leukemias, lymphomas, myeloproliferative neoplasms and multiple myeloma. However, mutations associated with these cancers are also found in a significant proportion of non-cancerous elderly individuals. To distinguish these pathogenic mutations from benign ones, the clinical term `clonal hematopoiesis of indeterminate potential' (CHIP) was introduced (Jaiswal \cite{j}). Individuals with CHIP carry a cancer-associated somatic mutation in their blood or bone marrow hematopoietic cells but are asymptomatic and are not diagnosed with a hematologic cancer. However, the presence of CHIP significantly increases the risk of developing myelodysplastic syndromes and leukemias (Mitchell \textit{et al.} \cite{mgj}). In some pathological conditions, the regulatory network of clonal hematopoiesis may become overwhelmed or dysfunctional, leading to diseases such as myelodysplasia or leukemia. These clonal bone marrow diseases, particularly in the elderly, are characterized by chronic cytopenias and morphologic dysplasia of hematopoietic cells, carrying a high risk of progression to acute myeloid leukemia. Ineffective clonal hematopoiesis, as seen in myelodysplasia or leukemia, often arises from an inflammatory environment that induces malignant clonal alterations. These alterations are driven by chromosomal abnormalities (e.g. del(5q), del(7q)) deletions/additions, specific mutations of the spliceosome (e.g. SF3B1, SRSF2), transcription factors (e.g. RUNX1, ETV6), NLRP3 inflammasome activation, overexpressed SMAD2/3 downstream mediators, TGF-$\beta$ signaling, epigenetic modifiers (e.g. TET2, DNMT3A 5, IDH1/2, ASXL1) and RNA splicing factors, all of which contribute to MDS pathogenesis by inhibiting cell maturation. The complexity of these genetic alterations makes it challenging to interpret them mathematically, as each patient with altered clonal hematopoiesis is different. Still, in this manuscript, we aim to propose a mathematical model for altered clonal hematopoiesis, in order to bring the field closer to a better understanding of leukemogenesis.

This paper focuses on leukemias, with a particular emphasis on chronic myeloid leukemia (CML). \textbf{Leukemias} are a group of hematologic malignancies characterized by the malignant proliferation of immature blood cells in the bone marrow. This uncontrolled growth impedes the maturation of healthy hematopoietic cells, leading to an accumulation of dysfunctional and often undifferentiated leukocytes. Leukemias can manifest in acute or chronic forms, depending on the rate of disease progression and the maturity of the leukemic cells. They are broadly classified into four main types: acute myeloid leukemia (AML), acute lymphoblastic leukemia (ALL), chronic myeloid leukemia (CML) and chronic lymphocytic leukemia (CLL) (Arber \textit{et al.} \cite{aoh}). Acute leukemias involve rapidly proliferating, immature cells that require immediate treatment, while chronic leukemias are characterized by a slower accumulation of mature but malignant cells. Approximately $30\%$ of leukemia cases are classified as CML (Lin \textit{et al.} \cite{lms}). The clinical and biological manifestations of CML are caused by a specific mutation: the $t(9;22)(q34;q11)$ translocation, which results in a truncated chromosome $22$ known as the Philadelphia chromosome. This translocation fuses the BCR gene from chromosome $22$ with the ABL1 gene from chromosome $9$, creating the BCR-ABL1 fusion gene. The resultant BCR-ABL1 RNA transcripts are translated into fusion proteins with elevated tyrosine kinase activity, conferring a survival advantage to leukemic cells and leading to clonal dominance in both the bone marrow and peripheral blood. With 5-year overall survival rates ranging from $85\%$ to $95\%$, the life expectancy of CML patients is now nearly equivalent to that of the general population across all ages (Breccia \textit{et al.} \cite{bes}). This remarkable improvement is largely due to the advent of tyrosine kinase inhibitors (TKIs), as well as advances in diagnostic techniques that enable early diagnosis and detection of relapse through the measurement of BCR-ABL1 transcript levels (Arbore \textit{et al.} \cite{agd}). However, not all patients respond well to therapy. CML can be medically classified into three phases: chronic, accelerated and acute or blast phase. While the chronic phase is the most common form and the majority of patients remain in this phase, approximately $1-1.5\%$ per year of these patients progress to the acute/blast phase, which is clinically and prognostically more severe than acute leukemias (Hehlmann \cite{h}). Most of these patients succumb within the first two months of diagnosis. Therefore, predicting who will progress to the blast phase and who will respond to treatment is crucial in clinical practice. Currently, a universally accepted predictive model does not exist. To develop a highly predictive model, it is essential to understand and mathematically describe the clonal hematopoiesis that leads to CML. This knowledge is critical for advancing clinical artificial intelligence-based predictive models.

\subsection{The Healthy-Leukemic Stem Cell System}

\noindent Mackey and Glass \cite{mg} introduced a mathematical model for blood production, employing differential equations that incorporate sigmoid or Hill functions to account for the self-regulating nature of hematopoiesis. Building upon this framework, Parajdi \textit{et al.} \cite{ppbt} formulated the following mathematical model:
\begin{equation}
\begin{cases}
x^{\prime }(t)=\frac{ax(t)}{1+b_{1}x(t)+b_{2}y(t)}-c x(t) \\[2mm]
y^{\prime }(t)=\frac{Ay(t)}{1+B( x(t)+y(t)) }-C y(t)\ .%
\end{cases}
\label{s}
\end{equation}

Here, $x(t)$ and $y(t)$ represent the populations of healthy and malignant/leukemic HSCs at time $t$, respectively. The parameters $a$ and $A$ correspond to the growth rates, while $c$ and $C$ signify the cell death rates or apoptosis rates. Additionally, the sensitivity parameters $b_1$, $b_2$ and $B$ play a crucial role in governing the self-regulation process. The complexity of biochemical and biophysical processes at the cellular level is represented in the mathematical model by the parameters $a$, $b_{1}$, $b_{2}$, $c$, $A$, $B$ and $C$, whose determination requires laboratory and clinical investigations (see Michor \cite{mi}, Michor \textit{et al.} \cite{mhi} and Foo \textit{et al.} \cite{fdc}). These mathematical models use a single Hill-type feedback on stem cells, a simplification that highlights the critical role of niche-mediated feedback in hematopoiesis and leukemogenesis, as supported by Dingli and Michor \cite{dm}, Cucuianu and Precup \cite{cp}, Michor \textit{et al.} \cite{mhi} and Foo \textit{et al.} \cite{fdc}. Future studies integrating additional feedback mechanisms for other cell types (e.g. progenitor or differentiated cells) could further refine and expand these models.

It is assumed that, for both stem cell populations, the growth rate exceeds the death rate, denoted as $a > c$ and $A > C$. Furthermore, the relative advantage of leukemic stem cells, characterized by their reduced sensitivity to the bone marrow microenvironment compared to healthy stem cells, is expressed through the relationships $b_1 \geq b_2 > B$. When $b_{1}=b_{2}$, this scenario was investigated by Dingli and Michor \cite{dm} (also explored by Cucuianu and Precup \cite{cp}), in order to describe the time competition between healthy and malignant HSCs. In this particular scenario, the system features only two non-zero steady states, represented as $E_{1}(d, 0)$ and $E_{2}(0, D)$, where $d$ and $D$ denote the homeostatic amounts of healthy and malignant stem cells, determined as follows:
\begin{equation*}
d = \frac{1}{b_{1}}\left(\frac{a}{c}-1\right)\ \ \ \text{and}\ \ \ D = \frac{
1}{B}\left(\frac{A}{C}-1\right).
\end{equation*}

The scenario where $b_1 > b_2$ was initially introduced and analyzed in \cite{p, ppbt}. In this context, the possibility of a third steady state $E_{3}(x^{\ast}, y^{\ast})$ emerges, where:\ \ $x^{\ast} = \frac{b_{1}}{b_{1}-b_{2}}d - \frac{b_{2}}{b_{1}-b_{2}}D$\ \ and\ \ $y^{\ast} = \frac{b_{1}}{b_{1}-b_{2}}(D-d)$.
Under the assumption of $b_1 > b_2$, we observe that both $x^{\ast}$ and $y^{\ast}$ are positive if and only if $d < D < (b_1/b_2)d$. The scenarios are delineated as follows: when $D < d$, it corresponds to the healthy hematologic state; for $d < D < (b_1/b_2)d$, we are in the chronic phase of the disease; and $(b_1/b_2)d < D$ signifies the acute/blast phase. In the healthy state (first scenario), the system (\ref{s}) exhibits a unique stable equilibrium denoted as $E_{1}(d, 0)$. Transitioning to the chronic phase (second scenario), the stable equilibrium takes the form of $E_{3}(x^{\ast}, y^{\ast})$, where both $x^{\ast}$ and $y^{\ast}$ are positive. In the acute/blast phase (third scenario), the stable equilibrium is represented by $E_{2}(0, D)$.

\section{The Mathematical Model} \label{Section:2}

\noindent In this section, we will introduce a mathematical model inspired by the works of Molina-Peña \textit{et al.} \cite{mta} and Parajdi \textit{et al.} \cite{ppbt}. This model is designed to aid in understanding and describing the progression of chronic myeloid leukemia.
\vspace{-0.4cm}
\begin{figure}[H]
  \centering
\includegraphics[width=0.75\textwidth, height = 0.55\textwidth]{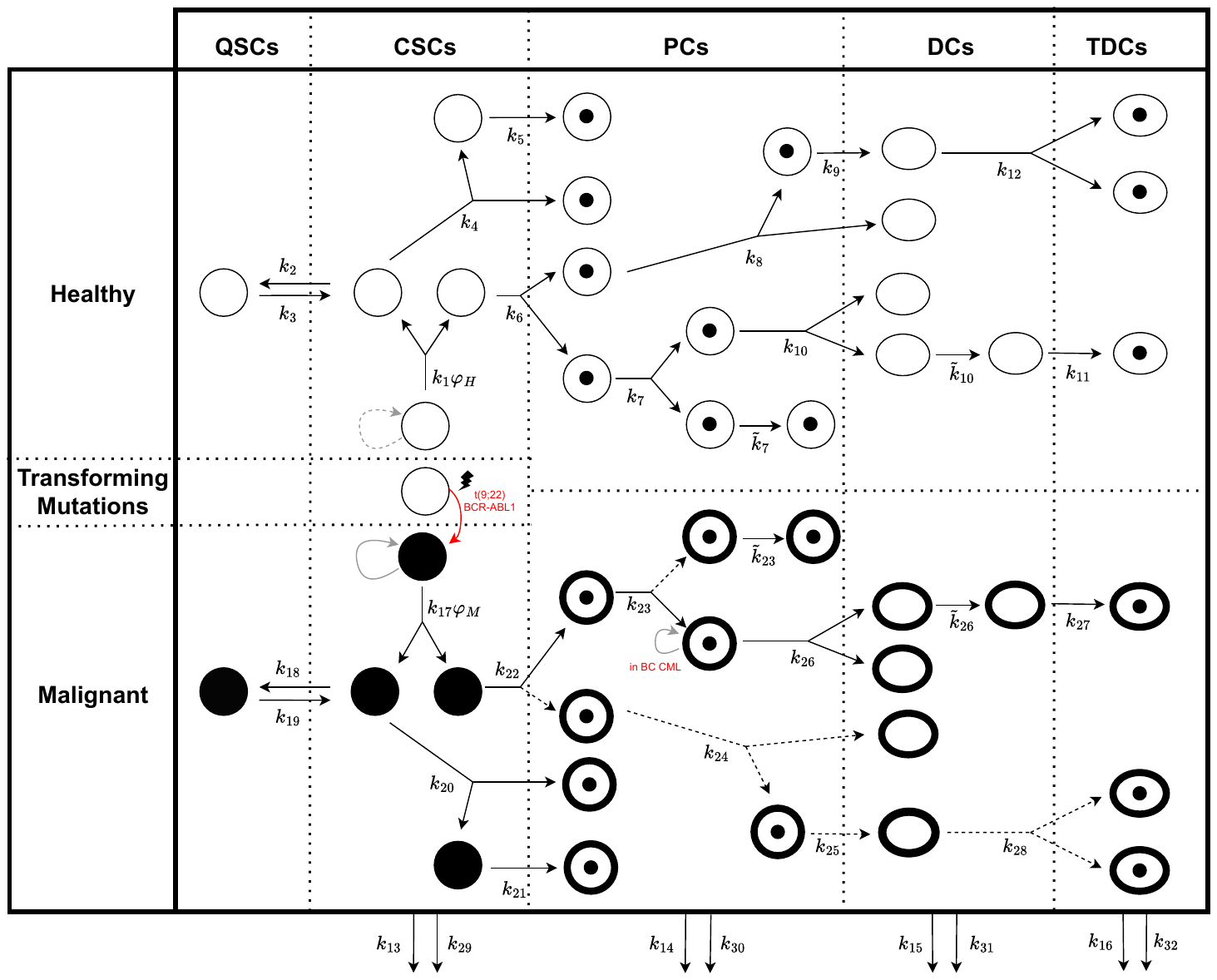}
\vspace{-0.6cm}
  \caption{\small{The proposed model is represented using compartments, where a specific rate constant influences each cellular event denoted as $k_{j}$. Both healthy and malignant cycling stem cells (HCSCs/MCSCs) possess the ability to self-renew, which is indicated by the rate constants multiplied by the corresponding self-regulatory functions ($k_{1}\varphi_{H})/(k_{17}\varphi_{M}$), respectively. Additionally, they can enter a resting phase with rate constants ($k_{2})/(k_{18}$), during which they become healthy/malignant quiescent stem cells (HQSCs/MQSCs). After a period of time, the quiescent stem cells can reactivate, denoted by rate constants ($k_{3})/(k_{19}$), and return to an active state as healthy/malignant cycling stem cells. The healthy/malignant cycling stem cells can give rise to intermediate healthy/malignant progenitor cells (HPCs/MPCs) through asymmetric ($k_{4})/(k_{20}$) division and symmetric ($k_{6})/(k_{22}$) differentiation. Additionally, they can undergo direct differentiation ($k_{5})/(k_{21}$) into intermediate healthy/malignant progenitor cells. These intermediate healthy/malignant progenitor cells, in turn, can proliferate and give rise to later healthy/malignant progenitor cells (HPCs/MPCs) through symmetric ($k_{7})/(k_{23}$) and asymmetric ($k_{8})/(k_{24}$) divisions. Similarly, they can give rise to healthy/malignant differentiated cells (HDCs/MDCs) through asymmetric ($k_{8}$)/($k_{24}$) division. Healthy/malignant later progenitor cells can directly differentiate ($\widetilde{k}_{7})/(\widetilde{k}_{23}$) into other types of later healthy/malignant progenitor cells. Additionally, they can undergo direct differentiation ($k_{9}$)/($k_{25}$) into healthy/malignant differentiated cells, and they can give rise to healthy/malignant differentiated cells through symmetric ($k_{10})/(k_{26}$) differentiation. The healthy/malignant differentiated cells can undergo direct differentiation ($\widetilde{k}_{10})/(\widetilde{k}_{26}$) into other types of healthy/malignant differentiated cells. Furthermore, they can give rise to healthy/malignant terminally differentiated cells (HTDCs/MTDCs) through direct ($k_{11})/(k_{27}$) and symmetric ($k_{12})/(k_{28}$) differentiation. Each cell type, except for quiescent stem cells, whether healthy or malignant, has a death rate represented by the rate constants ($k_{13})/(k_{29}$) for cycling stem cells, ($k_{14})/(k_{30}$) for progenitor cells, ($k_{15})/(k_{31}$) for differentiated cells and ($k_{16})/(k_{32}$) for terminally differentiated cells.
}}
  \label{fig:mygraph}
\end{figure}

It allows us to explore the transition process from healthy hematopoiesis to the chronic and accelerated-acute stages of myeloid leukemia. In the following, we will consider the cellular division events involved in hematopoiesis, as illustrated in Figure \ref{fig:mygraph}, as pseudo-chemical reactions mediated by the rate constants: $k_{j}$, where $j=1,\ldots,32$ and $\widetilde{k}_{j}$, where $j=7,10,23,26$.

\textbf{Clinical interpretation} of Figure \ref{fig:mygraph}: all of the cell types outlined in the proposed mathematical model are based on hematopoiesis, which allows each cellular compartment to be recognized and isolated using cell surface markers (flow cytometry) or functional analysis (clonogenic assays). Biologically, healthy hematopoiesis can be disrupted by the acquisition of the $t(9;22)$ mutation, which results in the BCR::ABL1 transcript and fusion protein with strong tyrosine kinase activity leading to survival advantage. This is typical of the chronic phase of CML. Further genetic mutations at the stem/progenitor compartment level can lead to disease progression to blast-phase CML through an intermediate accelerated-phase CML. The accelerated and blast phases are characterized by a block in differentiation and the proliferation of malignant, early progenitor cells known as blasts. The percentage of blasts in the bone marrow distinguishes the two phases ($>20\%$ in the acute/blast phase, $10\%$ to $20\%$ in the accelerated phase). \\

The healthy cycling stem cells (HCSCs) can increase in number through symmetric self-renewal (see Morrison and Kimble \cite{mk}, Dingli \textit{et al.} \cite{dtm}, and Barile \textit{et al.} \cite{bbf}):
\begin{equation}
  \label{R1}
  HCSC \xrightarrow{\text{$k_{1}\varphi_{H}$}} 2HCSC \tag{$R_{1}$},
\end{equation}
where $\varphi_{H} = \frac{1}{1+b_{1}HCSC+b_{2}MCSC}$. \\

The healthy cycling stem cells (HCSCs) can also enter a resting phase called the $G_{0}$ phase, during which they transition into healthy quiescent stem cells (HQSCs). After a certain period of time, the quiescent stem cells can reactivate and return to the active state as healthy cycling stem cells (see Mackey \cite{m}, Foo \textit{et al.} \cite{fdc}, Fuchs and Chen \cite{fc}):

\begin{equation}
  \label{R2}
  HCSC \underset{k_{3}}{\stackrel{k_{2}}{\rightleftarrows}} HQSC \tag{$R_{2}$}.
\end{equation}

Healthy cycling stem cells (HCSCs) can also undergo asymmetric division, giving rise to one healthy cycling stem cell (HCSC), thereby maintaining their numbers, and one intermediate healthy progenitor cell (HPC) destined for differentiation (see Morrison and Kimble \cite{mk}, Dingli \textit{et al.} \cite{dtm}, and Barile \textit{et al.} \cite{bbf}):

\begin{equation}
  \label{R3}
  HCSC \xrightarrow{\text{$k_{4}$}} HCSC + HPC \tag{$R_{3}$}.
\end{equation}

Healthy cycling stem cells (HCSCs) can also directly differentiate into intermediate healthy progenitor cells (HPCs) (see Cheng \textit{et al.} \cite{czc}, and Barile \textit{et al.} \cite{bbf}) and undergo symmetric differentiation, giving rise to two intermediate healthy progenitor cells (HPCs) (see Morrison and Kimble \cite{mk}, Dingli \textit{et al.} \cite{dtm}, and Barile \textit{et al.} \cite{bbf}):
\begin{equation}
  \label{R4}
  HCSC \xrightarrow{\text{$k_{5}$}} HPC \tag{$R_{4}$}.
\end{equation}
\begin{equation}
  \label{R5}
  HCSC \xrightarrow{\text{$k_{6}$}} 2HPC \tag{$R_{5}$}.
\end{equation}

Intermediate healthy progenitor cells (HPCs) can undergo symmetric division, leading to the generation of two later healthy progenitor cells (see Lee and Vasioukhin \cite{lv}, and Zhang \textit{et al.} \cite{zdm}). Furthermore, later healthy progenitor cells (HPCs) can directly differentiate into other types of later healthy progenitor cells (see Carr \cite{c}):
\begin{equation}
  \label{R6}
  HPC \xrightarrow{\text{$k_{7}$}} 2HPC \tag{$R_{6}$}.
\end{equation}
\begin{equation}
  \label{R66}
  HPC \xrightarrow{\text{$\widetilde{k}_{7}$}} HPC \tag{$\widetilde{R}_{6}$}.
\end{equation}

Intermediate healthy progenitor cells (HPCs) can also undergo asymmetric division, giving rise to one later healthy progenitor cell (HPC) and one healthy differentiated cell (HDC) (see Lee and Vasioukhin \cite{lv}, and Gómez-López \textit{et al.} \cite{glp}):
\begin{equation}
  \label{R7}
  HPC \xrightarrow{\text{$k_{8}$}} HPC + HDC \tag{$R_{7}$}.
\end{equation}

Later healthy progenitor cells (HPCs) can directly differentiate into healthy differentiated cells (HDCs) and can also undergo symmetric differentiation, leading to the generation of two healthy differentiated cells (HDCs) (see Lee and Vasioukhin \cite{lv}, and Zhang \textit{et al.} \cite{zdm}):
\begin{equation}
  \label{R8}
  HPC \xrightarrow{\text{$k_{9}$}} HDC \tag{$R_{8}$},
\end{equation}
\begin{equation}
  \label{R9}
  HPC \xrightarrow{\text{$k_{10}$}} 2HDC \tag{$R_{9}$}.
\end{equation}

Healthy differentiated cells (HDCs) can directly differentiate into other types of healthy differentiated cells. They can also terminally differentiate and undergo symmetric differentiation into healthy terminally differentiated cells (HTDCs), which do not have the capacity to proliferate (see Duffy \textit{et al.} \cite{dwm}, Riffelmacher and Simon \cite{rs}, and Grinenko \textit{et al.} \cite{get}):
\begin{equation}
  \label{R99}
  HDC \xrightarrow{\text{$\widetilde{k}_{10}$}} HDC \tag{$\widetilde{R}_{9}$},
\end{equation}
\begin{equation}
  \label{R10}
  HDC \xrightarrow{\text{$k_{11}$}} HTDC \tag{$R_{10}$},
\end{equation}
\begin{equation}
  \label{R11}
  HDC \xrightarrow{\text{$k_{12}$}} 2HTDC \tag{$R_{11}$}.
\end{equation}

Each of the healthy cell types can undergo cell death due to differentiation, apoptosis and other elimination mechanisms (see Domen \cite{d}, Alenzi \textit{et al.} \cite{aaa}, and Cisneros \textit{et al.} \cite{cdc}):
\begin{equation}
  \label{R12}
  HCSC \xrightarrow{\text{$k_{13}$}} 0 \tag{$R_{12}$},
\end{equation}
\begin{equation}
  \label{R13}
  HPC \xrightarrow{\text{$k_{14}$}} 0 \tag{$R_{13}$},
\end{equation}
\begin{equation}
  \label{R14}
  HDC \xrightarrow{\text{$k_{15}$}} 0 \tag{$R_{14}$},
\end{equation}
\begin{equation}
  \label{R15}
  HTDC \xrightarrow{\text{$k_{16}$}} 0 \tag{$R_{15}$}.
\end{equation}

Similar pseudo-chemical reactions can be considered for the malignant cells (see Dingli \textit{et al.} \cite{dtm}, Lee and Vasioukhin \cite{lv}, Foo \textit{et al.} \cite{fdc}, Hu and Li \cite{hl}, and Vuelta \textit{et al.} \cite{vgh}), we have:
\begin{equation}
  \label{R16}
  MCSC \xrightarrow{\text{$k_{17}\varphi_{M}$}} 2MCSC \tag{$R_{16}$},
\end{equation}
where $\varphi_{M} = \frac{1}{1+B(HCSC+MCSC)}$.
\begin{equation}
  \label{R17}
   MCSC \underset{k_{19}}{\stackrel{k_{18}}{\rightleftarrows}} MQSC \tag{$R_{17}$},
\end{equation}
\begin{equation}
  \label{R18}
  MCSC \xrightarrow{\text{$k_{20}$}} MCSC + MPC \tag{$R_{18}$},
\end{equation}
\begin{equation}
  \label{R19}
  MCSC \xrightarrow{\text{$k_{21}$}} MPC \tag{$R_{19}$},
\end{equation}
\begin{equation}
  \label{R20}
  MCSC \xrightarrow{\text{$k_{22}$}} 2MPC \tag{$R_{20}$},
\end{equation}
\begin{equation}
  \label{R21}
  MPC \xrightarrow{\text{$k_{23}$}} 2MPC \tag{$R_{21}$},
\end{equation}
\begin{equation}
  \label{R211}
  MPC \xrightarrow{\text{$\widetilde{k}_{23}$}} MPC \tag{$\widetilde{R}_{21}$},
\end{equation}
\begin{equation}
  \label{R22}
  MPC \xrightarrow{\text{$k_{24}$}} MPC + MDC \tag{$R_{22}$},
\end{equation}
\begin{equation}
  \label{R23}
  MPC \xrightarrow{\text{$k_{25}$}} MDC \tag{$R_{23}$},
\end{equation}
\begin{equation}
  \label{R24}
  MPC \xrightarrow{\text{$k_{26}$}} 2MDC \tag{$R_{24}$},
\end{equation}
\begin{equation}
  \label{R244}
  MDC \xrightarrow{\text{$\widetilde{k}_{26}$}} MDC \tag{$\widetilde{R}_{24}$},
\end{equation}
\begin{equation}
  \label{R25}
  MDC \xrightarrow{\text{$k_{27}$}} MTDC \tag{$R_{25}$},
\end{equation}
\begin{equation}
  \label{R26}
  MDC \xrightarrow{\text{$k_{28}$}} 2MTDC \tag{$R_{26}$}.
\end{equation}
Similarly, as in the healthy case, each of the malignant cell types can undergo cell death (see Vivier \textit{et al.} \cite{vtb}, Alenzi \textit{et al.} \cite{aaa}, and Riether \textit{et al.} \cite{rso}):
\begin{equation}
  \label{R27}
  MCSC \xrightarrow{\text{$k_{29}$}} 0 \tag{$R_{27}$},
\end{equation}
\begin{equation}
  \label{R28}
  MPC \xrightarrow{\text{$k_{30}$}} 0 \tag{$R_{28}$},
\end{equation}
\begin{equation}
  \label{R29}
  MDC \xrightarrow{\text{$k_{31}$}} 0 \tag{$R_{29}$},
\end{equation}
\begin{equation}
  \label{R30}
  MTDC \xrightarrow{\text{$k_{32}$}} 0 \tag{$R_{30}$}.
\end{equation}

We observe that each pseudo-chemical reaction from (\ref{R2}) to (\ref{R15}) and (\ref{R17}) to (\ref{R30}), except the reactions (\ref{R1}) and (\ref{R16}), pushes the vector field in the direction of its reaction vector at a rate proportional to the product of its reactant concentrations. This method of rate construction is referred to as the \textit{mass-action} kinetics. Specifically, the rates of the pseudo-chemical reactions (\ref{R2}) - (\ref{R15}) and (\ref{R17}) - (\ref{R30}) are expressed as $k_{2}HCSC$,\ $k_{3}HQSC$,\ $k_{4}HCSC$,\ $k_{5}HCSC$,\ $k_{6}HCSC$,\ $k_{7}HPC$,\ $\widetilde{k}_{7}HPC$,\ $k_{8}HPC$,\ $k_{9}HPC$,\ $k_{10}HPC$,\ $\widetilde{k}_{10}HDC$,\ $k_{11}HDC$,\ $k_{12}HDC$,\ $k_{13}HCSC$,\ $k_{14}HPC$,\ $k_{15}HDC$,\ $k_{16}HTDC$,\ $k_{18}MCSC$,\ $k_{19}MQSC$,\ $k_{20}MCSC$,\ $k_{21}MCSC$,\ $k_{22}MCSC$,\ $k_{23}MPC$,\ $\widetilde{k}_{23}MPC$,\ $k_{24}MPC$,\ $k_{25}MPC$,\ $k_{26}MPC$,\ $\widetilde{k}_{26}MDC$, $k_{27}MDC$, $k_{28}MDC$,\ $k_{29}MCSC$,\ $k_{30}MPC$,\ $k_{31}MDC$ and $k_{32}MTDC$, where $k_{j_{H}} > 0$ for 
 $j_{H}=2,\ldots,16$\ , $\widetilde{k}_{j_{H}} > 0$ for $j_{H}=7,10$ and $k_{j_{M}} > 0$ for $j_{M}=18,\ldots,32$\ , $\widetilde{k}_{j_{M}} > 0$ for $j_{M}=23,26$. These quantities $k_{j_{H}/j_{M}}$ and $\widetilde{k}_{j_{H}/j_{M}}$ are referred to as the \textit{rate constants}, also known as \textit{kinetic parameters}. In pseudo-chemical reactions (\ref{R66})/(\ref{R211}) and (\ref{R99})/(\ref{R244}), intermediate healthy/malignant progenitor cells (HPCs/MPCs) are depicted as directly differentiating into later healthy/malignant progenitor cells, and healthy/malignant differentiated cells (HDCs/MDCs) are represented as directly differentiating into other types of healthy/malignant differentiated cells. However, it's important to note that these reactions are introduced for descriptive purposes and do not impact the dynamic behavior of our model. For instance, a reaction where $HPC$ converts to $HPC$ (or $HDC$ converts to $HDC$) does not alter the system's dynamics because it doesn't change the concentrations of any species. The pseudo-chemical reactions (\ref{R66})/(\ref{R211}) and (\ref{R99})/(\ref{R244}) are included to conceptually capture biological processes, but they do not affect the overall outcomes or changes in our model or simulations.

Therefore, due to pseudo-chemical reactions (\ref{R1}) and (\ref{R16}), where the rates are given by
$$
k_{1}\varphi_{H} HCSC = k_{1}\frac{1}{1+b_{1}HCSC+b_{2}MCSC}HCSC\ \ \ \ \text{and} \ \ \ \ k_{17}\varphi_{M} MCSC = k_{17}\frac{1}{1+B(HCSC+MCSC)}MCSC
$$
the entire pseudo-chemical reaction network (\ref{R1}) to (\ref{R30}) is classified as a \textit{nonmass-action} reaction network. In a \textit{nonmass-action} reaction network, the reaction rates are not directly proportional to the concentrations of the reactants but instead depend on more complex relationships, such as nonlinear dependencies or regulatory mechanisms. The functions $\varphi_{H} = \frac{1}{1+b_{1}HCSC+b_{2}MCSC}$ and $\varphi_{M} = \frac{1}{1+B(HCSC+MCSC)}$ are referred to as self-regulatory functions. These functions model the competition for space and resources (e.g. nutrients and oxygen) that occurs between stem cells within the stem cell niche (see Parajdi \textit{et al.} \cite{ppbt}, and Parajdi \cite{p}). Here, by competition, we do not refer to the same concept as in population dynamics, such as predator-prey models. Instead, we refer to a property of the mathematical model, specifically its competitive nature, which, in the case of acute/blast disease, leads to the proliferation of leukemic cells and the elimination of normal cells. In this context, we are not describing direct cell-to-cell competition but rather a competition between two cell populations: healthy and mutant (leukemic) cells, driven by shared resources and space within the stem cell niche. Experimental studies have shown that leukemic stem cells can outcompete healthy stem cells for access to the niche, often by remodeling the microenvironment or releasing factors that inhibit normal hematopoiesis. These findings provide biological support for modeling such interactions as competitive. For a classification of mathematical models into cooperative and competitive systems, we refer to the work of Smith \cite{sh}. 

We have chosen to apply the self-regulatory function exclusively to cycling stem cells, both healthy and malignant, as competition for space and resources occurs predominantly within this compartment. This assumption aligns with biological evidence indicating that the stem cell niche is the primary site where such competition takes place. The cycling stem cell compartment is also where cellular proliferation is most tightly regulated to maintain homeostasis and prevent unbounded growth. Notably, the dynamics induced by this self-regulatory mechanism at the level of cycling stem cells are consistently reflected in downstream cell lines (i.e. progenitor cells, differentiated cells and terminally differentiated cells). For instance, Dingli and Michor \cite{dm} present a model that applies the self-regulatory function exclusively to the stem cell compartment in the context of normal and leukemic stem cells, as well as differentiated cells. Similarly, both Michor \textit{et al.} \cite{mhi} and Foo \textit{et al.} \cite{fdc} introduce self-regulatory mechanisms applied exclusively to the stem cell compartment (or cycling stem cell compartment), an approach closely aligned with our model. Additionally, Parajdi \cite{p} and Parajdi \textit{et al.} \cite{ppbt} propose a model with eight equations describing the dynamics of four cell lines, demonstrating that the self-regulatory dynamics originating in the stem cell compartment are mirrored across other cell types. These studies strongly support the idea that applying the self-regulatory function exclusively to cycling stem cells does not compromise the biological validity of the model. Instead, this choice reflects the key role of the stem cell niche in regulating cellular dynamics. Furthermore, our assumption prevents unbounded growth in downstream compartments, as they depend on the self-regulatory function introduced at the stem cell compartment level. Therefore, we believe that our choice to apply the self-regulatory function solely to the cycling stem cell compartment is both mathematically convenient and biologically plausible, ensuring that competition and growth regulation are accurately represented in the model (see also Yang \textit{et al.} \cite{ya} for related discussions). Additionally, it's worth mentioning that each cellular event in the model is associated with a kinetic parameter defined by either $k_{j}$ or $\widetilde{k}_{j}$. For instance, here, $k_{1}$ and $k_{17}$ represent the intrinsic reaction rate constants, which indicate the natural tendency of HCSCs and MCSCs to divide symmetrically, producing two HCSCs (\ref{R1}) and two MCSCs (\ref{R16}), respectively. HCSCs and MCSCs can divide not only symmetrically but also asymmetrically. In the case of healthy stem cells, asymmetric division results in the generation of one HCSC and one HPC (\ref{R3}). Similarly, in the case of malignant stem cells, asymmetric division leads to one MCSC and one MPC (\ref{R18}). Symmetric and asymmetric divisions are not exclusive to healthy and malignant stem cells; they also occur with healthy and malignant progenitor cells (see Morrison and Kimble \cite{mk}, Dingli \textit{et al.} \cite{dtm}, Lee and Vasioukhin \cite{lv}, Zhang \textit{et al.} \cite{zdm}, and Gómez-López \textit{et al.} \cite{glp}). In addition to these divisions, healthy and malignant progenitor cells can directly differentiate. For example, the intrinsic reaction rate constants denoted by $\widetilde{k}_{7}$ and $\widetilde{k}_{23}$ represent the natural tendency of later HPCs and MPCs to directly differentiate into other types of later HPCs and MPCs, respectively. In the case of later healthy progenitor cells, this direct differentiation results in the generation of another type of later HPC but much more differentiated (\ref{R66}). Similarly, in the case of later malignant progenitor cells, direct differentiation yields another type of later MPC but much more differentiated (\ref{R211}). Direct differentiation is not limited to healthy and malignant progenitor cells, it can also take place in healthy and malignant differentiated cells (see Carr \cite{c}).

At a certain point in time, a mutation (Philadelphia chromosome (Ph)) occurs inside the healthy hematopoietic stem cell population (see Howard \textit{et al.} \cite{hhb}), leading to the gradual formation of two distinct populations: healthy stem cells and malignant or leukemic stem cells. This mutation, known as $t(9;22)$ (a reciprocal translocation of the ABL gene from chromosome $9$ to chromosome $22$, next to the BCR gene) (see Young \textit{et al.} \cite{y}), causes the entire cell population from the bone marrow to differentiate into two categories: healthy cells and malignant cells. In the following, to simplify the notation in the pseudo-chemical reaction network (\ref{R1})-(\ref{R30}), we will denote the healthy and malignant cells at time $t$ as $x_{i}(t)$ and $y_{i}(t)$, respectively, where $i$ ranges from $0$ to $4$. Please refer to the table below for a clear representation:

\begin{table}[H]
\centering
\begin{tabular}{|c|c|c|c|c|c|c|c|c|c|}
\hline
HCSCs & HQSCs & HPCs & HDCs & HTDCs & MCSCs & MQSCs & MPCs & MDCs & MTDCs \\ \hline
$x_{0}(t)$ & $x_{1}(t)$ & $x_{2}(t)$ & $x_{3}(t)$ & $x_{4}(t)$ & $y_{0}(t)$ & $y_{1}(t)$ & $y_{2}(t)$ & $y_{3}(t)$ & $y_{4}(t)$ \\ \hline
\end{tabular}
\end{table}

Using the aforementioned notations, we can express the complete pseudo-chemical reaction network (\ref{R1})-(\ref{R30}) as a nonlinear system of differential equations, which is:
\vspace{-0.2cm}

\begin{scriptsize}
\begin{equation*}
\begin{array}{llll}
\dot{x}_{0}(t)=\left(\frac{k_{1}}{1+b_{1}x_{0}+b_{2}y_{0}} - k_5 - k_6 -k_{13}\right)x_{0}-k_2x_0+k_{3}x_{1} &
\text{(HCSCs)} & \dot{y}_{0}(t)=\left(\frac{k_{17}}{1+B(x_{0}+y_{0})} -k_{21}-k_{22}-k_{29}\right)y_{0}-k_{18}y_{0}+k_{19}y_{1} & \text{(MCSCs)} \\[8pt]
\dot{x}_{1}(t)=k_{2}x_{0}-k_{3}x_{1} & \text{(HQSCs)} & \dot{y}_{1}(t)
=k_{18}y_{0}-k_{19}y_{1} & \text{(MQSCs)} \\[8pt]
\dot{x}_{2}(t)=
(k_{4}+k_{5}+2k_{6})x_{0}+k_7x_{2}-(k_{9}+k_{10}+k_{14})x_{2} & \text{(HPCs)} & \dot{y}_{2}(t)
=(k_{20}+k_{21}+2k_{22})y_{0}+k_{23}y_{2}-(k_{25}+k_{26}+k_{30})y_{2}  & \text{(MPCs)} \\[8pt]
\dot{x}_{3}(t)=(k_{8}+k_{9}+2k_{10})x_{2}-(k_{11}+k_{12}+k_{15})x_{3} & \text{(HDCs)} & \dot{y}_{3}(t)
=(k_{24}+k_{25}+2k_{26})y_{2}-(k_{27}+k_{28}+k_{31})y_{3} & \text{(MDCs)} \\[8pt]
\dot{x}_{4}(t)=(k_{11}+2k_{12})x_{3}-k_{16}x_{4} & \text{(HTDCs)} & \dot{y}_{4}(t)
=(k_{27}+2k_{28})y_{3}-k_{32}y_{4} & \text{(MTDCs).}
\end{array}
\end{equation*}
\end{scriptsize}

\noindent we will denote this nonlinear system by $(2)$. Here, the parameters $k_{j_{H}}$ where $j_{H}=1,\ldots,12$ and $k_{j_{M}}$ where $j_{M} = 17,\ldots,28$ represent the intrinsic reaction rate constants indicating the natural tendency of healthy and malignant cells to undergo direct differentiation, to divide either symmetrically or asymmetrically and to enter or leave a resting phase. For example, $k_{7}$ and $k_{23}$ are the natural tendency of intermediate HPCs and intermediate MPCs to divide symmetrically, producing two later HPCs and two later MPCs, respectively. The parameters $k_{j_{\overline{H}}}$ where $j_{\overline{H}}=13,\ldots,16$ and $k_{j_{\overline{M}}}$ where $j_{\overline{M}}=29,\ldots,32$ represent the intrinsic reaction rate constants that indicate the cell death rates (or cell turnover rates) for healthy and malignant cells. For example, $k_{14}$ and $k_{30}$ are the cell death rates for HPCs and MPCs, respectively; and $b_{1}, b_{2}, B$ are the sensitivity parameters that govern the self-limiting process at the level of HCSCs and MCSCs, respectively. For the meaning of the other intrinsic reaction rate constants in the system, see the description provided in Figure \ref{fig:mygraph}. Furthermore, we assume that for both healthy and malignant stem cells, the intrinsic reaction rate constants $k_{1}$ and $k_{17}$ (representing growth rates) are greater than the sum of intrinsic reaction rate constants $k_5+k_6+k_{13}$ and $k_{21}+k_{22}+k_{29}$ (indicative of division, differentiation and death rates). This can be expressed as:
\begin{equation*}
k_{1}>k_5+k_6+k_{13}\ \ \ \text{and}\ \ \ k_{17}>k_{21}+k_{22}+k_{29}.
\end{equation*}

The eventual advantage of malignant/leukemic stem cells of being less sensitive to the microenvironment than healthy stem cells is expressed by $b_{1} \geq b_{2} > B$. If we consider the extended system $(2)$ in the case when $b_{1}=b_{2}$, we will have only two non-zero steady states:
\vspace{-0.2cm}

\begin{footnotesize}
\begin{equation*}
\begin{aligned}
\bigg(&r, r \frac{k_2}{k_3}, r \frac{k_4+k_5+2k_6}{-k_7+k_9+k_{10}+k_{14}}, r \frac{(k_4+k_5+2k_6)(k_8+k_9+2k_{10})}{(-k_7+k_9+k_{10}+k_{14})(k_{11}+k_{12}+k_{15})}, r \frac{(k_4+k_5+2k_6)(k_8+k_9+2k_{10})(k_{11}+2k_{12})}{(-k_7+k_9+k_{10}+k_{14})(k_{11}+k_{12}+k_{15})k_{16}},0,0,0,0,0\bigg)\ \ \ \ \text{and} \\[8pt]
\bigg(&0,0,0,0,0,R, R \frac{k_{18}}{k_{19}}, R \frac{k_{20}+k_{21}+2k_{22}}{-k_{23}+k_{25}+k_{26}+k_{30}}, R \frac{(k_{20}+k_{21}+2k_{22})(k_{24}+k_{25}+2k_{26})}{(-k_{23}+k_{25}+k_{26}+k_{30})(k_{27}+k_{28}+k_{31})}, \\
&R \frac{(k_{20}+k_{21}+2k_{22})(k_{24}+k_{25}+2k_{26})(k_{27}+2k_{28})}{(-k_{23}+k_{25}+k_{26}+k_{30})(k_{27}+k_{28}+k_{31})k_{32}}\bigg)
\end{aligned}
\end{equation*}
\end{footnotesize}

\noindent where $r$ and $R$ represent the homeostatic amounts of healthy and malignant stem cells, and they are given by:
\begin{align}   \label{eq:star}
        r = \frac{k_{1}-k_{5}-k_{6}-k_{13}}{b_{1} (k_{5}+k_{6}+k_{13})} > 0\ \ \ \text{and}\ \ \ 
        R = \frac{k_{17}-k_{21}-k_{22}-k_{29}}{B (k_{21}+k_{22}+k_{29})} > 0.
        \tag{$\ast$}
\end{align}
In this paper, we consider that $b_{1}>b_{2}$. Therefore, under this assumption, in addition to the two non-zero steady states mentioned above, a new steady state $(\bar{x}_{0}^{\ast},\bar{x}_{1}^{\ast},\bar{x}_{2}^{\ast},\bar{x}_{3}^{\ast},\bar{x}_{4}^{\ast},\bar{y}_{0}^{\ast},\bar{y}_{1}^{\ast},\bar{y}_{2}^{\ast},\bar{y}_{3}^{\ast},\bar{y}_{4}^{\ast})$ could also exist, where all components are positive and non-zero. This makes the new model able to differentiate between chronic and accelerated-acute phases in chronic myeloid leukemia. Next, we continue to analyze the nonlinear system $(2)$, by determining the steady states and establishing both their local and global stability. \\

\textbf{(a.)} \textbf{Steady States}.  A steady state (or an equilibrium) is a constant solution, i.e. a solution for which $dx_i/dt = dy_i/dt =  0,$ where $i = 0,\ldots,4$. Hence, the steady states are obtained by solving an algebraic system derived from the nonlinear system of differential equations mentioned above, considering all the derivatives to be equal to zero. This algebraic system has four distinct steady states, denoted as follows:
\vspace{-0.2cm}

\begin{footnotesize}
\begin{equation*}
    \begin{aligned}
E_{0}\big(&0,0,0,0,0,0,0,0,0,0\big), \\[8pt]
E_{1}\bigg(&r, r \frac{k_2}{k_3}, r \frac{k_{4}+k_{5}+2k_{6}}{-k_7+k_9+k_{10}+k_{14}}, r \frac{(k_4 + k_5 + 2k_6)(k_8 + k_9 + 2k_{10})}{(-k_7 + k_9 + k_{10} + k_{14})(k_{11} + k_{12} + k_{15})},  r \frac{(k_{4} + k_{5} + 2k_{6})(k_{8} + k_{9} + 2k_{10})(k_{11} + 2k_{12})}{(-k_7 + k_9 + k_{10} + k_{14})(k_{11} + k_{12} + k_{15})k_{16}},0,0,0,0,0\bigg), \\[8pt]
E_{2}\bigg(&0,0,0,0,0,R, R \frac{k_{18}}{k_{19}},  R \frac{k_{20} + k_{21} + 2k_{22}}{-k_{23} + k_{25} + k_{26} + k_{30}}, R \frac{(k_{20} + k_{21} + 2k_{22})(k_{24} + k_{25} + 2k_{26})}{(-k_{23} + k_{25} + k_{26} + k_{30})(k_{27} + k_{28} + k_{31})}, \\
&R \frac{(k_{20} + k_{21} + 2k_{22})(k_{24} + k_{25} + 2k_{26})(k_{27} + 2k_{28})}{(-k_{23} + k_{25} + k_{26} + k_{30})(k_{27} + k_{28} + k_{31})k_{32}}\bigg), \\[8pt]
E_{3}\bigg(&\overline{x}^*, \overline{x}^*\frac{k_2}{k_3}, \overline{x}^* \frac{k_4+k_5+2k_6}{-k_7+k_9+k_{10}+k_{14}}, \overline{x}^* \frac{(k_4+k_5+2k_6)(k_8+k_9+2k_{10})}{(-k_7+k_9+k_{10}+k_{14})(k_{11}+k_{12}+k_{15})},  \overline{x}^* \frac{(k_4+k_5+2k_6)(k_8+k_9+2k_{10})(k_{11}+2k_{12})}{(-k_7+k_9+k_{10}+k_{14})(k_{11}+k_{12}+k_{15})k_{16}}, \\
&\overline{y}^*, \overline{y}^* \frac{k_{18}}{k_{19}},  \overline{y}^* \frac{k_{20}+k_{21}+2k_{22}}{-k_{23}+k_{25}+k_{26}+k_{30}},
\overline{y}^* \frac{(k_{20}+k_{21}+2k_{22})(k_{24}+k_{25}+2k_{26})}{(-k_{23}+k_{25}+k_{26}+k_{30})(k_{27}+k_{28}+k_{31})},   \overline{y}^* \frac{(k_{20}+k_{21}+2k_{22})(k_{24}+k_{25}+2k_{26})(k_{27}+2k_{28})}{(-k_{23}+k_{25}+k_{26}+k_{30})(k_{27}+k_{28}+k_{31})k_{32}}\bigg)
    \end{aligned}
\end{equation*}
\end{footnotesize}

\noindent    where 
    \begin{align*}
        r = \frac{k_{1}-k_{5}-k_{6}-k_{13}}{b_{1} (k_{5}+k_{6}+k_{13})}~,~
        R = \frac{k_{17}-k_{21}-k_{22}-k_{29}}{B (k_{21}+k_{22}+k_{29})}
    \end{align*}
    and
    \begin{align*}
        \overline{x}^{\ast} = \frac{(k_{5}+k_{6}+k_{13})(k_{22}+k_{29}-k_{17}+k_{21}) b_{2}+B (k_{1}-k_{5}-k_{6}-k_{13})(k_{21}+k_{22}+k_{29})}{(k_{5}+k_{6}+k_{13}) (k_{21}+k_{22}+k_{29}) (b_{1}-b_{2}) B},
\\[8pt]  \overline{y}^{\ast} = \frac{-(k_{5}+k_{6}+k_{13})(k_{22}+k_{29}-k_{17}+k_{21}) b_{1}-B (k_{1}-k_{5}-k_{6}-k_{13})(k_{21}+k_{22}+k_{29})}{(k_{5}+k_{6}+k_{13}) (k_{21}+k_{22}+k_{29}) (b_{1}-b_{2}) B}.
    \end{align*}
Direct calculation leads to:
\begin{align} \label{eq:starstar}
    \overline{x}^* = \frac{b_2}{b_1 - b_2}\left(\frac{b_1}{b_2}r-R\right)\ \ \text{and} \ \ \ \overline{y}^* = \frac{b_1}{b_1 - b_2}\left(R-r\right). \tag{$\ast\ast$}
\end{align}
Here, we assume that all steady states (or equilibrium points) are non-negative. To ensure this, we make the following assumptions about the intrinsic reaction rate constants:
\begin{align*}
    k_{9}+k_{10}+k_{14} > k_{7}\ \ \ \text{and}\ \ \  k_{25}+k_{26}+k_{30} > k_{23}.
\end{align*}
It is easy to observe, based on equation (\ref{eq:starstar}), that under the assumption $b_1 > b_2$, both $\overline{x}^*$ and $\overline{y}^*$ are positive if and only if $r<R<(b_1/b_2)r$.

In the following steps, we will prove the local stability for the nonlinear system $(2)$. To accomplish this, we simplify the nonlinear system of equations $(2)$, making it more computationally manageable while ensuring consistent mathematical outcomes. \\

\noindent Next, we proceed to rewrite the nonlinear system $(2)$ by setting:
\vspace{-0.2cm}

\begin{small}
\begin{equation*}
\begin{array}{llll}
a_0 = k_1, ~c_0 = k_5+k_6+k_{13}&\text{(HCSCs)} & A_0 = k_{17} , ~C_0 = k_{21}+k_{22}+k_{29} & \text{(MCSCs)}\\[8pt]
a_1=k_2,~c_1=k_3& \text{(HQSCs)}  & A_1
=k_{18},~C_1=k_{19} & \text{(MQSCs)} \\[8pt]
a_2 = k_4+k_5+2k_6, ~c_2=-k_7+k_9+k_{10}+k_{14}& \text{(HPCs)} & A_2 = k_{20}+k_{21}+2k_{22}, ~C_2 = -k_{23}+k_{25}+k_{26}+k_{30} & \text{(MPCs)}\\[8pt]
a_3=k_{8}+k_{9}+2k_{10}, ~c_3=k_{11}+k_{12}+k_{15}& \text{(HDCs)} & A_3 = k_{24}+k_{25}+2k_{26}, ~C_3 = k_{27}+k_{28}+k_{31} & \text{(MDCs)}\\[8pt]
a_4=k_{11}+2k_{12},~c_4=k_{16}& \text{(HTDCs)} & A_4 
=k_{27}+2k_{28}, ~C_4 = k_{32}& \text{(MTDCs).}
\end{array}%
\end{equation*}
\end{small}

\noindent Therefore, the new nonlinear system will be:
\begin{equation*}
\begin{array}{llll}
\dot{x}_{0}(t)=\left(\frac{a_{0}}{1+b_{1}x_{0}+b_{2}y_{0}}-a_{1}-c_{0} \right)x_{0}+c_{1}x_{1} &
\text{(HCSCs)} & \dot{y}_{0}(t)=\left(\frac{A_{0}}{1+B(x_{0}+y_{0})}-A_{1}-C_{0} \right)y_{0}+C_{1}y_{1} & \text{(MCSCs)} \\[8pt]
\dot{x}_{1}(t)=a_{1}x_{0}-c_{1}x_{1} & \text{(HQSCs)} & \dot{y}_{1}(t)
=A_{1}y_{0}-C_{1}y_{1} & \text{(MQSCs)} \\[8pt]
\dot{x}_{2}(t)=a_{2}x_{0}-c_{2}x_{2} & \text{(HPCs)} & \dot{y}_{2}(t)
=A_{2}y_{0}-C_{2}y_{2} & \text{(MPCs)} \\[8pt]
\dot{x}_{3}(t)=a_{3}x_{2}-c_{3}x_{3} & \text{(HDCs)} & \dot{y}_{3}(t)
=A_{3}y_{2}-C_{3}y_{3} & \text{(MDCs)} \\[8pt]
\dot{x}_{4}(t)=a_{4}x_{3}-c_{4}x_{4} & \text{(HTDCs)} & \dot{y}_{4}(t)
=A_{4}y_{3}-C_{4}y_{4} & \text{(MTDCs).}
\end{array}%
\vspace{0.1cm}
\end{equation*}
we will denote this simplified nonlinear system by $(3)$. 

Therefore, the steady states (or the equilibrium points) of the simplified nonlinear system $(3)$ can be rewritten as:
\begin{align*}
      \tilde{E_{0}} &(0,0,0,0,0,0,0,0,0,0), \\[8pt]
         \tilde{E_{1}} & \Big(d, d \frac{a_1}{c_1}, d \frac{a_2}{c_2}, d \frac{a_2a_3}{c_2c_3}, d \frac{a_2a_3a_4}{c_2c_3c_4} ,0,0,0,0,0 \Big), \\[8pt]
       \tilde{E_{2}} & \Big(0,0,0,0,0,D, D \frac{A_1}{C_1}, D \frac{A_2}{C_2}, D \frac{A_2A_3}{C_2C_3}, D \frac{A_2A_3A_4}{C_2C_3C_4} \Big), \\[8pt]
     \tilde{E_{3}} & \Big(x^*, x^*\frac{a_1}{c_1}, x^*\frac{a_2}{c_2}, x^*\frac{a_2a_3}{c_2c_3}, x^*\frac{a_2a_3a_4}{c_2c_3c_4}, y^*, y^*\frac{A_1}{C_1}, y^*\frac{A_2}{C_2}, y^*\frac{A_2A_3}{C_2C_3}, y^*\frac{A_2A_3A_4}{C_2C_3C_4} \Big),
    \end{align*}
    where 
    \begin{align*}
        d = \frac{a_0 - c_0}{b_1c_0},\ \ D = \frac{A_0-C_0}{BC_0},
    \end{align*}
    and
    \begin{align*}
        x^* = \frac{BC_0a_0-BC_0c_0-A_0b_2c_0+C_0b_2c_0}{BC_0(b_1-b_2)c_0},\ \  y^* = -\frac{BC_0a_0-BC_0c_0-A_0b_1c_0+C_0b_1c_0}{BC_0c_0(b_1-b_2)} .
    \end{align*}
Direct calculation leads to:
\begin{align*}
    x^* = \frac{b_2}{b_1 - b_2}\left(\frac{b_1}{b_2}d-D\right),\ \  y^* = \frac{b_1}{b_1 - b_2}\left(D-d\right).
\end{align*}

It is easy to observe that under the assumption that $b_1 > b_2$, both $x^*$ and $y^*$ are positive if and only if:
\begin{align*}
    d<D<\frac{b_1}{b_2}d.
\end{align*}
Thus, in addition to the non-zero steady states $\tilde{E_1}$ and $\tilde{E_2}$, a positive steady state $\tilde{E_3}$ appears, contrary to the case where $b_1 = b_2$. After changing the variables, both inequalities $k_1>k_5+k_6+k_{13}$ and $k_{17}>k_{21}+k_{22}+k_{29}$ can be expressed as:
$$
a_{0} > c_{0}\ \ \  \text{and}\ \ \ A_{0} > C_{0}.
$$
It's important to highlight that $E_{0}$ is equivalent to $\tilde{E_{0}}$, $E_{1}$ is equivalent to $\tilde{E_{1}}$, $E_{2}$ is equivalent to $\tilde{E_{2}}$ and finally $E_{3}$ is equivalent to $\tilde{E_{3}}$. \\

\textbf{(b).} \textbf{Local Stability.} The local stability of the simplified model, system $(3)$, is mathematically equivalent to that of the original model, system $(2)$. Thus, we study the stability of the steady states of the simplified system $(2)$ using the standard first approximation Lyapunov's method (for details, see Kaplan and Glass \cite{kg}, Coddington and Levinson \cite{cl}, and Jones \textit{et al.} \cite{jps}). According to this method, an equilibrium is considered \textit{asymptotically stable} if the Jacobian matrix is a Hurwitz matrix, which means that Re $\lambda < 0$ for all its characteristic roots $\lambda$, and is \textit{unstable} if Re $\lambda > 0$ for at least one of its characteristic roots. \\

$\rhd$ For the steady state $\tilde{E_0}$, the characteristic equation of the Jacobian matrix $J(\tilde{E_0})$, evaluated at $\tilde{E_0}$ for the simplified system $(3)$, is as follows:
\begin{equation*}
\begin{aligned}
    &(\lambda+c_2)(\lambda+c_3)(\lambda+c_4)(\lambda+C_2)(\lambda+C_3)(\lambda+C_4)
    (\lambda^2 + (-A_0+C_0+A_1+C_1)\lambda - C_1(A_0 -C_0))\\
    &(\lambda^2 + (a_1 + c_1 - a_0 + c_0)\lambda - c_1(a_0-c_0)) = 0.
\end{aligned}
\end{equation*}
Using this characteristic equation, we can compute the eigenvalues: $\lambda_1 = -c_2 < 0$, $\lambda_2 = -c_3 < 0$, $\lambda_3 = -c_4 < 0$, $\lambda_4 = -C_2 < 0$, $\lambda_5 = -C_3 < 0$ and $\lambda_6 = -C_4 < 0$. We can observe that there are two quadratic polynomials in this characteristic equation. The product of the roots of the first quadratic polynomial, $\lambda^2 + (-A_0+C_0+A_1+C_1)\lambda - C_1(A_0 -C_0)$, can be expressed as follows: 
\begin{align*}
    \lambda_7\lambda_8 = -C_1(A_0 - C_0) < 0,
\end{align*}
based on the assumption that $A_0 > C_0$. A similar conclusion arises when considering the second quadratic polynomial, $\lambda^2 + (a_1 + c_1 - a_0 + c_0)\lambda - c_1(a_0-c_0)$, i.e. $\lambda_9\lambda_{10} = -c_1(a_0 - c_0) < 0$, based on the assumption that $a_0 > c_0$. Consequently, the signs of $\lambda_7$ and $\lambda_8$ (similarly for $\lambda_{9}$ and $\lambda_{10}$) are opposite, leading to two real eigenvalues in this characteristic equation. Thus, the steady state $\tilde{E_0}$ is unstable. 

It's important to note that the same conclusion holds, i.e. the same inequalities hold, even when we work with the initial mathematical model $(2)$. Therefore, the eigenvalues that correspond to the characteristic equation of the Jacobian matrix $J(E_0)$, evaluated at $E_0$ for system $(2)$, are:
\begin{equation*}
\begin{aligned}
    \lambda_1 &= -c_2 = k_7-k_9-k_{10}-k_{14} < 0\ \ \ \text{based on the assumption that $k_9+k_{10}+k_{14} > k_{7}$} \\  \lambda_2 &= -c_3 = -k_{11}-k_{12}-k_{15} < 0,\ \ \lambda_3 = -c_4 = -k_{16} < 0, \\ 
    \lambda_4 &= -C_2 = k_{23}-k_{25}-k_{26}-k_{30} < 0\ \ \ \text{based on the assumption that $k_{25}+k_{26}+k_{30} > k_{23}$} \\ 
    \lambda_5 &= -C_3 = -k_{27}-k_{28}-k_{31} < 0,\ \    \lambda_6 = -C_4 = -k_{32} < 0, \\
    \lambda_7\lambda_8 &= -C_1(A_0 - C_0) = -k_{19}(k_{17}-k_{21}-k_{22}-k_{29}) < 0\ \ \ 
    \text{based on the assumption that $k_{17} > k_{21}+k_{22}+k_{29}$}\ \ \text{and} \\
    \lambda_9\lambda_{10} &= -c_1(a_0 - c_0) = -k_{3}(k_{1}-k_{5}-k_{6}-k_{13}) < 0 \ \ \ 
    \text{based on the assumption that $k_{1} > k_{5}+k_{6}+k_{13}$}.
\end{aligned}
\end{equation*}
Thus, we can conclude that the steady state $E_{0}$ is also unstable. \\

$\rhd$ For the steady state $\tilde{E_1}$, the characteristic equation of the Jacobian matrix $J(\tilde{E_1})$, evaluated at $\tilde{E_1}$ for the simplified system $(3)$, is as follows:
\begin{small}
\begin{equation*}
\begin{aligned}
    &\frac{1}{(d b_1 + 1)^2 (Bd+1)}(\lambda+c_2)(\lambda+c_3)(\lambda+c_4)(\lambda+C_2)(\lambda+C_3)(\lambda+C_4)\\
    &\Bigg((Bd+1)\lambda^2 + \Big(B(C_0+A_1+C_1)d-A_0+C_0+A_1+C_1\Big)\lambda + C_1(BdC_0-A_0+C_0) \Bigg)\\
    &\Bigg((d b_1 + 1)^2\lambda^2 + \Big(b_1^2(a_1+c_1+c_0)d^2 + 2b_1(a_1+c_1+c_0)d+a_1+c_1-a_0+c_0\Big)\lambda + c_1(d^2b_1^2c_0+2db_1c_0-a_0+c_0) \Bigg)
     = 0.
\end{aligned}
\end{equation*}
\end{small}

\noindent Using this characteristic equation, we can compute the eigenvalues $\lambda_1 = -c_2 < 0$, $\lambda_2 = -c_3 < 0$, $\lambda_3 = -c_4 < 0$, $\lambda_4 = -C_2 < 0$, $\lambda_5 = -C_3 < 0$ and $\lambda_6 = -C_4 < 0$. We can observe that there are two quadratic polynomials in this characteristic equation. To prove that all the roots (eigenvalues) of this characteristic equation are negative, we can show that the coefficients of both quadratic polynomials are all positive. Here, we use the following statement: the roots of a quadratic equation will have negative real parts if and only if all the coefficients have the same sign. Let's start with the first quadratic polynomial, which has the following coefficients:
\begin{align*}
    \eta_{11}&=Bd+1 > 0, \\
    \eta_{12}&=B(C_0+A_1+C_1)d-A_0+C_0+A_1+C_1 = BC_0(d-D)+(A_1+C_1)(Bd+1) > 0 \ \ \text{if}\ \ D<d, \\
    \eta_{13}&= C_1(BdC_0-A_0+C_0) = C_1BC_0(d-D) > 0\ \ \text{if}\ \ D<d.
\end{align*}
Consider the second quadratic polynomial, which has the following coefficients:
\begin{align*}
    \eta_{21}&=(db_1+1)^2 > 0,\\ \eta_{22}&=b_1^2(a_1+c_1+c_0)d^2 + 2b_1(a_1+c_1+c_0)d+a_1+c_1-a_0+c_0 \\
    &= b_1^2(a_1+c_1+c_0)d^2 +b_1c_0d + (a_1+c_1)(2b_1d + 1) > 0, \\
    \eta_{23}&=c_1(d^2b_1^2c_0+2db_1c_0-a_0+c_0) = c_1(d^2b_1^2c_0 + db_1c_0) > 0.
\end{align*}
Since the coefficients for the first quadratic equation $\eta_{11}, \eta_{12}$ and $\eta_{13}$ are all positive if we assume $D<d$, then by the statement we mentioned previously, the roots $\lambda_7$ and $\lambda_8$ of this quadratic equation have negative real parts if and only if $D < d$, similarly, $\lambda_9$ and $\lambda_{10}$ have negative real parts because all the coefficients $\eta_{21},\ \eta_{22}$ and $\eta_{23}$ of the second quadratic polynomial (or equation) are positive. Thus, the steady state $\tilde{E_1}$ is asymptotically stable if and only if $D<d$. On the contrary, if $D>d$, then the equilibrium $\tilde{E_1}$ is unstable.

Here, it's important to note that the same conclusion holds, i.e. the same inequalities hold, even when we work with the initial mathematical model $(2)$:
\begin{align*}
    \eta_{11}&=Bd+1 = Br+1 > 0, \\
    \eta_{12}&=BC_{0}(d-D)+(A_1+C_1)(Bd+1) = B(k_{21}+k_{22}+k_{29})(r-R)+(k_{18}+k_{19})(Br+1) > 0\ \ \text{if}\ \ R<r, \\
    \eta_{13}&= C_1BC_0(d-D) = k_{19}B(k_{21}+k_{22}+k_{29})(r-R) > 0\ \ \text{if}\ \ R<r, \\
    \eta_{21}&= (db_1+1)^2 = (rb_1+1)^2 > 0, \\
    \eta_{22}&= b_1^2(a_1+c_1+c_0)d^2 +b_1c_0d + (a_1+c_1)(2b_1d + 1) \\
    &= b_1^2(k_2+k_3+k_5+k_6+k_{13})r^2 + b_1(k_5+k_6+k_{13})r + (k_2 + k_3)(2b_1r+1) > 0, \\
    \eta_{23}&= c_1(d^2b_1^2c_0 + db_1c_0) = k_3r^2b_1^2(k_5+k_6+k_{13}) + k_3rb_1(k_5+k_6+k_{13}) > 0.
\end{align*}
Therefore, the eigenvalues that correspond to the characteristic equation of the Jacobian matrix $J(E_1)$, evaluated at $E_1$ for system $(2)$, are:
\begin{equation*}
\begin{aligned}
    \lambda_1 &= -c_2 = k_7-k_9-k_{10}-k_{14} < 0\ \ \ \text{based on the assumption that $k_9+k_{10}+k_{14} > k_{7}$} \\  \lambda_2 &= -c_3 = -k_{11}-k_{12}-k_{15} < 0,\ \ \lambda_3 = -c_4 = -k_{16} < 0, \\ 
    \lambda_4 &= -C_2 = k_{23}-k_{25}-k_{26}-k_{30} < 0\ \ \ \text{based on the assumption that $k_{25}+k_{26}+k_{30} > k_{23}$} \\ 
    \lambda_5 &= -C_3 = -k_{27}-k_{28}-k_{31} < 0,\ \    \lambda_6 = -C_4 = -k_{32} < 0, \\
    \lambda_7 & \ \text{and}\ \lambda_8\ \text{have negative real parts if and only if $R < r$, while} \\
    \lambda_9 & \ \text{and}\ \lambda_{10}\ \text{have negative real parts because all the coefficients of the second quadratic equation are positive}.
\end{aligned}
\end{equation*}
Thus, we can conclude that the steady state $E_{1}$ is also asymptotically stable if and only if $R<r$. On the contrary, if $R>r$, then the equilibrium $E_{1}$ is unstable.\\
 
$\rhd$ Similarly, for the steady state $\tilde{E_2}$, the characteristic equation of the Jacobian matrix $J(\tilde{E_2})$, evaluated at $\tilde{E_2}$ for the simplified system $(3)$, is as follows:
\begin{small}
\begin{equation*}
\begin{aligned}
    &\frac{1}{(D b_2 + 1) (BD+1)^2}(\lambda+c_2)(\lambda+c_3)(\lambda+c_4)(\lambda+C_2)(\lambda+C_3)(\lambda+C_4) \\ &\Bigg((b_2D+1)\lambda^2 + \Big(b_2(c_0+a_1+c_1)D -a_0 +c_0 + a_1+c_1\Big)\lambda + c_1(b_2Dc_0-a_0+c_0) \Bigg) \\
    &\Bigg((DB + 1)^2\lambda^2 + \Big(B^2(A_1+C_1+C_0)D^2 +
    2B(A_1+C_1+C_0)D + A_1+C_1-A_0+C_0\Big)\lambda + C_1(D^2B^2C_0+2DBC_0-A_0+C_0) \Bigg)
     = 0.
\end{aligned}
\end{equation*}
\end{small}

\noindent Using this characteristic equation, we can compute the eigenvalues $\lambda_1 = -c_2 < 0$, $\lambda_2 = -c_3 < 0$, $\lambda_3 = -c_4 < 0$, $\lambda_4 = -C_2 < 0$, $\lambda_5 = -C_3 < 0$ and $\lambda_6 = -C_4 < 0$. Similarly to the case of the steady state $\tilde{E_1}$, we can observe that there are two quadratic polynomials in this characteristic equation. To prove that all the roots (eigenvalues) of this characteristic equation have negative real parts, we can show that the coefficients of both quadratic polynomials are all positive. Let's start with the first quadratic polynomial, which has the following coefficients:
\begin{align*}
    \eta_{11}&= b_2D+1 > 0, \\
    \eta_{12}&= b_2(c_0+a_1+c_1)D -a_0 +c_0 +a_1 +c_1 = c_0b_2\left(D-\frac{b_1}{b_2}d\right) + (a_1+c_1)(b_2D+1) > 0 \ \ \text{if}\ \ D>\frac{b_1}{b_2}d, \\
    \eta_{13}&= c_1(b_2Dc_0 - a_0 + c_0) = c_1c_0b_2\left(D-\frac{b_1}{b_2}d\right) > 0\ \ \text{if}\ \ D>\frac{b_1}{b_2}d.
\end{align*}
Consider the second quadratic polynomial, which has the following coefficients:
\begin{align*}
    \eta_{21}&=(DB+1)^2 > 0,\\ \eta_{22}&=B^2(A_1+C_1+C_0)D^2+
    2B(A_1+C_1+C_0)D + A_1+C_1-A_0+C_0 \\
    &= B^2(A_1+C_1+C_0)D^2 + DBC_0 + (A_1+C_1)(2BD+1) > 0, \\
    \eta_{23}&= C_1(D^2B^2C_0+2DBC_0-A_0+C_0) = C_1(D^2B^2C_0 + DBC_0) > 0.
\end{align*}
Therefore, the eigenvalues $\lambda_7$ and $\lambda_8$ have negative real parts if and only if $D > (b_1/b_2)d$, while $\lambda_9$ and $\lambda_{10}$ have negative real parts because all the coefficients $\eta_{21},\ \eta_{22}$ and $\eta_{23}$ of the second quadratic polynomial (or equation) are positive. Thus, the steady state $\tilde{E_2}$ is asymptotically stable if and only if $D>(b_1/b_2)d$. On the contrary, if $D<(b_1/b_2)d$, then the equilibrium $\tilde{E_2}$ is unstable.

Here, it's important to note that the same conclusion holds, i.e. the same inequalities hold, even when we work with the initial mathematical model $(2)$:
\begin{align*}
    \eta_{11}&= b_2D+1 = b_2R+1 > 0, \\
    \eta_{12}&= c_0b_2\left(D-\frac{b_1}{b_2}d\right) + (a_1+c_1)(b_2D+1) = (k_5+k_6+k_{13})b_2\left(R-\frac{b_1}{b_2}r\right) + (k_2+k_3)(b_2R+1) > 0 \ \ \text{if}\ \ R>\frac{b_1}{b_2}r, \\
    \eta_{13}&= c_1c_0b_2\left(D-\frac{b_1}{b_2}d\right) = k_3(k_5+k_6+k_{13})b_2\left(R-\frac{b_1}{b_2}r\right) > 0\ \ \text{if}\ \ R>\frac{b_1}{b_2}r, \\
    \eta_{21}&=(DB+1)^2 = (RB+1)^2 > 0,\\ \eta_{22}&= B^2(A_1+C_1+C_0)D^2 + DBC_0 + (A_1+C_1)(2BD+1) \\
    &= B^2(k_{18}+k_{19}+k_{21}+k_{22}+k_{29})R^2 + RB(k_{21}+k_{22}+k_{29}) + (k_{18}+k_{19})(2BR+1) > 0, \\
    \eta_{23}&= C_1(D^2B^2C_0 + DBC_0) = k_{19}(R^2B^2(k_{21}+k_{22}+k_{29}) + RB(k_{21}+k_{22}+k_{29})) > 0.
\end{align*}
Therefore, the eigenvalues that correspond to the characteristic equation of the Jacobian matrix $J(E_2)$, evaluated at $E_2$ for system $(2)$, are:
\begin{equation*}
\begin{aligned}
    \lambda_1 &= -c_2 = k_7-k_9-k_{10}-k_{14} < 0\ \ \ \text{based on the assumption that $k_9+k_{10}+k_{14} > k_{7}$} \\  \lambda_2 &= -c_3 = -k_{11}-k_{12}-k_{15} < 0,\ \ \lambda_3 = -c_4 = -k_{16} < 0, \\ 
    \lambda_4 &= -C_2 = k_{23}-k_{25}-k_{26}-k_{30} < 0\ \ \ \text{based on the assumption that $k_{25}+k_{26}+k_{30} > k_{23}$} \\ 
    \lambda_5 &= -C_3 = -k_{27}-k_{28}-k_{31} < 0,\ \    \lambda_6 = -C_4 = -k_{32} < 0, \\
    \lambda_7 & \ \text{and}\ \lambda_8\ \text{have negative real parts if and only if $R > \frac{b_1}{b_2}r$ while} \\
    \lambda_9 & \ \text{and}\ \lambda_{10}\ \text{have negative real parts because all the coefficients of the second quadratic equation are positive}.
\end{aligned}
\end{equation*}
Thus, we can conclude that the steady state $E_{2}$ is also asymptotically stable if and only if $R>(b_1/b_2)r$. On the contrary, if $R<(b_1/b_2)r$, then the equilibrium $E_{2}$ is unstable. \\

$\rhd$ When it comes to the analysis of the steady state $\tilde{E_3}$, the characteristic equation of the Jacobian matrix $J(\tilde{E_3})$, evaluated at $\tilde{E_3}$ for the simplified system $(3)$, looks a bit complicated and is provided in the \textit{Supplementary Notes} and \textit{Maple Supplementary file}: $\textbf{SupplMaple.mw}$. We still want to apply the same method to show that all the eigenvalues of the characteristic equation have negative real parts. Here, treating the characteristic equation as a polynomial of $\lambda$ with degree $10$ and denoting each eigenvalue from this equation as $\lambda_i$, where $i$ ranges from $1$ to $10$, the characteristic equation is:
\begin{align*}
  \frac{1}{(BD+1)^2(db_{1}+1)^2(b_{1}-b_{2})}\Big((\lambda+c_2)(\lambda+c_3)(\lambda+c_4)(\lambda+C_2)(\lambda+C_3)(\lambda+C_4)P(\lambda)\Big) = 0
\end{align*}
where $P(\lambda)$ is a degree $4$ polynomial (check the \textit{Maple Supplementary file} for full formulation of $P(\lambda)$):
\begin{align*}
   P(\lambda) = \mu_0 \lambda^4 + \mu_1 \lambda^3 + \mu_2 \lambda^2 + \mu_3 \lambda + \mu_4.
\end{align*}
Using this characteristic equation, we can compute the first six eigenvalues $\lambda_1 = -c_2 < 0$, $\lambda_2 = -c_3 < 0$, $\lambda_3 = -c_4 < 0$, $\lambda_4 = -C_2 < 0$, $\lambda_5 = -C_3 < 0$ and $\lambda_6 = -C_4 < 0$. Next, we focus on the polynomial $P(\lambda)$. First, we note:
\begin{align*}
    \mu_0 &= (BD+1)^2(b_1-b_2)(db_1+1)^2 > 0\ \ \text{if and only if}\ \ b_1>b_2\ .
\end{align*}
Based on the Routh-Hurwitz criterion for the 4th-order polynomials, if we can prove that $\mu_1, ..., \mu_4$ are all positive and that $\mu_1\mu_2-\mu_0\mu_3 > 0$ and $\mu_1\mu_2\mu_3 - \mu_1^2\mu_4 - \mu_0\mu_3^2 > 0$, then we can conclude that the eigenvalues $\lambda_7$, $\lambda_8$, $\lambda_9$ and $\lambda_{10}$ have negative real parts, indicating that $\tilde{E_3}$ is asymptotically stable. To begin, we rewrite:
\begin{align*}
    a_0 := db_1c_0+c_0\ \ \text{and}\ \  A_0 := BC_0D+C_0
\end{align*}
and plug them into $\mu_i$, where $i=1,\ldots,4$ to reduce the number of variables. 

Let's start by checking the sign of $\mu_1$.
\vspace{-0.3cm}

\begin{small}
\begin{align*}
  \mu_1 =&\ (d b_{1}+1) (B D +1) \Big[\ (-Db_2+db_1)\Big(BDb_1c_0 +b_1c_0\Big) + (D-d)\Big(BdC_0b_1^2 + BC_0b_1\Big) + (b_1-b_2)\Big(BDdA_1b_1 + BDdC_1b_1 \\
  & + BDda_1b_1 + BDdb_1c_1 + BDA_1 + BDC_1 + BDa_1 + BDc_1 + dA_1b_1 + dC_1b_1 + da_1b_1 + db_1c_1 + A_1 + C_1 + a_1 + c_1\Big)\Big]\ \\ 
  =&\ (d b_{1}+1) (B D +1) \Big[\ (-Db_2+db_1)\nu_{11} + (D-d)\nu_{12} + (b_1-b_2)\nu_{13}\Big]\ . &&
\end{align*}
\end{small}

\noindent Thus, $\ \mu_1 > 0\ $ if and only if $\ b_1 > b_2$ and $d<D<(b_1/b_2)d$.\\

Similarly, we can show that $\mu_2,\ \mu_3$ and $\mu_4$ are all positive under the same assumptions. For full details on how to obtain the factor form of $\mu_{1},\ \mu_{2}$ and $\mu_{3}$, please refer to the \textit{Supplementary Notes}. These coefficients $\mu_{1}$, $\mu_{2}$ and $\mu_{3}$ were obtained using the Maple \textit{software}. Next, we have:
\vspace{-0.2cm}

\begin{small}
\begin{flalign*}
\mu_2 =&\ (d b_{1}+1) (B D +1)\Big[\ (-Db_2+db_1)\Big(BDA_1b_1c_0 + BDC_1b_1c_0 + BDb_1c_0c_1 + A_1b_1c_0 + C_1b_1c_0 + b_1c_0c_1\Big)  \\&+  (D-d)\Big(BC_0b_1c_0(-Db_2+db_1) + BdC_0C_1b_1^2 + BdC_0a_1b_1^2
+ BdC_0b_1^2c_1 + BC_0C_1b_1 + BC_0a_1b_1 + BC_0b_1c_1\Big) 
  \\& + (b_1-b_2)\Big(BDdA_1a_1b_1 + BDdA_1b_1c_1 + BDdC_1a_1b_1 + BDdC_1b_1c_1 +BDA_1a_1 + BDA_1c_1 + BDC_1a_1 + BDC_1c_1 \\ & +  dA_1a_1b_1 + dA_1b_1c_1 + dC_1a_1b_1 + 
dC_1b_1c_1 + A_1a_1 + A_1c_1 + C_1a_1 + C_1c_1\Big)\Big]\ \\
=&\ (d b_{1}+1) (B D +1)\Big[\ (-Db_2+db_1)\nu_{21} + (D-d)\nu_{22} + (b_1-b_2)\nu_{23}\Big]\ , &&
\end{flalign*}
\end{small}
\vspace{-0.4cm}
\begin{small}
\begin{flalign*}
\mu_3 =&\ b_{1} (d b_{1}+1) (B D +1)\Big[\ (-Db_2 + db_1)\Big(BDA_1c_0c_1 + (D-d)(BC_0C_1c_0 + BC_0c_0c_1)  + BDC_1c_0c_1 + A_1c_0c_1 + C_1c_0c_1\Big) \\
 & + (D-d)\Big( BdC_0C_1a_1b_1 + BdC_0C_1b_1c_1 + BC_0C_1a_1 + BC_0C_1c_1\Big)\Big]\ \\
 =&\ b_{1} (d b_{1}+1) (B D +1)\Big[\ (-Db_2 + db_1)\nu_{31} + (D-d)\nu_{32}\Big]\ , &&
\end{flalign*}
\end{small}
\vspace{-0.4cm}
\begin{small}
\begin{flalign*}
&\mu_4 = b_{1} B C_{0} c_{0} (d b_{1}+1) (D -d) (-D b_{2}+d b_{1}) (B D+1) C_{1} c_{1}\ .&&
\end{flalign*}
\end{small}

\noindent Thus, $\ \mu_2,\ \mu_3,\ \mu_4 > 0 $ if and only if $\ b_1 > b_2$ and $d<D<(b_1/b_2)d$.
\vspace{0.3cm}
 
In addition, to show the steady state $\tilde{E_3}$ is asymptotically stable, we need to verify two more conditions. The first condition:
\begin{small}
\begin{align*}
    \mu_1\mu_2 - \mu_0\mu_3 > 0
\end{align*}
\end{small}

\noindent if and only if $\ b_{1}>b_{2}\ $ and $\ d<D<(b_1/b_2)d$. The second condition:
\begin{small}
\begin{align*}
\mu_1\mu_2\mu_3 - \mu_1^2\mu_4 - \mu_0\mu_3^2 > 0 
\end{align*}
\end{small}

\noindent if and only if $\ b_{1} > b_{2}\ $ and $\ d<D<(b_{1}/b_{2})d$. 

Thus, the steady state $\tilde{E_3}$ is asymptotically stable if and only if $b_1 > b_2\ $ and $\ d<D<(b_1/b_2)d$. The proof of the two additional conditions, $\mu_{1}\mu_{2} - \mu_{0}\mu_{3}$ and  $\mu_1\mu_2\mu_3 - \mu_1^2\mu_4 - \mu_0\mu_3^2$, being positive, can be found in the \textit{Supplementary Notes}.
\vspace{0.3cm}

$\blacktriangleright$ It's important to note that the same conclusion holds (based on the Routh-Hurwitz criterion for the 4th-order polynomials), i.e. the same inequalities hold, even when we work with the initial mathematical model $(2)$. We have:
\vspace{-0.2cm}

\begin{small}
\begin{flalign*}
\mu_0 =&\ (BR+1)^2(b_1-b_2)(rb_1+1)^2 > 0\ \ \text{if and only if}\ \ b_1>b_2, && 
  \end{flalign*}
\end{small}
\vspace{-0.4cm}
\begin{small}
\begin{flalign*}
\mu_1 =&\ (r b_{1}+1) (B R +1) \Big[\ (-Rb_2+rb_1)\Big(BRb_1(k_5+k_6+k_{13}) +b_1(k_5+k_6+k_{13})\Big) + (R-r)\Big(Br(k_{21}+k_{22}+k_{29})b_1^2 \\
&+ B(k_{21}+k_{22}+k_{29})b_1\Big) + (b_1-b_2)\Big(BRrk_{18}b_1 + BRrk_{19}b_1 + BRrk_{2}b_1 + BRrb_1k_{3} + BRk_{18} + BRk_{19} + BRk_{2} \\
&+ BRk_{3} + rk_{18}b_1 + rk_{19}b_1 + rk_{2}b_1 + rb_1k_{3} + k_{18} + k_{19} + k_{2} + k_{3}\Big)\Big]\ \\ 
  =&\ (r b_{1}+1) (B R +1) \Big[\ (-Rb_2+rb_1)\overline{\nu}_{11} + (R-r)\overline{\nu}_{12} + (b_1-b_2)\overline{\nu}_{13}\Big]\ > 0\ \ \text{if and only if}\ \ b_1>b_2\ \ \text{and}\ \ r<R<\frac{b_1}{b_2}r, &&
  \end{flalign*}
\end{small}
\vspace{-0.4cm}
\begin{small}
\begin{flalign*}
\mu_2 =&\ (r b_{1}+1) (B R +1)\Big[\ (-Rb_2+rb_1)\Big(BRk_{18}b_1(k_5+k_6+k_{13}) + BRk_{19}b_1(k_5+k_6+k_{13}) + BRb_1(k_5+k_6+k_{13})k_{3} \\
&+ k_{18}b_1(k_5+k_6+k_{13}) + k_{19}b_1(k_{5}+k_{6}+k_{13}) + b_1(k_5+k_6+k_{13})k_{3}\Big) +  (R-r)\Big(B(k_{21}+k_{22}+k_{29})b_1(k_{5}+k_{6}+k_{13})\\
&\ (-Rb_2+rb_1) + Br(k_{21}+k_{22}+k_{29})k_{19}b_1^2 + Br(k_{21}+k_{22}+k_{29})k_{2}b_1^2
+ Br(k_{21}+k_{22}+k_{29})b_1^2k_{3} + B(k_{21}+k_{22}+k_{29})k_{19}b_1 \\
&+ B(k_{21}+k_{22}+k_{29})k_{2}b_1 + B(k_{21}+k_{22}+k_{29})b_1k_{3}\Big) + (b_1-b_2)\Big(BRrk_{18}k_{2}b_1 + BRrk_{18}b_1k_{3} + BRrk_{19}k_{2}b_1 + BRrk_{19}b_1k_{3} \\
&+ BRk_{18}k_{2} + BRk_{18}k_{3} + BRk_{19}k_{2} + BRk_{19}k_{3} +  rk_{18}k_{2}b_1 + rk_{18}b_1k_{3} + rk_{19}k_{2}b_1 + 
rk_{19}b_1k_{3} + k_{18}k_{2} + k_{18}k_{3} + k_{19}k_{2} + k_{19}k_{3}\Big)\Big]\ \\
=&\ (r b_{1}+1) (B R +1)\Big[\ (-Rb_2+rb_1)\overline{\nu}_{21} + (R-r)\overline{\nu}_{22} + (b_1-b_2)\overline{\nu}_{23}\Big]\ > 0\ \ \text{if and only if}\ \ b_1>b_2\ \ \text{and}\ \ r<R<\frac{b_1}{b_2}r, &&
  \end{flalign*}
\end{small}
\vspace{-0.4cm}
\begin{small}
\begin{flalign*}
\mu_3 =&\ b_{1} (r b_{1}+1) (B R +1)\Big[\ (-Rb_2 + rb_1)\Big(BRk_{18}(k_5+k_6+k_{13})k_{3} + (R-r)(B(k_{21}+k_{22}+k_{29})k_{19}(k_{5}+k_{6}+k_{13}) \\
&+ B(k_{21}+k_{22}+k_{29})(k_5+k_6+k_{13})k_{3})  + BRk_{19}(k_5+k_6+k_{13})k_{3} + k_{18}(k_5+k_6+k_{13})k_{3} + k_{19}(k_5+k_6+k_{13})k_{3}\Big) \\
&+ (R-r)\Big( Br(k_{21}+k_{22}+k_{29})k_{19}k_{2}b_1 + Br(k_{21}+k_{22}+k_{29})k_{19}b_1k_{3} + B(k_{21}+k_{22}+k_{29})k_{19}k_{2} + B(k_{21}+k_{22}+k_{29})k_{19}k_{3}\Big)\Big]\ \\ 
 =&\ b_{1} (r b_{1}+1) (B R +1)\Big[\ (-Rb_2 + rb_1)\overline{\nu}_{31} + (R-r)\overline{\nu}_{32}\Big]\ > 0\ \ \text{if and only if}\ \ r<R<\frac{b_1}{b_2}r, && 
   \end{flalign*}
\end{small}
\vspace{-0.5cm}
\begin{small}
\begin{flalign*}
\mu_4 =&\ b_{1} B (k_{21}+k_{22}+k_{29}) (k_{5}+k_{6}+k_{13}) (r b_{1}+1) (R - r) (-R b_{2}+r b_{1}) (B R + 1) k_{19} k_{3} > 0\ \ \text{if and only if}\ \ r<R<\frac{b_1}{b_2}r, &&
\end{flalign*}
\end{small}
\vspace{-0.4cm}
\begin{small}
\begin{align*}
\mu_1\mu_2 - \mu_0\mu_3 > 0\ \ \text{if and only if}\ \ b_1>b_2\ \ \text{and}\ \ r<R<\frac{b_1}{b_2}r
\end{align*}
and
\end{small}
\begin{small}
\begin{align*}
\mu_1\mu_2\mu_3 - \mu_1^2\mu_4 - \mu_0\mu_3^2 > 0\ \ \text{if and only if}\ \ b_{1} > b_{2}\ \ \text{and}\ \ r<R<\frac{b_{1}}{b_{2}}r.
\end{align*}
\end{small}

The proof of the positivity of the two additional conditions, $\mu_{1}\mu_{2} - \mu_{0}\mu_{3}$ and $\mu_1\mu_2\mu_3 - \mu_1^2\mu_4 - \mu_0\mu_3^2$, in the case of the initial mathematical model $(2)$, can be found in the \textit{Supplementary Notes}. This proof is applicable even if we work with the simplified mathematical model $(3)$ or the initial mathematical model $(2)$.

Therefore, the eigenvalues that correspond to the characteristic equation of the Jacobian matrix $J(E_3)$, evaluated at $E_3$ for system $(2)$, are:
\begin{equation*}
\begin{aligned}
    \lambda_1 &= -c_2 = k_7-k_9-k_{10}-k_{14} < 0\ \ \ \text{based on the assumption that $k_9+k_{10}+k_{14} > k_{7}$} \\  \lambda_2 &= -c_3 = -k_{11}-k_{12}-k_{15} < 0,\ \ \lambda_3 = -c_4 = -k_{16} < 0, \\ 
    \lambda_4 &= -C_2 = k_{23}-k_{25}-k_{26}-k_{30} < 0\ \ \ \text{based on the assumption that $k_{25}+k_{26}+k_{30} > k_{23}$} \\ 
    \lambda_5 &= -C_3 = -k_{27}-k_{28}-k_{31} < 0,\ \    \lambda_6 = -C_4 = -k_{32} < 0, \\
    \lambda_7 &,\ \lambda_8,\ \lambda_9\ \text{and}\ \lambda_{10}\ \text{have negative real parts if and only if  $\ b_1>b_2\ $ and $\ r<R<\frac{b_1}{b_2}r\ $}.
\end{aligned}
\end{equation*}
Thus, we can conclude that the steady state $E_{3}$ is also asymptotically stable if and only if $b_{1}>b_{2}$ and $r<R<(b_1/b_2)r$.

\begin{figure}[H]
  \centering
\includegraphics[width=0.75\textwidth, height = 0.1\textwidth]{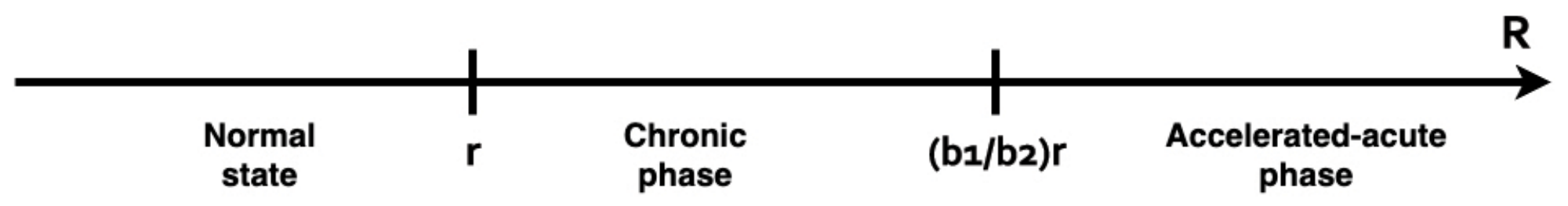}
  \caption{This diagram illustrates the transition from healthy hematopoiesis to the chronic and accelerated-acute stages in myeloid leukemia. In this context, $r$ and $R$ represent the homeostatic quantities of healthy and malignant cycling hematopoietic stem cells, respectively, as defined by equation (\ref{eq:star}). Values of $R$ less than $r$ correspond to the healthy hematopoietic state; values of $R$ between $r$ and $(b_{1}/b_{2})r$ correspond to the chronic phase of leukemia; values of $R$ greater than $(b_{1}/b_{2})r$ characterize the accelerated-acute phase of the disease.} \label{fig:2}
\end{figure}

The above analysis reveals a qualitative change in the behavior of system $(2)$, specifically, a change in its stable equilibrium as the parameter $R$ varies. The specific values of this parameter at which the stable equilibrium changes, known as bifurcation points, are found at $R = r$ and $R = (b_1/b_2)r$. The bifurcation analysis of our mathematical model $(2)$ is illustrated by Figure \ref{fig:2}.

A condition such as $R = r$ or $R = (b_1/b_2)r$ is highly physiologically unstable. Small variations in the parameters can swiftly transition the system between the healthy state and the chronic phase in the case of $R = r$ or between the chronic phase and the accelerated-acute phase in the case of $R = (b_1/b_2)r$. Additionally, from a medical perspective, the situations where $R = r$ and $R = (b_1/b_2)r$ are practically undetectable. The scenario where $R = 0$ could exist in the context of healthy individuals, where the rates $k_{17}$ and $k_{21}+k_{22}+k_{29}$ are equal. However, in the context of a leukemic case, specifically the chronic phase of Chronic Myeloid Leukemia (CML), it is evident that $R > 1$. When a single mutation occurs, we can characterize it by setting the homeostatic amount of leukemic stem cells as $R=1$. 

Therefore, based on this analysis of the local asymptotic stability of the steady states of system $(2)$, we obtain the following result:
\begin{theorem}
    (a) If $R<r$, then $E_{1}$ is the only one steady state that is locally asymptotically stable. \\[5pt]
(b) If $b_{1}>b_{2}$ and $r<R<(b_1/b_2)r$, then $E_{3}$ is the only one steady state  that is locally asymptotically stable. \\[5pt]
(c) If $R>(b_1/b_2)r$, then $E_{2}$ is the only one steady state that is locally asymptotically stable.
\end{theorem} \label{th:2.2}
\begin{remark}
    In all of the three cases of the previous Theorem, the steady state $E_{0}$ is unstable, as can be shown based on the assumptions $k_{1}>k_{5}+k_{6}+k_{13}$ and $k_{17}>k_{21}+k_{22}+k_{29}$.
\end{remark}

\textbf{(c).} \textbf{Global Stability.} It's important to note that the local stability we have analyzed here corresponds to the global stability of the respective steady state of the system $(2)$. Consequently, we can formulate the following remark:

\begin{remark}
    For any positive saturated solution $u=(x_0,x_1,x_2,x_3,x_4,y_0,y_1,y_2,y_3,y_4)$ of system $(2)$ one has:
\begin{equation*}
\begin{aligned}
(i)\ \ u( t) &\rightarrow E_{1} \ \ \text{as\
\ }t\rightarrow +\infty ,\ \ \text{in case \ }R<r; \\
\text{\textit{(ii)}\ \ }u( t) &\rightarrow E_{3} \ \ \text{as\ \ }t\rightarrow +\infty ,\ \ \text{in case \ }%
r<R<( b_{1}/b_{2}) r; \\
\text{\textit{(iii)}\ \ }u( t) &\rightarrow E_{2} \ \ \text{%
as\ \ }t\rightarrow +\infty ,\ \ \text{in case \ }( b_{1}/b_{2})
r<R.
\end{aligned}
\end{equation*}
\end{remark}

The \textit{proof} of the global stability statement follows a similar approach to the one presented for system (\ref{s}) in the paper by Parajdi and Precup \cite{pp}, with additional and more detailed explanations.

\section{Numerical Simulation of the Model} \label{Section:3}

\subsection{Parameter Estimation}

\noindent The model's parameters are influenced by numerous biophysical and biochemical mechanisms, making precise estimation nearly impractical. Instead, one can anticipate obtaining parameters from an estimation procedure with confidence intervals. In the context of our qualitative analysis, precise parameter estimation is not crucial because it does not alter the behavior of our model. Specifically, it does not affect the stability of the steady states, as demonstrated earlier in Section 2. However, parameter estimation becomes essential when the model is used for real-time predictions and applications to individual patients.

According to the paper of Foo \textit{et al.} \cite{fdc}, in the equilibrium state, the total number of healthy and malignant stem cells (cycling plus quiescent) is $x_{0}^{\ast}+x_{1}^{\ast} = 10^6$ for healthy cells and $y_{0}^{\ast}+y_{1}^{\ast} = 10^7$ for malignant cells. Moreover, the number of healthy progenitor cells is $x_{2}^{\ast} = 10^8$, the number of healthy differentiated cells is $x_{3}^{\ast}=10^{10}$ and the number of healthy terminally differentiated cells is $x_{4}^{\ast} = 10^{12}$. Next, we assume that healthy cycling stem cells undergo symmetric division (self-renewal) every $36$ days (according to Barile \textit{et al.} \cite{bbf}, and Wilson \textit{et al.} \cite{wlowb}), undergo asymmetric division every $200$ days (according to Barile \textit{et al.} \cite{bbf}) and are eliminated every $333$ days (according to Foo \textit{et al.} \cite{fdc}). Consequently, we can obtain the symmetric self-renewal rate, the asymmetric division rate as well as the death (apoptosis) rate of healthy cycling stem cells (per capita, per day) as follows: $k_{1} = 1/36 \simeq 0.028\ $, $k_{4} = 1/200 = 0.005\ $ and $\ k_{13} = 1/333 \simeq 0.003$. Here, $k_{1} = (0.028/0.005)k_{4} = 5.6k_{4}$, which means $k_{1}$ is $5.6$ times greater than $k_{4}$. According to the paper of Barile \textit{et al.} \cite{bbf}, healthy cycling HSCs exhibit symmetric self-renewal every $22$ to $36$ days (at a rate of 2.8--4.6\%), undergo asymmetric division every $167$ to $10,000$ days (at a rate of 0.01--0.6\%) and engage in direct differentiation every $200$ to $10,000$ days (at a rate of 0.01--0.5\%). Additionally according to the paper of Barile \textit{et al.} \cite{bbf}, we can compute the range of $k_{6}$, in which the healthy cycling HSCs undergo symmetric differentiation by fixing: the rate at which healthy cycling HSCs undergo asymmetric division $k_{4} = 0.005$, the rate at which healthy cycling HSCs undergo direct differentiation $k_{5} = 1/1000 = 0.001$ and the total differentiation rate $\alpha$ between 0.7--1.1\%. Using the formula where the total differentiation rate $\alpha = 2\rho + \gamma + \mu$, where $\rho\ (:= k_{6})$ represents the rate at which healthy cycling HSCs undergo symmetric differentiation, $\gamma\ (:=k_{4})$ represents the rate at which healthy cycling HSCs undergo asymmetric division and $\mu\ (:=k_{5})$ represents the rate at which healthy cycling HSCs undergo direct differentiation, we can compute a reference interval of $\rho$, which is 0.05-0.25\%. Therefore, the healthy cycling HSCs undergo symmetric differentiation every $400$ to $2,000$ days. The rate at which healthy cycling HSCs undergo asymmetric division ($k_{4}$), direct differentiation ($ k_{5}$), symmetric differentiation ($k_{6}$) and the total differentiation rate ($\alpha$) were obtained from an extensive dataset obtained from laboratory experiments conducted on mice (see Barile \textit{et al.} \cite{bbf}). In our simulations, we assume that $k_{6} = 1/400 = 0.0025$. At first approximation, data obtained from laboratory mice may appear consistent with those observed in humans. However, several studies have highlighted significant differences in hematopoietic properties between species (see the work of Parekh and Crooks \cite{pc}). For example, while hematopoiesis originates in the yolk sac during the embryonic period in both humans and mice (see Palis and Yoder \cite{py}), multi-omics analyses and integrative cross-species transcriptome data have shown that the development of the granulocytic series and myelopoiesis occurs earlier in primates, including humans, compared to mice (see Du \textit{et al.} \cite{dlgll}). These findings highlight important species-specific differences in hematopoietic processes.

In our study, we focus on the global properties of the hematopoietic system, as reflected in the model parameters, rather than on specific biochemical or biophysical characteristics that may differ between species. Consequently, while we acknowledge the differences in hematopoietic properties between humans and mice, we believe that the conclusions of our study remain robust and are not significantly affected by the use of parameter data derived from both species. Our approach emphasizes the overarching principles of the system, which are sufficiently general to provide meaningful insights despite the interspecies variability. Regarding the lifespan of HSCs, the study by Barile \textit{et al.} \cite{bbf} proposes two distinct apoptosis ranges: one at every $200$ to $1000$ days (at a rate of 0.1--0.5\%), and another at every $25$ to $50$ days (at a rate of 2--4\%). The parameter $G=k_{3}/(k_{2}+k_{3})$ specifies the frequency of healthy cycling stem cells when the system is in equilibrium. If we assume that $k_{2}=k_{18}$ and $k_{3}=k_{19}$, then $G$ also represents the equilibrium fraction of malignant cycling stem cells. In our simulations, we have chosen specific values for $k_{2},\ k_{3},\ k_{18}$ and $k_{19}$, such as $k_{2}=k_{18}=0.0001$ and $k_{3}=k_{19}=0.0009$, resulting in an initial equilibrium frequency of $90\%$ for both healthy and malignant cycling stem cells, corresponding to $G=0.9$. However, based on Foo's work \cite{fdc}, we can choose different values for $k_{2}=k_{18}=0.0009$ and $k_{3}=k_{19}=0.0001$, leading to an initial equilibrium frequency of $10\%$ for both healthy and malignant cycling stem cells, denoted by $G=0.1$. Alternatively, we can choose values like $k_{2}=k_{18}=0.0005$ and $k_{3}=k_{19}=0.0005$, resulting in an initial equilibrium frequency of $50\%$ for both healthy and malignant cycling stem cells, with $G=0.5$. The rates of symmetric and asymmetric divisions among intermediate healthy progenitor cells can vary depending on the specific cell type, microenvironment and the body's demand for different blood cell types. These cells may divide over a range of time intervals, spanning from hours to several days. The exact frequencies of symmetric and asymmetric divisions for intermediate healthy progenitor cells are influenced by various factors, including the type of cell, the presence of growth factors, cytokines and the body's requirements for specific blood cell production. Hence, in our simulations, we assume that $k_{7} = (2 \times 5.6)k_{8} = 11.2k_{8}$, which means $k_{7}$ is $11.2$ times greater than $k_{8}$.

In the equilibrium state, for healthy quiescent cells, we have $k_{2}x_{0}^{\ast}-k_{3}x_{1}^{\ast}=0$, which leads to $x_{0}^{\ast}/x_{1}^{\ast} = k_{3}/k_{2} = 9$. Therefore, using the fact that the total number of healthy stem cells is $x_{0}^{\ast}+x_{1}^{\ast} = 10^6$, we can compute the number of healthy quiescent stem cells, which is $x_{1}^{\ast} = 10^5$. Returning to the equation $x_{0}^{\ast}+x_{1}^{\ast} = 10^6$, now that we know the number of healthy quiescent stem cells, we can compute the number of healthy cycling stem cells, which is $x_{0}^{\ast} = 9 \times 10^5$. Similarly, in the equilibrium state for the malignant quiescent stem cells, we have $k_{18}y_{0}^{\ast}-k_{19}y_{1}^{\ast}=0$, which leads to $y_{0}^{\ast}/y_{1}^{\ast} = k_{19}/k_{18} = 9$. Using the fact that the total number of malignant stem cells is $y_{0}^{\ast}+y_{1}^{\ast} = 10^7$, we can compute the number of malignant quiescent stem cells, which is $y_{1}^{\ast} = 10^6$. Now that we know the number of malignant quiescent stem cells, we can use the equation $y_{0}^{\ast}+y_{1}^{\ast} = 10^7$ to compute the number of malignant cycling stem cells, which is $y_{0}^{\ast} = 9 \times 10^6$. Next, in the equilibrium state, for healthy progenitor cells, we have $(k_{4}+k_{5}+2k_{6})x_{0}^{\ast}-(-k_{7}+k_{9}+k_{10}+k_{14})x_{2}^{\ast}=0$ which leads to $(k_{4}+k_{5}+2k_{6})/(-k_{7}+k_{9}+k_{10}+k_{14}) = x_{2}^{\ast}/x_{0}^{\ast} = 10^3/9$, for healthy differentiated cells, we have $(k_{8}+k_{9}+2k_{10})x_{2}^{\ast}-(k_{11}+k_{12}+k_{15})x_{3}^{\ast} = 0$, which leads to $(k_{8}+k_{9}+2k_{10})/(k_{11}+k_{12}+k_{15}) = x_{3}^{\ast}/x_{2}^{\ast} = 10^2$ and finally for healthy terminally differentiated cells, we have $(k_{11}+2k_{12})x_{3}^{\ast}-k_{16}x_{4}^{\ast} = 0$, which leads to $(k_{11}+2k_{12})/k_{16} = x_{4}^{\ast}/x_{3}^{\ast} = 10^2$. Therefore, if we assume the death rates of healthy cells to be $k_{13} = 0.003$, $k_{14} = 0.008$, $k_{15}=0.05$, $k_{16} = 1$ (see Foo \textit{et al.} \cite{fdc}), then we have:

$\diamond$ $k_{11}+2k_{12} = 10^2$, if we suppose that $k_{11} > k_{12}$ (most of the healthy and malignant differentiated cells directly differentiate into terminally differentiated cells instead of undergoing symmetric differentiation) then we assume $k_{11}=2k_{12}$ and using $k_{11}+2k_{12} = 10^2$ we obtain $k_{12} = 25$ and $k_{11} = 50$;

$\diamond$ $(k_{8}+k_{9}+2k_{10})/(k_{11}+k_{12}+k_{15}) = 10^2$ thus we obtain that $k_{8} + k_{9} + 2k_{10} = 7505$. Therefore, if we suppose that $k_{9} > k_{10}$ (most of the healthy and malignant progenitor cells directly differentiate into differentiated cells instead of undergoing symmetric differentiation) then we assume $k_{9} = 2k_{10}$ we have $k_{8}+4k_{10} = 7505$;

$\diamond$ $(k_{4}+k_{5}+2k_{6})/(-k_{7}+k_{9}+k_{10}+k_{14}) = 10^3/9$ we obtain that $-k_{7} + k_{9} + k_{10} = -0.007901\ $ then if we use $k_{7} = 11.2k_{8}$ and $k_{9} = 2k_{10}$ we have $-11.2k_{8} + 3k_{10} = -0.007901$.

Therefore, from this two equations $k_{8}+4k_{10} = 7505$ and $-11.2k_{8} + 3k_{10} = -0.007901$ we can compute $k_{8} = 471.0257658$, $\ k_{10} = 1758.493559$ and using the fact that $k_{7} = 11.2k_{8}$ and $k_{9} = 2k_{10}$ we can compute $k_{7} = 5275.488577$, $\ k_{9} = 3516.987118$. The ratio $b_1/b_2$ allows for the possibility that the malignant cycling stem cells are less sensitive to the bone marrow microenvironment than the healthy cycling stem cells. The parameter $b_1$, which represents the sensitivity of healthy cycling stem cells to the bone marrow microenvironment, can be estimated from the expression of the homeostatic amount of healthy cycling stem cells $r$, as follows: $b_1 = (k_1-k_5-k_6-k_{13})/(r(k_5+k_6+k_{13})) = 3.675213676 \times 10^{-6}$, where $r = x_{0}^{\ast} = 9 \times 10^5$.

The parameters $B$, $b_{2}$ and $k_{j_{M}}$ where $j_{M} = 17\ldots32$, can vary from patient to patient, and so the parameter $R$ and the length of the chronicity interval, denoted as [$r, (b_{1}/b_{2}) r$]. For our numerical simulations, we have chosen a value of $2$ for $b_{1}/b_{2}$, resulting in $b_{2} = 1.837606838 \times 10^{-6}$. Additionally, following the approach in Dingli and Michor \cite{dm}, we assume that parameter $B$ is approximately half of $b_{2}$, resulting in $B = 9.188034190 \times 10^{-7}$. Regarding the parameters $k_{j_{M}}$, where $j_{M} = 17\ldots32$, in our simulations, we assume that $k_{17} = 2k_{1}$, $k_{18} = k_{2}$, $k_{19} = k_{3}$, $k_{20} = 2k_{4}$, $k_{21} = 2k_{5}$, $k_{22} = 2k_{6}$, $k_{23} = 2k_{7}$, $k_{24} = 2k_{8}$, $k_{25} = 2k_{9}$, $k_{26} = 2k_{10}$, $k_{27} = k_{11}$, $k_{28} = k_{12}$ and $k_{29} > k_{13}$, $k_{30} = 2k_{14}$, $k_{31} = k_{15}$, $k_{32} = k_{16}$.

\subsection{Numerical Simulations}

\noindent We perform numerical simulations of the nonlinear system $(2)$ to investigate the behavior of various cell types in cloned hematopoiesis, including healthy and malignant cycling stem cells, quiescent stem cells, progenitor cells, differentiated cells and terminally differentiated cells, across three distinct scenarios: $R < r$ (healthy states), $r < R < (b_{1}/b_{2})r$ (chronic state) and $(b_{1}/b_{2})r < R$ (accelerated-acute state). We restrict our simulations to the scenario where $k_{17} > k_{1}$, $k_{18} = k_{2}$, $k_{19} = k_{3}$, $k_{20} > k_{4}$, $k_{21} > k_{5}$, $k_{22} > k_{6}$, $k_{23} > k_{7}$, $k_{24} > k_{8}$, $k_{25} > k_{9}$, $k_{26} > k_{10}$, $k_{27} = k_{11}$, $k_{28} = k_{12}$ and $k_{30} > k_{14}$, $k_{31} = k_{15}$, $k_{32} = k_{16}$. In our simulations, we explore the transition process from healthy hematopoiesis to the chronic and accelerated-acute stages. We assume that this transition is primarily influenced by the parameter $k_{29}$. Therefore, we conduct simulations with varying values of this parameter. It's important to note that the transition between these three stages can also be influenced by the modification of the homeostatic amounts of healthy and malignant stem cells $r$ and $R$, which, in turn, corresponds to the variation of parameters $k_{1}$, $k_{5}$, $k_{6}$, $k_{13}$, $b_{1}$, $k_{17}$, $k_{21}$, $k_{22}$, $k_{29}$ and $B$. Furthermore, we make the assumption that malignant cycling stem cells exhibit lower sensitivity to environmental crowding compared to healthy cycling stem cells, specifically, $b_{1} > b_{2} > B$. Additionally, we ensure that the two conditions for non-negative equilibrium points are satisfied, which are $k_{9}+k_{10}+k_{14}-k_{7} > 0$ (equivalent to $0.000100 > 0$) and $k_{25}+k_{26}+k_{30}-k_{23} > 0$ (equivalent to $0.000204 > 0$). This ensures that the equilibrium points remain non-negative. Next, we check that $k_{1}-k_{5}-k_{6}-k_{13} > 0$ (equivalent to $0.0215 > 0$) and $k_{17}-k_{21}-k_{22}-k_{29} > 0$ (equivalent to $0.0240 > 0$). We can point out that similar simulations can be performed with different sets of conditions on the parameters.
\vspace{0.2cm}

\begin{figure}[th!]
\centering
\subfigure[~Healthy state]{\includegraphics[width=5.25cm]{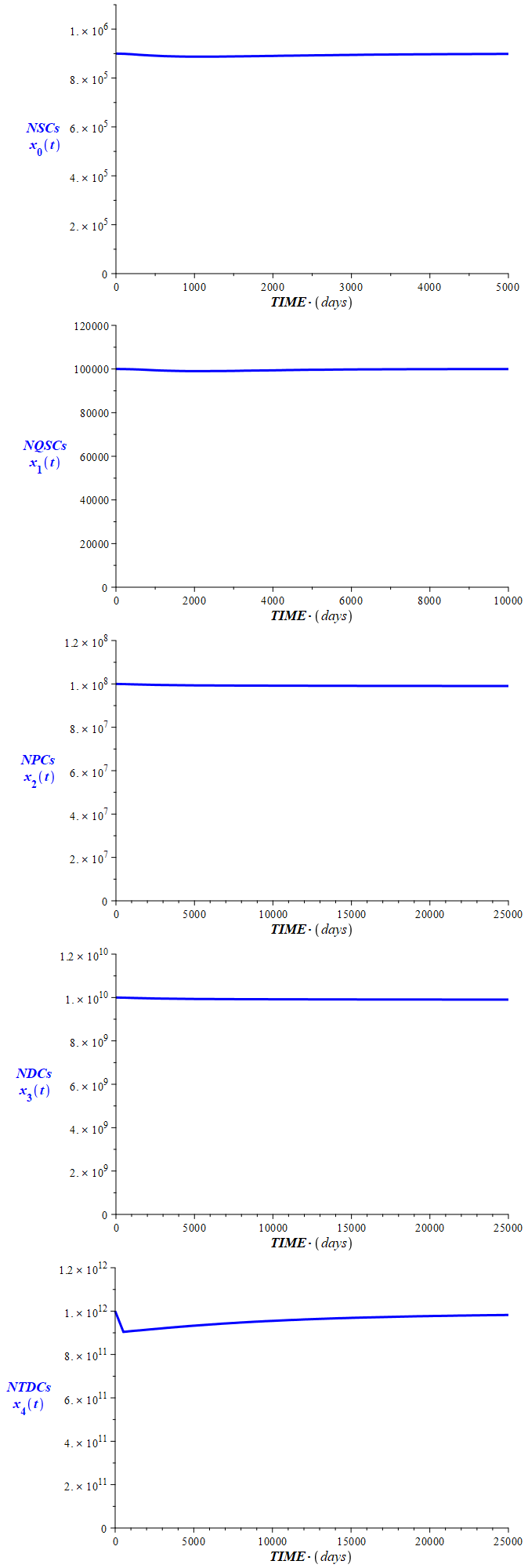}} %
\subfigure[~Chronic
phase]{\includegraphics[width=5.55cm]{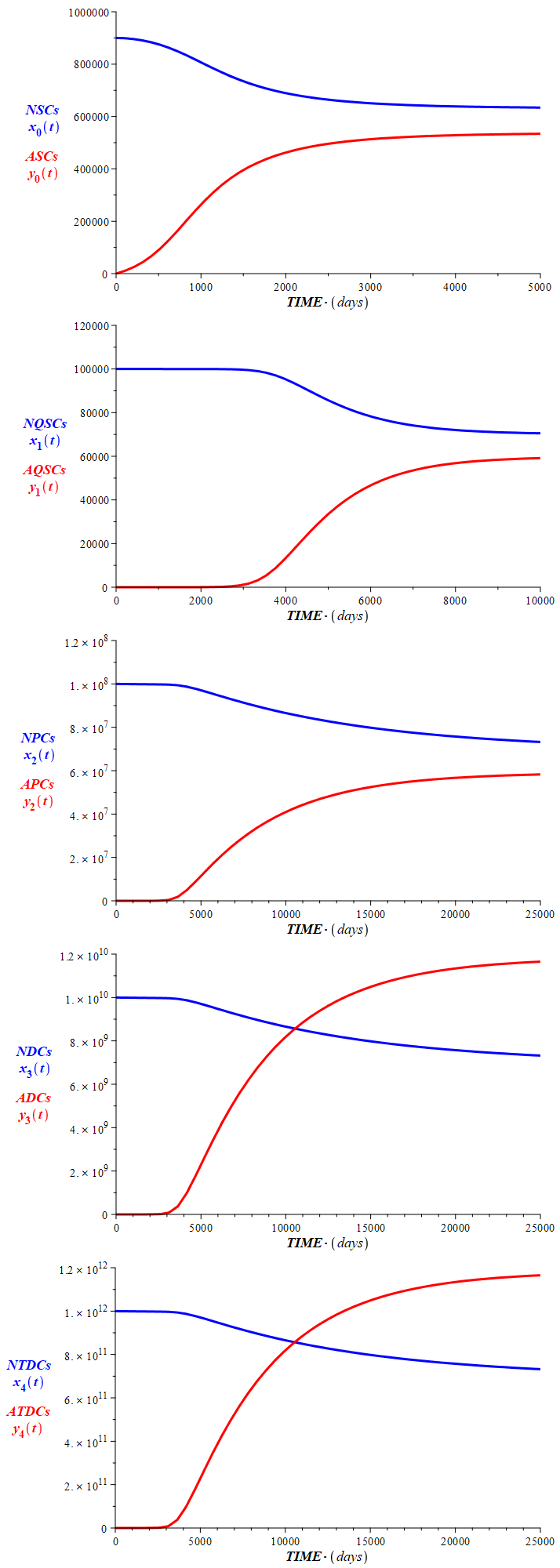}}
\subfigure[~Accelerated-acute
phase]{\includegraphics[width=5.55cm]{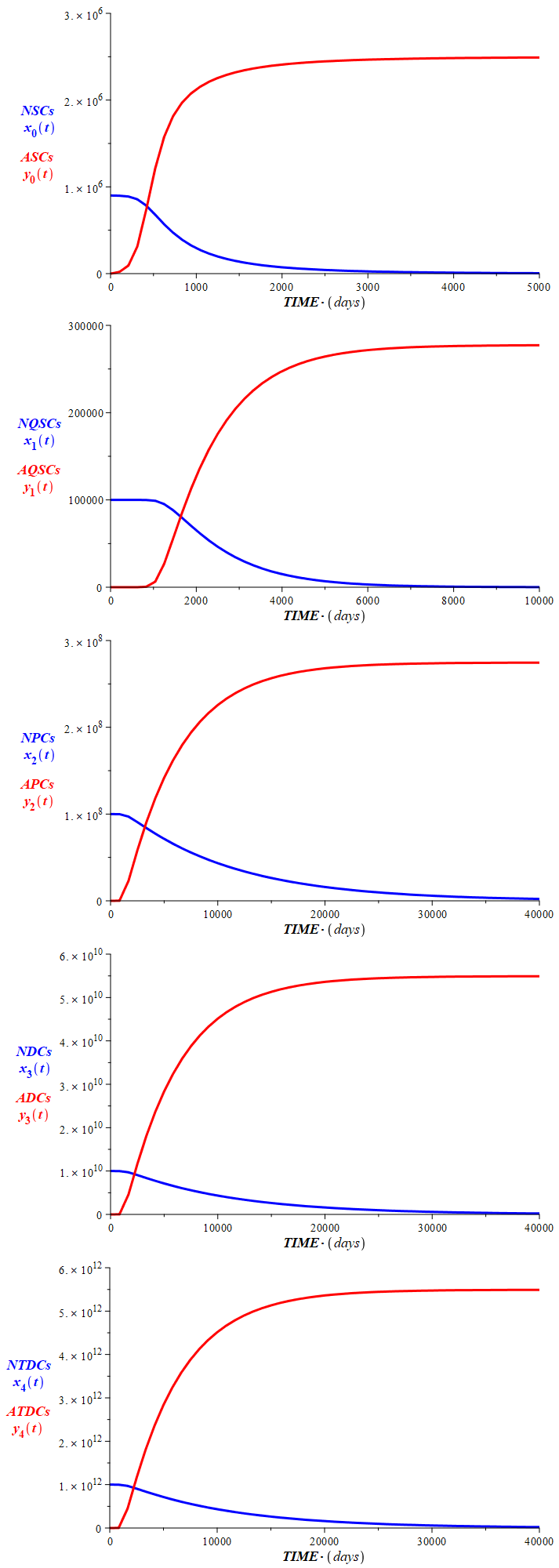}}
\caption{{\protect\small Behavior of healthy and malignant
(leukemic) cell populations. Initial conditions: (\textbf{a}) healthy cell populations: cycling stem cells $x_{0}(0)=9\times10^5$, quiescent stem cells $x_{1}(0)=10^5
$, progenitor cells $x_{2}(0)=10^8$, differentiated cells $x_{3}(0)=10^{10}$, terminally differentiated cells $x_{4}(0)=10^{12}$; (\textbf{b}) - (\textbf{c}) healthy and malignant cells: cycling stem cells $x_{0}(0)=9\times10^5$ and $y_{0}(0)=1$, quiescent stem cells $x_{1}(0)=10^5$ and $y_{1}(0)=1$, progenitor cells $x_{2}(0)=10^8$ and $y_{2}(0)=1$, differentiated cells $x_{3}(0)=10^{10}$ and $y_{3}(0)=1$, terminally differentiated cells $x_{4}(0)=10^{12}$ and $y_{4}(0)=1$.
}}
\label{fig.3}
\end{figure}

Figure \ref{fig.3}(a) shows the time-dependent behavior of various cell populations, including healthy: cycling stem cells ($T=5000$ days), quiescent stem cells ($T=10000$ days), progenitor cells ($T=25000$ days), differentiated cells ($T=25000$ days), terminally differentiated cells ($T=25000$ days). These simulations are based on the parameter values provided in the first column of Table \ref{tab.2}, which correspond to the healthy hematopoietic state where $R<r$. Therefore, the healthy cycling stem cell population, denoted as $x_{0}(t)$ (blue solid line), tends to approach the value $r = 899999.9998$. Similarly, the healthy quiescent stem cell population, $x_{1}(t)$ (blue solid line), tends to reach the value $r(k_{2}/k_{3}) = 99999.99998$. The healthy progenitor cell population, $x_{2}(t)$ (blue solid line), tends to approach the value $r[\,(k_4+k_5+2k_6)/(-k_7+k_9+k_{10}+k_{14})]\, = 9.899999998 \times 10^7$. The healthy differentiated cell population, $x_{3}(t)$ (blue solid line), tends to stabilize at the value $r[\,((k_{4}+k_{5}+2k_{6})(k_{8}+k_{9}+2k_{10}))/((-k_{7}+k_{9}+k_{10}+k_{14})(k_{11}+k_{12}+k_{15}))]\, = 9.900000000 \times 10^9$. Lastly, the healthy terminally differentiated cell population, $x_{4}(t)$ (blue solid line), tends to reach the value $r[\,((k_{4}+k_{5}+2k_{6})(k_{8}+k_{9}+2k_{10})(k_{11}+2k_{12}))/((-k_{7}+k_{9}+k_{10}+k_{14})(k_{11}+k_{12}+k_{15})k_{16})]\, = 9.900000000 \times 10^{11}$. All of these dynamics occur in the absence of any malignant cell populations. Notably, in the context of healthy hematopoiesis ($R<r$), mutations that arise over time are at times naturally eradicated by random events. In this scenario, the population of malignant cells tends towards zero. Biologically, this phenomenon of mutant or malignant cell extinction due to random events has been explained and demonstrated in several studies involving stem cell lineages (see Jilkine and Gutenkunst \cite{jg}, Driessens \textit{et al.} \cite{dbcsb}, Klein and Simons \cite{ks}, Lopez-Garcia \textit{et al.} \cite{lksw} and Snippert \textit{et al.} \cite{sfseb}). However, it's important to emphasize that the natural extinction of malignant cells is a rare phenomenon, and consequently, it cannot be guaranteed. In cases where this does not occur, leukemia may advance.
\vspace{0.2cm}

\begin{table}[H]
\scriptsize
\centering
\begin{tabular}{|c|ccc|c|c|}
\hline
\textbf{Figure}                  & $3(a)$        & $3(b)$                                      & $3(c)$   &  $\mathbf{Source/Reference}$  & $\mathbf{Assumed/Calculated}$    \\ \hline
$\mathbf{k_{1}}$               & $0.028$       & $0.028$                                     & $0.028$    & see Barile \textit{et al.} \cite{bbf} and Wilson \textit{et al.} \cite{wlowb}   &  assumed value  \\
$\mathbf{b_{1}\times 10^{-6}}$ & $3.675213676$ & $3.675213676$                               & $3.675213676$      & see Parajdi \textit{et al.} \cite{ppbt} and Neiman \cite{n}       &  calculated value \\
$\mathbf{b_{2}\times 10^{-6}}$ & $1.837606838$ & $1.837606838$                               & $1.837606838$      &  see Parajdi \textit{et al.} \cite{ppbt} and Neiman \cite{n}      &  calculated value \\
$\mathbf{k_{2}=k_{18}}$              & $0.0001$      & $0.0001$                                    & $0.0001$     & see Foo \textit{et al.} \cite{fdc}     &  assumed value \\
$\mathbf{k_{3}=k_{19}}$              & $0.0009$      & $0.0009$                                    & $0.0009$     & see Foo \textit{et al.} \cite{fdc}     &  assumed value \\
$\mathbf{k_{4}}$               & $0.005$       & $0.005$                                     & $0.005$    &  see Barile \textit{et al.} \cite{bbf}   &  assumed value  \\
$\mathbf{k_{5}}$              & $0.001$       & $0.001$                                     & $0.001$      &  see Barile \textit{et al.} \cite{bbf}   &  assumed value  \\
$\mathbf{k_{6}}$              & $0.0025$      & $0.0025$                                    & $0.0025$     &  see Barile \textit{et al.} \cite{bbf}   & calculated value \\
$\mathbf{k_{7}}$              & $5275.488577$ & $5275.488577$                               & $5275.488577$    &  see Foo \textit{et al.} \cite{fdc}      &  calculated value\\
$\mathbf{k_{8}}$               & $471.0257658$ & $471.0257658$                               & $471.0257658$   &  see Foo \textit{et al.} \cite{fdc}      &  calculated value \\
$\mathbf{k_{9}}$              & $3516.987118$ & $3516.987118$                               & $3516.987118$    &  see Foo \textit{et al.} \cite{fdc}      &  calculated value \\
$\mathbf{k_{10}}$              & $1758.493559$ & $1758.493559$                               & $1758.493559$   &  see Foo \textit{et al.} \cite{fdc}      &  calculated value \\
$\mathbf{k_{11}=k_{27}}$              & $50$          & $50$                                        & $50$         &  see Foo \textit{et al.} \cite{fdc}             &  calculated value \\
$\mathbf{k_{12}=k_{28}}$              & $25$          & $25$                                        & $25$         &  see Foo \textit{et al.} \cite{fdc}             &  calculated value \\
$\mathbf{k_{13}}$              & $0.003$       & $0.003$                                     & $0.003$      & see Foo \textit{et al.} \cite{fdc}     & assumed value \\
$\mathbf{k_{14}}$              & $0.008$       & $0.008$                                     & $0.008$      & see Foo \textit{et al.} \cite{fdc}     & assumed value \\
$\mathbf{k_{15}=k_{31}}$             & $0.05$        & $0.05$                                      & $0.05$       & see Foo \textit{et al.} \cite{fdc}     & assumed value \\
$\mathbf{k_{16}=k_{32}}$              & $1$           & $1$                                         & $1$          & see Foo \textit{et al.} \cite{fdc}     & assumed value \\
$\mathbf{k_{17}}$              & $0.056$       & $0.056$                                     & $0.056$      & see Barile \textit{et al.} \cite{bbf} and Dingli, Michor \cite{dm}     & assumed value \\
$\mathbf{B\times 10^{-7}}$     & $9.188034190$ & $9.188034190$                               & $9.188034190$      &  see Parajdi \textit{et al.} \cite{ppbt} and Dingli, Michor \cite{dm}     &  calculated value \\
$\mathbf{k_{20}}$              & $0.010$       & $0.010$                                     & $0.010$      & see Foo \textit{et al.} \cite{fdc}     & assumed value \\
$\mathbf{k_{21}}$               & $0.002$       & $0.002$                                     & $0.002$      & see Foo \textit{et al.} \cite{fdc}     & assumed value \\
$\mathbf{k_{22}}$               & $0.0050$      & $0.0050$                                    & $0.0050$     & see Foo \textit{et al.} \cite{fdc}     & assumed value \\
$\mathbf{k_{23}}$              & $10550.97715$ & $10550.97715$                               & $10550.97715$    & see Foo \textit{et al.} \cite{fdc}     & assumed value \\
$\mathbf{k_{24}}$              & $942.0515316$ & $942.0515316$                               & $942.0515316$    & see Foo \textit{et al.} \cite{fdc}     & assumed value \\
$\mathbf{k_{25}}$               & $7033.974236$ & $7033.974236$                               & $7033.974236$   & see Foo \textit{et al.} \cite{fdc}     & assumed value \\
$\mathbf{k_{26}}$              & $3516.987118$ & $3516.987118$                               & $3516.987118$    & see Foo \textit{et al.} \cite{fdc}     & assumed value \\
$\mathbf{k_{29}}$               & $0.025$       & $0.02$                                      & $0.01$       & see Parajdi \textit{et al.} \cite{ppbt}       &  assumed value \\
$\mathbf{k_{30}}$               & $0.016$       & $0.016$                                     & $0.016$      & see Parajdi \textit{et al.} \cite{ppbt}       &  assumed value \\ \hline
\textbf{\textit{S--S}}                 & $E_{1}$       & $E_{3}$ & $E_{2}$   &        &   \\ \hline
\end{tabular}
\caption{{\protect\small Parameter values for simulations. \textit{S-S} = steady state. An assumed value indicates that the value was taken from a specific reference. A calculated value means that it was determined by solving a mathematical inverse problem at a steady state.}} \label{tab.2}
\end{table}

Figure \ref{fig.3}(b) shows the time-dependent behavior of various cell populations, including healthy and malignant: cycling stem cells ($T=5000$ days), quiescent stem cells ($T=10000$ days), progenitor cells ($T=25000$ days), differentiated cells ($T=25000$ days) and terminally differentiated cells ($T=25000$ days). These simulations are based on the parameter values provided in the second column of Table \ref{tab.2}, which correspond to the chronic state where $r<R<(b_{1}/b_{2})r$. Thus, the healthy and malignant cycling stem cell populations, denoted as $x_{0}(t)$ (blue solid line) and $y_{0}(t)$ (red solid line) tend toward $\overline{x}^* = 631007.7523$ and $\overline{y}^* = 537984.4956$, respectively. Similarly, the healthy and malignant quiescent stem cell populations, $x_{1}(t)$ (blue solid line) and $y_{1}(t)$ (red solid line), tend to reach the values $\overline{x}^*(k_{2}/k_{3}) = 70111.97248$ and $\overline{y}^*(k_{18}/k_{19}) = 59776.05507$, respectively. The healthy and malignant progenitor cell populations, $x_{2}(t)$ (blue solid line) and $y_{2}(t)$ (red solid line), tend to approach the values $\overline{x}^*[\,(k_4+k_5+2k_6)/(-k_7+k_9+k_{10}+k_{14})]\, = 6.941085275 \times 10^7$ and $\overline{y}^*[\,(k_{20}+k_{21}+2k_{22})/(-k_{23}+k_{25}+k_{26}+k_{30})]\, = 5.801793578\times10^7$, respectively. The healthy and malignant differentiated cell populations, $x_{3}(t)$ (blue solid line) and $y_{3}(t)$ (red solid line) tend to stabilize at the values $\overline{x}^*[\,((k_4+k_5+2k_6)(k_8+k_9+2k_{10}))/((-k_7+k_9+k_{10}+k_{14})(k_{11}+k_{12}+k_{15}))]\, = 6.941085276 \times 10^9$ and $\overline{y}^*[\,((k_{20}+k_{21}+2k_{22})(k_{24}+k_{25}+2k_{26}))/((-k_{23}+k_{25}+k_{26}+k_{30})(k_{27}+k_{28}+k_{31}))]\, = 1.160358716\times10^{10}$, respectively. Lastly, the healthy and malignant terminally differentiated cell populations, $x_{4}(t)$ (blue solid line) and $y_{4}(t)$ (red solid line), tend to reach the value $\overline{x}^*[\,((k_4+k_5+2k_6)(k_8+k_9+2k_{10})(k_{11}+2k_{12}))/((-k_7+k_9+k_{10}+k_{14})(k_{11}+k_{12}+k_{15})k_{16})]\, = 6.941085276\times10^{11}$ and $\overline{y}^*[\,((k_{20}+k_{21}+2k_{22})(k_{24}+k_{25}+2k_{26})(k_{27}+2k_{28}))/((-k_{23}+k_{25}+k_{26}+k_{30})(k_{27}+k_{28}+k_{31})k_{32})]\, = 1.160358716\times10^{12}$, respectively.
\vspace{0.2cm}

Figure \ref{fig.3}(c) shows the time-dependent behavior of various cell populations, including healthy and malignant: cycling stem cells ($T=5000$ days), quiescent stem cells ($T=10000$ days), progenitor cells ($T=40000$ days), differentiated cells ($T=40000$ days) and terminally differentiated cells ($T=40000$ days). These simulations are based on the parameter values provided in the third column of Table \ref{tab.2}, which correspond to the accelerated-acute state where $R>(b_{1}/b_{2})r$. Therefore, the healthy cycling stem cell population, denoted as $x_{0}(t)$ (blue solid line), tends towards $0$, while the malignant cycling stem cell population, denoted as $y_{0}(t)$ (red solid line), tends to the value $R = 2.496853625 \times 10^6$. Similarly, the healthy quiescent stem cell population, $x_{1}(t)$ (blue solid line), tends toward $0$, while the malignant quiescent stem cell population, $y_{1}(t)$ (red solid line), tends to reach the value $R(k_{18}/k_{19}) = 277428.1806$. The healthy progenitor cell population, $x_{2}(t)$ (blue solid line), tends towards $0$, while the malignant progenitor cell population, $y_{2}(t)$ (red solid line), tends to approach the value $R[\,(k_{20}+k_{21}+2k_{22})/(-k_{23}+k_{25}+k_{26}+k_{30})]\, = 2.692685282\times10^8$. The healthy differentiated cell population, $x_{3}(t)$ (blue solid line), tends towards $0$, while the malignant differentiated cell population, $y_{3}(t)$ (red solid line), tends to stabilize at the value $R[\,((k_{20}+k_{21}+2k_{22})(k_{24}+k_{25}+2k_{26}))/((-k_{23}+k_{25}+k_{26}+k_{30})(k_{27}+k_{28}+k_{31}))]\, = 5.385370564\times10^{10}$. Lastly, the healthy terminally differentiated cell population, $x_{4}(t)$ (blue solid line), tends towards $0$, while the malignant terminally differentiated cell population, $y_{4}(t)$ (red solid line), tends to reach the value $R[\,((k_{20}+k_{21}+2k_{22})(k_{24}+k_{25}+2k_{26})(k_{27}+2k_{28}))/((-k_{23}+k_{25}+k_{26}+k_{30})(k_{27}+k_{28}+k_{31})k_{32})]\, = 5.385370564\times10^{12}$.
\vspace{0.2cm}

Figure \ref{fig.3}(a)–(c) illustrates a parallelism in the behaviors of both healthy and malignant (leukemic) cell populations across all four cell compartments.

\section{Discussion and Conclusions} \label{Section:4}

\noindent To highlight the conclusions derived from the preceding discussions, it can be affirmed that the relationship $\ R < r\ $ characterizes the \textit{healthy hematopoiesis}; the relationship $\ r<R<(b_{1}/b_{2})r\ $ signifies the \textit{chronic phase} of leukemia; while, the inequality $\ R>(b_{1}/b_{2})r\ $ characterizes the \textit{accelerated-acute phase} of the disease. In the case where $%
R<r,$ the healthy cell populations exhibit the following behaviors: cycling stem cells $x_{0}\left( t\right)$ approach the equilibrium abundance $r$ (healthy homeostatic state); quiescent stem cells $x_{1}\left(t \right)$ approach the equilibrium abundance $r(k_{2}/k_{3})$; progenitor cells $x_{2}\left(t \right)$ approach the equilibrium abundance $r[\, (k_{4}+k_{5}+2k_{6})/(-k_{7}+k_{9}+k_{10}+k_{14})]\, $; differentiated cells $x_{3}\left(t \right)$ approach the equilibrium abundance $r[\, (k_{4}+k_{5}+2k_{6})(k_8+k_9+2k_{10})/(-k_{7}+k_{9}+k_{10}+k_{14})(k_{11}+k_{12}+k_{15})]\, $; and terminally differentiated cells $x_{4}\left(t \right)$ approach the equilibrium abundance $r[\, (k_{4}+k_{5}+2k_{6})(k_8+k_9+2k_{10})(k_{11}+2k_{12})/(-k_{7}+k_{9}+k_{10}+k_{14})(k_{11}+k_{12}+k_{15})k_{16}]\, $, while all malignant (leukemic) cell populations (cycling stem cells $y_{0}\left( t\right)$; quiescent stem cells $y_{1}\left(t \right)$; progenitor cells $y_{2}\left(t \right)$; differentiated cells $y_{3}\left(t \right)$; and terminally differentiated cells $y_{4}\left(t \right)$) tend to zero. 

In the case where $r<R<(b_{1}/b_{2})r,$ both healthy and malignant (leukemic) cycling stem cell populations $x_{0}\left( t\right)$ and $y_{0}\left( t\right)$ approach their equilibrium abundances $\overline{x}^{\ast}$ and $\overline{y}^{\ast},$ respectively; healthy and malignant quiescent stem cell populations $x_{1}(t)$ and $y_{1}(t)$ approach their equilibrium abundances $\overline{x}^{\ast}(k_2/k_3)$ and $\overline{y}^{\ast}(k_{18}/k_{19})$, respectively; healthy and malignant progenitor cell populations $x_{2}(t)$ and $y_{2}(t)$ approach their equilibrium abundances $\overline{x}^*[\,(k_4+k_5+2k_6)/(-k_7+k_9+k_{10}+k_{14})]\,$ and $\overline{y}^*[\,(k_{20}+k_{21}+2k_{22})/(-k_{23}+k_{25}+k_{26}+k_{30})]\,$, respectively; healthy and malignant differentiated cell populations $x_{3}(t)$ and $y_{3}(t)$ approach their equilibrium abundances $\overline{x}^*[\,((k_4+k_5+2k_6)(k_8+k_9+2k_{10}))/((-k_7+k_9+k_{10}+k_{14})(k_{11}+k_{12}+k_{15}))]\,$ and $\overline{y}^*[\,((k_{20}+k_{21}+2k_{22})(k_{24}+k_{25}+2k_{26}))/((-k_{23}+k_{25}+k_{26}+k_{30})(k_{27}+k_{28}+k_{31}))]\,$, respectively; and healthy and malignant terminally differentiated cell populations $x_{4}(t)$ and $y_{4}(t)$ approach their equilibrium abundances $\overline{x}^*[\,((k_4+k_5+2k_6)(k_8+k_9+2k_{10})(k_{11}+2k_{12}))/((-k_7+k_9+k_{10}+k_{14})(k_{11}+k_{12}+k_{15})k_{16})]\,$ and $\overline{y}^*[\,((k_{20}+k_{21}+2k_{22})(k_{24}+k_{25}+2k_{26})(k_{27}+2k_{28}))/((-k_{23}+k_{25}+k_{26}+k_{30})(k_{27}+k_{28}+k_{31})k_{32})]\,$, respectively. 

Finally, if $R>(b_{1}/b_{2})r,$ the malignant cell populations exhibit the following behaviors: cycling stem cells $y_{0}(t)$ become dominant, approaching their equilibrium abundance $R$ (leukemic homeostatic state); quiescent stem cells $y_{1}(t)$ approach the equilibrium abundance $R(k_{18}/k_{19})$; progenitor cells $y_{2}(t)$ approach the equilibrium abundance $R[\,(k_{20}+k_{21}+2k_{22})/(-k_{23}+k_{25}+k_{26}+k_{30})]\,$; differentiated cells $y_{3}(t)$ approach the equilibrium abundance $R[\,((k_{20}+k_{21}+2k_{22})(k_{24}+k_{25}+2k_{26}))/((-k_{23}+k_{25}+k_{26}+k_{30})(k_{27}+k_{28}+k_{31}))]\,$; and terminally differentiated cells $y_{4}(t)$ approach the equilibrium abundance $R[\,((k_{20}+k_{21}+2k_{22})(k_{24}+k_{25}+2k_{26})(k_{27}+2k_{28}))/((-k_{23}+k_{25}+k_{26}+k_{30})(k_{27}+k_{28}+k_{31})k_{32})]\,$, while all healthy cell populations (cycling stem cells $x_{0}(t)$; quiescent stem cells $x_{1}(t)$; progenitor cells $x_{2}(t)$; differentiated cells $x_{3}(t)$; and terminally differentiated cells $x_{4}(t)$) tend to zero.

A condition such as $R=r$ or $R=(b_{1}/b_{2})r$ is physiologically highly unstable. Small variations in the kinetic parameters can switch the healthy state into the chronic leukemic state and vice versa, if $R=r$, and can switch the chronic state into the accelerated-acute phase and vice versa, if $R=(b_{1}/b_{2})r$. Additionally, from a medical perspective, situations where $R=r$ and $R=(b_{1}/b_{2})r$ are practically undetectable.

In terms of the system's biological growth parameters, the hematological states are characterized as follows:
\vspace{0.2cm}

\begin{equation*}
\begin{array}{ll}
\frac{1}{B}\left( \frac{k_{17}}{k_{21}+k_{22}+k_{29}}-1\right) < \frac{1}{b_{1}}\left( \frac{k_{1}}{k_{5}+k_{6}+k_{13}}-1\right) & \left( \text{healthy state}\right) , \\[8pt]
\frac{1}{b_{1}}\left( \frac{k_{1}}{k_{5}+k_{6}+k_{13}}-1\right) < \frac{1}{B}\left( \frac{k_{17}}{k_{21}+k_{22}+k_{29}}-1\right) < \frac{1}{b_{2}}\left( \frac{k_{1}}{k_{5}+k_{6}+k_{13}}-1\right) \  & \left( \text{%
chronic phase}\right) , \\[8pt]
\frac{1}{b_{2}}\left( \frac{k_{1}}{k_{5}+k_{6}+k_{13}}-1\right) <\frac{1}{B}\left( \frac{k_{17}}{k_{21}+k_{22}+k_{29}}-1\right) \  & \left( \text{accelerated-acute phase}\right) .%
\end{array}%
\end{equation*}
\vspace{0.2cm}

A diagram illustrating the transition from healthy hematopoiesis to chronic and accelerated-acute stages in CML is presented in Figure \ref{fig:2}. Note that the length of the interval $\left[ r,\left( b_{1}/b_{2}\right) r\right]
$\ corresponding to the chronic phase is $\left(b_{1}/b_{2}-1\right) r$\
and depends on the ratio $b_{1}/b_{2}.$ The larger the ratio $\ b_{1}/b_{2}$%
\ is, the larger interval in which a patient's parameter $R$\ can lie in the chronic phase. It's important to recognize that the patient-related parameter $b_{1}/b_{2}$ represents the ratio of contributions of the healthy and malignant (leukemic) stem cells to the reduction of the intrinsic reaction rate $k_{1}$ (or the nonrestrictive growth rate) of the healthy cell population, with the goal of restoring homeostasis. According to our model, values of parameter $R$ close to $r$\ correspond to the early stages of the disease, while values of $R$ close to $\left(b_{1}/b_{2}\right) r$\ indicate advanced stages of the disease moving towards the accelerated-acute phase.

Based on parameter estimation, we highlight that leukemic cycling stem cells exhibit significantly higher intrinsic reaction rate for symmetric self-renewal, asymmetric division, direct differentiation and symmetric differentiation compared to healthy cycling stem cells (twice as high, see the Table \ref{tab.2}). Similarly, leukemic progenitor cells also exhibit higher intrinsic reaction rates for symmetric and asymmetric divisions, direct differentiation, symmetric differentiation and apoptosis (death) compared to healthy progenitor cells (twice as high, see the Table \ref{tab.2}). Our analysis indicated that the progression of leukemic pathology is largely determined by the ratio of death (apoptosis) rates between healthy and leukemic cycling stem cells. Specifically, when the death rate of mutant cycling hematopoietic stem cells (HSCs) is approximately eight times higher than that of healthy cycling HSCs, the disease progresses to a stage where all mutant cell lines are eliminated, assuming mutant cells were initially present. Based on our numerical simulations, we observe that when the death rate of leukemic cycling stem cells ranges between $5 \times k_{13}$ and $7.667 \times k_{13}$, the disease evolved towards a stabilization of both healthy and leukemic cell lines, marking the chronic phase of chronic myeloid leukemia. Finally, the disease progression shifts towards the accelerated-acute phase when the death rate of mutant cycling HSCs falls below $5 \times k_{13}$, a scenario clinically resembling acute/blast phase of chronic myeloid leukemia.

The analysis of our mathematical model and numerical simulations has provided insights and conclusions that are challenging to discern, explain, or quantify solely based on the raw data generated by laboratory studies or clinical observations. Our findings from the mathematical model and simulations lead to the following assertions:

(A) The intrinsic reaction rate of symmetric self-renewal in mutant cycling HSCs serves as a predictive factor for the initiation of the accelerated-acute phase. A heightened intrinsic reaction rate of symmetric self-renewal, compared to healthy cycling stem cells, anticipates an earlier onset of the accelerated-acute state (for related topics, see the work of Shahriyari and Komarova \cite{sk}).

(B) The death rate, along with the intrinsic reaction rate for the direct and symmetric differentiation of leukemic cycling stem cells, serves as a predictive factor for the overall progression of the disease, influencing transitions between different phases of chronic myeloid leukemia.

Consequently, clinical judgment, treatment planning and therapeutic research for leukemic diseases could be influenced or based on the significance of these two factors. Regarding treatment, our mathematical model suggests that a more targeted approach to controlling the proliferation rate and the death rate, along with the intrinsic reaction rate for the direct and symmetric differentiation of mutant cycling HSCs, may be beneficial in achieving a coexisting phase of the two cell populations. This approach should consider patients' symptoms and quality of life rather than adopting an aggressive strategy aimed at eradicating all leukemic cells. Considering recent research and viewpoints, as presented by Gerlinger \textit{et al.} \cite{grhle} and Enriquez-Navas \textit{et al.} \cite{ekdhs}, adopting a strategy of coexistence and collaboration with cancer might be a more intelligent approach than attempting complete eradication.

While the actual clinical experience of the transition of leukemic disease from one phase to another is undeniably complex, it likely involves additional parameters not considered in this study. This omission may occur either because certain factors were challenging/impossible to incorporate into our mathematical model or because they remain unknown to medical researchers. Nevertheless, the application of mathematical models and numerical simulations in the biomedical field can unveil aspects and insights worthy of thorough exploration, in conjunction with clinical practice and ongoing research. These interdisciplinary collaborations serve as a foundation for the essential thought process that will subsequently pinpoint research directions aimed at refining treatments for extensively studied pathologies, such as leukemia. Although the complexities of real-world clinical scenarios may surpass the scope of mathematical modeling, the insights gained contribute valuable perspectives that can inform and enhance medical practice and research.

Although our paper focuses on chronic myeloid leukemia (CML), we acknowledge that the mechanisms of disease progression described are similar to those modeled for other myeloid diseases. The method and conclusions can also be adapted to other types of myeloid leukemias, such as acute myeloid leukemia (AML) (see the work of Cucuianu and Precup \cite{cp}). Furthermore, this mathematical model is versatile and can be applied to other types of leukemia, including lymphoid leukemias, such as acute lymphoid leukemia (ALL) (see the work of Badralexi \textit{et al.} \cite{bhm}).

\subsection*{Code Availability}
\noindent The \textit{Maple Supplementary file}, named \textbf{SupplMaple.mw}, is available at the following link: \url{https://github.com/lorandparajdi/A-Math-Model-of-Clonal-Hematopoiesis-in-CML}

\subsection*{Acknowledgement}
\noindent This work of the first author was supported by the project `The Development of Advanced and Applicative Research Competencies in the Logic of STEAM + Health' (POCU/993/6/13/153310), a project co-financed by the European Social Fund through the Romanian Operational Programme `Human Capital', 2014-2020.

\noindent The authors of this paper, L.G. Parajdi, D. Kegyes and C. Tomuleasa, received funding through a grant from the Academy of Romanian Scientists for the years 2023-2024. D. Kegyes and C. Tomuleasa are funded by a national grant of the Romanian Research Ministry—PNRR 2024-2026 (PNRR/2022/C9/MCID/18, Contract No. 760278/26.03.2024).

\noindent The authors express their gratitude to Professor Dr. Radu Precup for his invaluable assistance in the preparation of this work and thank the referees for their valuable remarks and suggestions, which significantly improved it.

\subsection*{Conflicts of Interest}
\noindent The authors declare that they have no conflict of interest.

\newpage
\section*{Supplementary Notes}
\vspace{0.3cm}

The characteristic equation of the Jacobian matrix $J(\tilde{E}_3)$, evaluated at $\tilde{E}_3$ for the simplified system (3), can be found in the \textit{Maple Supplementary file}: $\textbf{SupplMaple.mw}$, which is accessible via the following link: \\
\url{https://github.com/lorandparajdi/A-Math-Model-of-Clonal-Hematopoiesis-in-CML}

\subsection*{Proofs of $\mu_1,\ \mu_2$ and $\mu_3$ are positive}

Here, $\mu_1,\ \mu_2$ and $\mu_3$ are the coefficients of the characteristic equation $P(\lambda) = \mu_0 \lambda^4 + \mu_1 \lambda^3 + \mu_2 \lambda^2 + \mu_3 \lambda + \mu_4$ corresponding to the Jacobian matrix $J(\tilde{E}_3)$. These coefficients are available in the same \textit{Maple Supplementary file}.
\vspace{0.3cm}

$\blacktriangleright$ To show $\mu_1 >0$, it is equivalent to show that:
\begin{small}
\begin{flalign*}
  \hat{\mu}_1 :=&\ \frac{\mu_1}{(d b_{1}+1) (B D +1)} 
  \\=& -B D^{2} b_{1} b_{2} c_{0}+B D d A_{1} b_{1}^
{2}-B D d A_{1} b_{1} b_{2}+B D d C_{0} b_{1}^{2}
+B D d C_{1} b_{1}^{2}-B D d C_{1} b_{1} b_{2}+B
D d a_{1} b_{1}^{2}-B D d a_{1} b_{1} b_{2} + 
BD d b_{1}^{2} c_{0} \\& + B D d b_{1}^{2} c_{1} - B D d b_{1} b_{2} c_{1}-B d^{2} C_{0} b_{1}^{2}+B D A_{1} b_{1}-B D A_{1} b_{2}+B D C_{0} b_{1}+
B D C_{1} b_{1}-B D C_{1} b_{2} + B D a_
{1} b_{1}-B D a_{1} b_{2}+B D b_{1} c_{1} \\& -B
D b_{2} c_{1} -B d C_{0} b_{1}-D b_{1} b_{2} c_{0}
+d A_{1} b_{1}^{2}-d A_{1} b_{1} b_{2}+d C_{1} b_{1}^{2}-d C_{1} b_
{1} b_{2} + d a_{1} b_{1}^{2}-d a_{1} b_{1} b_{2}+d b_{1}^{2} c_{0}+
db_{1}^{2} c_{1}-d b_{1} b_{2} c_{1}+A_{1} b_{1} \\& -A_{1} b_{2}+C_{1} b_
{1} - C_{1} b_{2}+a_{1} b_{1}-a_{1} b_{2}+b_{1} c_{1}-b_{2} c_{1} &&
\end{flalign*}
\end{small}

\noindent is positive. We can rewrite $\hat{\mu}_1$ and get:
\begin{small}
\begin{flalign*}
  \hat{\mu}_1 =&\ BDb_1c_0(-Db_2+db_1) + BDdA_1b_1(b_1-b_2) + BdC_0b_1^2(D-d) + BDdC_1b_1(b_1-b_2) + BDda_1b_1(b_1-b_2) +  BDdb_1c_1(b_1-b_2) \\& + BDA_1(b_1-b_2) + BC_0b_1(D-d) + BDC_1(b_1-b_2) + BDa_1(b_1-b_2) + BDc_1(b_1-b_2) + b_1c_0(-Db_2+db_1) + dA_1b_1(b_1-b_2) \\& + dC_1b_1(b_1-b_2) + da_1b_1(b_1-b_2) + db_1c_1(b_1-b_2) + A_1(b_1-b_2) + C_1(b_1-b_2) + a_1(b_1-b_2) + c_1(b_1-b_2) \\
  =&\ (-Db_2+db_1)(BDb_1c_0 +b_1c_0) + (D-d)(BdC_0b_1^2 + BC_0b_1) + (b_1-b_2)(BDdA_1b_1 + BDdC_1b_1 + BDda_1b_1 + BDdb_1c_1 + BDA_1 \\& + BDC_1 + BDa_1 + BDc_1 + dA_1b_1 + dC_1b_1 + da_1b_1 + db_1c_1 + A_1 + C_1 + a_1 + c_1)\ . &&
\end{flalign*}
\end{small}

\noindent Thus, $\hat{\mu}_1$ is positive under the assumptions $b_{1}>b_{2}$ and $d<D<(b_{1}/b_{2})d$. This completes the proof of $\mu_{1}>0$ when $b_{1}>b_{2}$ and $d<D<(b_{1}/b_{2})d$.
\vspace{0.3cm}

$\blacktriangleright$ Similarly, to show $\mu_2 >0$, it is equivalent to show that:
\begin{small}
\begin{flalign*}
\hat{\mu}_2 :=&\ \frac{\mu_2}{(d b_{1}+1) (B D +1)}
\\ =&-B D^{2} A_{1} b_{1} b_{2} c_{0}-B D^{2} C_
{0} b_{1} b_{2} c_{0}-B D^{2} C_{1} b_{1} b_{2} c_{0}-B
D^{2} b_{1} b_{2} c_{0} c_{1}+B D d A_{1} a_{1}
b_{1}^{2}-B D d A_{1} a_{1} b_{1} b_{2}+B D d A_
{1} b_{1}^{2} c_{0} +B D d A_{1} b_{1}^{2} c_{1} \\& -B D d A_{1} b_{1} b_{2} c_{1}+B D d C_{0} C_{1} b_{1}^{2}+
B D d C_{0} a_{1} b_{1}^{2}+B D d C_{0} b_{1}^{2}
c_{0} + B D d C_{0} b_{1}^{2} c_{1} + B D d C_{0} b_
{1} b_{2} c_{0} + B D d C_{1} a_{1} b_{1}^{2}-B D d
C_{1} a_{1} b_{1} b_{2} \\& +B D d C_{1} b_{1}^{2} c_{0}+B
D d C_{1} b_{1}^{2} c_{1}-B D d C_{1} b_{1} b_{2}
c_{1}+B D d b_{1}^{2} c_{0} c_{1}-B d^{2} C_{0} C_{1} b_
{1}^{2} - B d^{2} C_{0} a_{1} b_{1}^{2} - B d^{2} C_{0} b_{1}^{2} c_
{0}-B d^{2} C_{0} b_{1}^{2} c_{1}+B D A_{1} a_{1} b_{1} \\& -B
D A_{1} a_{1} b_{2}+B D A_{1} b_{1} c_{1}-BDA_{1} b_{2} c_{1}+B DC_{0} C_{1} b_{1} + BDC_{0} a_{1} b_{1}+B DC_{0} b_{1} c_{1}+B
DC_{1} a_{1} b_{1}-B DC_{1} a_{1} b_{2}+B
DC_{1} b_{1} c_{1}\\& -B DC_{1} b_{2} c_{1}-B d C_
{0} C_{1} b_{1}-B d C_{0} a_{1} b_{1}-B d C_{0} b_{1} c_{1}-D A_{1} b_{1} b_{2} c_{0}-DC_{1} b_{1} b_{2} c_{0}-
Db_{1} b_{2} c_{0} c_{1}+d A_{1} a_{1} b_{1}^{2}-d A_{1}
a_{1} b_{1} b_{2}+d A_{1} b_{1}^{2} c_{0} \\& + d A_{1} b_{1}^{2} c_{1}-d
A_{1} b_{1} b_{2} c_{1}+d C_{1} a_{1} b_{1}^{2}-d C_{1} a_{1} b_{1}
b_{2}+d C_{1} b_{1}^{2} c_{0}+d C_{1} b_{1}^{2} c_{1}-d C_{1} b_{1}
b_{2} c_{1}+d b_{1}^{2} c_{0} c_{1}+A_{1} a_{1} b_{1}-A_{1} a_{1} b_
{2}+A_{1} b_{1} c_{1}\\& -A_{1} b_{2} c_{1}+C_{1} a_{1} b_{1}-C_{1} a_{1}
b_{2}+C_{1} b_{1} c_{1}-C_{1} b_{2} c_{1}
\end{flalign*}
\end{small}

\noindent is positive. 

We can rewrite $\hat{\mu}_2$ and get:
\begin{small}
\begin{flalign*}
\hat{\mu}_2 =&\ BDA_1b_1c_0(-Db_2+db_1) - BDC_0b_1b_2c_0(D-d) + BDC_1b_1c_0(-Db_2+db_1) + BDb_1c_0c_1(-Db_2+db_1) +  BDdA_1a_1b_1(b_1-b_2) \\& +
BDdA_1b_1c_1(b_1-b_2) + BdC_0C_1b_1^2(D-d) + BdC_0a_1b_1^2(D-d) + BdC_0b_1^2c_0(D-d) + BdC_0b_1^2c_1(D-d) + 
BDdC_1a_1b_1(b_1-b_2) \\& + BDdC_1b_1c_1(b_1-b_2) + BDA_1a_1(b_1-b_2) + BDA_1c_1(b_1-b_2) + BC_0C_1b_1(D-d) + BC_0a_1b_1(D-d) + BC_0b_1c_1(D-d) \\& + BDC_1a_1(b_1-b_2) + BDC_1c_1(b_1-b_2) + A_1b_1c_0(-Db_2+db_1) + C_1b_1c_0(-Db_2+db_1) + b_1c_0c_1(-Db_2+db_1) + dA_1a_1b_1(b_1-b_2) \\& + dA_1b_1c_1(b_1-b_2) + dC_1a_1b_1(b_1-b_2) + dC_1b_1c_1(b_1-b_2) + A_1a_1(b_1-b_2) + A_1c_1(b_1-b_2) + C_1a_1(b_1-b_2) + C_1c_1(b_1-b_2)
\end{flalign*}
\end{small}
\vspace{-1.2cm}
\begin{small}
\begin{flalign*}
\ \ \ \ \ \\=&\ (-Db_2+db_1)(BDA_1b_1c_0 + BDC_1b_1c_0 + BDb_1c_0c_1 + A_1b_1c_0 + C_1b_1c_0 + b_1c_0c_1)  +  (D-d)(-BDC_0b_1b_2c_0 + BdC_0C_1b_1^2 + BdC_0a_1b_1^2 \\& + BdC_0b_1^2c_0 
+ BdC_0b_1^2c_1 + BC_0C_1b_1 + BC_0a_1b_1 + BC_0b_1c_1) 
 +  (b_1-b_2)(BDdA_1a_1b_1 + BDdA_1b_1c_1 + BDdC_1a_1b_1 + BDdC_1b_1c_1 \\& +BDA_1a_1 + BDA_1c_1 + BDC_1a_1 + BDC_1c_1 + dA_1a_1b_1 + dA_1b_1c_1 + dC_1a_1b_1 + 
dC_1b_1c_1 + A_1a_1 + A_1c_1 + C_1a_1 + C_1c_1)\\=&\ (-Db_2+db_1)(BDA_1b_1c_0 + BDC_1b_1c_0 + BDb_1c_0c_1 + A_1b_1c_0 + C_1b_1c_0 + b_1c_0c_1)  + (D-d)\Big(BC_0b_1c_0(-Db_2+db_1) + BdC_0C_1b_1^2 \\& + BdC_0a_1b_1^2
+ BdC_0b_1^2c_1 + BC_0C_1b_1 + BC_0a_1b_1 + BC_0b_1c_1\Big) 
 + (b_1-b_2)(BDdA_1a_1b_1 + BDdA_1b_1c_1 + BDdC_1a_1b_1 + BDdC_1b_1c_1 \\& +BDA_1a_1 + BDA_1c_1 + BDC_1a_1 + BDC_1c_1 + dA_1a_1b_1 + dA_1b_1c_1 + dC_1a_1b_1 + 
dC_1b_1c_1 + A_1a_1 + A_1c_1 + C_1a_1 + C_1c_1).
\end{flalign*}
\end{small}

\noindent Thus, $\hat{\mu}_2$ is positive under the assumptions $b_1>b_2$ and $d<D<(b_1/b_2)d$. This completes the proof of $\mu_2>0$ when $b_1>b_2$ and $d<D<(b_1/b_2)d$.
\vspace{0.3cm}

$\blacktriangleright$ To show $\mu_3>0$, similarly, it is equivalent to show that:
\begin{small}
\begin{flalign*}
\hat{\mu}_3 :=&\ \frac{\mu_3}{b_1(d b_{1}+1) (B D +1)}\\
=&-B D^{2} A_{1} b_{2} c_{0} c_{1}-B D^{2} C_
{0} C_{1} b_{2} c_{0}-B D^{2} C_{0} b_{2} c_{0} c_{1}-B
D^{2} C_{1} b_{2} c_{0} c_{1}+B D d A_{1} b_{1}
c_{0} c_{1}+B D d C_{0} C_{1} a_{1} b_{1}
+C_{1} B D d C_{0} b_{1} c_{0} \\&
+B D d C_{0} C_{1} b_{1} c_{1}+d b_
{2} C_{1} B D C_{0} c_{0}+B D d C_{0} b_{1} c_{0}
c_{1}+d b_{2} B D C_{0} c_{0} c_{1}+B D d C_{1}
b_{1} c_{0} c_{1}
-B d^{2} C_{0} C_{1} a_{1} b_{1}-B d^{2} C_{0}
C_{1} b_{1} c_{0} \\& -B d^{2} C_{0} C_{1} b_{1} c_{1}-B d^{2} C_{0}
b_{1} c_{0} c_{1}+C_{1} B D C_{0} a_{1}+C_{1} c_{1} B D C_{0}
-B d C_{0} C_{1} a_{1}-B d C_{0} C_{1} c_{1}-D A_{1} b_{2} c_{0} c_{1}-D C_{1} b_{2} c_{0} c_{1}+d A_{1} b_{1} c_{0} c_{1} \\& + C_{1} c_{1} d b_{1} c_{0}
\end{flalign*}
\end{small}

\noindent is positive. Now we can rearrange the terms in $\hat{\mu}_3$ and get the factored form:
\begin{small}
\begin{flalign*}
\hat{\mu}_3 =&\ BDA_1c_0c_1(-Db_2+db_1) + BDC_0C_1c_0(-Db_2+db_1) + BDC_0c_0c_1(-Db_2+db_1) + BDC_1c_0c_1(-Db_2+db_1) + BdC_0C_1a_1b_1(D-d) \\& + BdC_0C_1b_1c_1(D-d) - BdC_0C_1c_0(-Db_2+db_1) - BdC_0c_0c_1(-Db_2+db_1) + BC_0C_1a_1(D-d) + BC_0C_1c_1(D-d) + A_1c_0c_1 \\& (-Db_2+db_1) + C_1c_0c_1(-Db_2+db_1) = (D-d)( BdC_0C_1a_1b_1 + BdC_0C_1b_1c_1 + BC_0C_1a_1 + BC_0C_1c_1) + (-Db_2 + db_1)(BDA_1c_0c_1 \\& + BDC_0C_1c_0 + BDC_0c_0c_1 + BDC_1c_0c_1 - BdC_0C_1c_0 - BdC_0c_0c_1 + A_1c_0c_1 + C_1c_0c_1) 
\\=&\ (-Db_2 + db_1)\Big(BDA_1c_0c_1 + (D-d)(BC_0C_1c_0 + BC_0c_0c_1)  + BDC_1c_0c_1 + A_1c_0c_1 + C_1c_0c_1\Big) + (D-d)( BdC_0C_1a_1b_1 \\& + BdC_0C_1b_1c_1 + BC_0C_1a_1 + BC_0C_1c_1).
\end{flalign*}
\end{small}

\noindent Therefore, $\hat{\mu}_3 >0$ if we assume $d<D<(b_1/b_2)d$. Hence, we can say $\mu_3 >0$ when $d<D<(b_1/b_2)d$. 

Regarding $\mu_{0}$ and $\mu_{4}$, they are positive if and only if $b_{1} > b_{2}$ and $d<D<(b_1/b_2)d$, as shown in the main text of the article.

\subsection*{Proofs of $\mu_1\mu_2-\mu_0\mu_3$ and $\mu_1\mu_2\mu_3 - \mu_0\mu_3^2 - \mu_1^2\mu_4$ are positive}
\vspace{0.3cm}

We have two remaining conditions to verify in order to establish the asymptotic stability of $\tilde{E}_3$. Firstly, we need to confirm that $\mu_1\mu_2 - \mu_0\mu_3 > 0$. We start by rewriting $\mu_1, \mu_2$ and $\mu_3$ in an abbreviated way to simplify the calculation process. We denote:
\begin{small}
\begin{align*}
    \nu_{11} =&\ BDb_1c_0 + b_1c_0 \\
    \nu_{12} =&\ BdC_0b_1^2 + BC_0b_1 \\
    \nu_{13} =&\ BDdA_1b_1 + BDdC_1b_1 +  BDda_1b_1 +BDdb_1c_1 + BDA_1 + BDC_1 + BDa_1 + BDc_1 + dA_1b_1 +dC_1b_1 + da_1b_1 + db_1c_1 \\& + A_1 + C_1 + a_1 + c_1 \\
    \nu_{21} =&\ BDA_1b_1c_0 + BDC_1b_1c_0 + BDb_1c_0c_1 + A_1b_1c_0 + C_1b_1c_0 + b_1c_0c_1 \\
    \nu_{22} =&\ BC_0b_1c_0(-Db_2+db_1) + BdC_0C_1b_1^2 + BdC_0a_1b_1^2
+ BdC_0b_1^2c_1 + BC_0C_1b_1 + BC_0a_1b_1 + BC_0b_1c_1 \\
    \nu_{23} =&\ BDdA_1a_1b_1 + BDdA_1b_1c_1 + BDdC_1a_1b_1 + BDdC_1b_1c_1 +BDA_1a_1 + BDA_1c_1 + BDC_1a_1 + BDC_1c_1 + dA_1a_1b_1 \ \ \ \ \\& + dA_1b_1c_1 + dC_1a_1b_1 + 
dC_1b_1c_1 + A_1a_1 + A_1c_1 + C_1a_1 + C_1c_1 \\
    \nu_{31} =&\ BDA_1c_0c_1 + (D-d)(BC_0C_1c_0 + BC_0c_0c_1)  + BDC_1c_0c_1 + A_1c_0c_1 + C_1c_0c_1 \\
    \nu_{32} =&\ BdC_0C_1a_1b_1 + BdC_0C_1b_1c_1 + BC_0C_1a_1 + BC_0C_1c_1
\end{align*}
\end{small}
Then we can rewrite $\mu_{0},\ \mu_{1},\ \mu_{2},\ \mu_{3}$ and $\mu_{4}$ as below:
\begin{small}
\begin{align*}
    \mu_0 =&\ (BD+1)^2(b_1-b_2)(db_1+1)^2 \\[6pt]
    \mu_1 =&\ (d b_{1}+1) (B D +1) \hat{\mu}_1\\
    =&\ (d b_{1}+1) (B D +1) \Big[(-Db_2+db_1)(BDb_1c_0 +b_1c_0) + (D-d)(BdC_0b_1^2 + BC_0b_1) + (b_1-b_2)(BDdA_1b_1 + BDdC_1b_1 \\& +  BDda_1b_1 + BDdb_1c_1 + BDA_1 + BDC_1 + BDa_1 + BDc_1 + dA_1b_1 + dC_1b_1 + da_1b_1 + db_1c_1 + A_1 + C_1 + a_1 + c_1)\Big]\\
    =&\ (d b_{1}+1) (B D +1)\Big[(-Db_2+db_1)\nu_{11} + (D-d)\nu_{12} + (b_1-b_2)\nu_{13}\Big]
\end{align*}
\end{small}
\vspace{-0.45cm}
\begin{small}
\begin{align*}
    \mu_2 =&\ (d b_{1}+1) (B D +1)\hat{\mu}_2\\
    =&\ (d b_{1}+1) (B D +1) \Big[(-Db_2+db_1)(BDA_1b_1c_0 + BDC_1b_1c_0 + BDb_1c_0c_1 + A_1b_1c_0 + C_1b_1c_0 + b_1c_0c_1) \\& + (D-d)\big(BC_0b_1c_0(-Db_2+db_1) + BdC_0C_1b_1^2 + BdC_0a_1b_1^2
+ BdC_0b_1^2c_1 + BC_0C_1b_1 + BC_0a_1b_1 + BC_0b_1c_1\big) \\&
 + (b_1-b_2)(BDdA_1a_1b_1 + BDdA_1b_1c_1 + BDdC_1a_1b_1 + BDdC_1b_1c_1 +BDA_1a_1 + BDA_1c_1 + BDC_1a_1 + BDC_1c_1 \ \ \ \\& + dA_1a_1b_1 + dA_1b_1c_1 + dC_1a_1b_1 + 
dC_1b_1c_1 + A_1a_1 + A_1c_1 + C_1a_1 + C_1c_1)\Big]\\
=&\ (d b_{1}+1) (B D +1) \Big[(-Db_2+db_1)\nu_{21} + (D-d)\nu_{22} + (b_1-b_2)\nu_{23}\Big]
\end{align*}
\end{small}
\vspace{-0.45cm}
\begin{small}
\begin{align*}
\mu_3 =&\ b_1(d b_{1}+1) (B D +1)\hat{\mu}_3\\
=&\ b_1(d b_{1}+1) (B D +1)\Big[ (-Db_2 + db_1)\big(BDA_1c_0c_1 + (D-d)(BC_0C_1c_0 + BC_0c_0c_1)  + BDC_1c_0c_1 + A_1c_0c_1 + C_1c_0c_1\big) \\& 
+(D-d)(BdC_0C_1a_1b_1 + BdC_0C_1b_1c_1 + BC_0C_1a_1 + BC_0C_1c_1) \Big]\\
     =&\ b_1(d b_{1}+1) (B D +1)\Big[ (-Db_2 + db_1) \nu_{31} + (D-d) \nu_{32}\Big] \\[6pt]
\mu_4 =&\ b_{1} B C_{0} c_{0} (d b_{1}+1) (D - d) (-D b_{2}+d b_{1}) (BD+1) C_{1} c_{1}.
\end{align*}
\end{small}

\noindent We know that $\ \mu_{0},\ \mu_{1},\ \mu_{2},\ \mu_{3},\ \mu_{4} > 0$ if we assume $b_{1} > b_{2}$ and $d < D < (b_{1}/b_{2})d\ $.
\vspace{0.3cm}

\noindent $\blacktriangledown$ Next, by rewriting the first condition, we can easily show that:
\begin{small}
\begin{align*}
    \mu_1\mu_2 - \mu_0\mu_3 =&\ (db_1+1)^2(BD+1)^2\Big((-Db_2+db_1)\nu_{11} + (D-d)\nu_{12} + (b_1-b_2)\nu_{13} \Big)\\
    &\ \Big((-Db_2+db_1)\nu_{21} + (D-d)\nu_{22} + (b_1-b_2)\nu_{23} \Big) - (BD+1)^3(db_1+1)^3(b_1-b_2)b_1\\
    &\ \Big((-Db_2+db_1)\nu_{31}+(D-d)\nu_{32} \Big)\\
    =&\ (db_1+1)^2(BD+1)^2\Big[(-Db_2+db_1)^2\nu_{11}\nu_{21}+(-Db_2+db_1)(D-d)\nu_{11}\nu_{22} \\& 
    + (-Db_2+db_1)(b_1-b_2)\nu_{11}\nu_{23}+(D-d)(-Db_2+db_1)\nu_{12}\nu_{21} +(D-d)^2 \nu_{12}\\
    &\ \nu_{22} + (b_1-b_2)(D-d)\nu_{12}\nu_{23} + (b_1-b_2)(-Db_2+db_1)\nu_{13}\nu_{21} +(b_1-b_2)(D-d)\\ &\ \nu_{13}\nu_{22}+(b_1-b_2)^2\nu_{13}\nu_{23} - (BD+1)(db_1+1)(b_1-b_2) b_1(-Db_2+db_1)(BDA_1c_0c_1 \\&
    + BDC_1c_0c_1+A_1c_0c_1+C_1c_0c_1)  - (BD+1)(db_1+1)(b_1-b_2)b_1(-Db_2+db_1)(D-d) (BC_0C_1c_0+BC_0c_0c_1) \\& 
    - (BD+1)(db_1+1)(b_1-b_2)b_1(D-d)(BdC_0C_1a_1b_1+BdC_0C_1b_1c_1+BC_0C_1a_1+BC_0C_1c_1) \Big] \\
    =&\ (db_1+1)^2(BD+1)^2\Big[(-Db_2+db_1)^2\nu_{11}\nu_{21}+(-Db_2+db_1)(D-d)\nu_{11}\nu_{22} + (D-d)\\
    &\ (-Db_2+db_1)\nu_{12}\nu_{21} +(D-d)^2 \nu_{12}\nu_{22} +
    (b_1-b_2)(-Db_2+db_1)\nu_{13}\nu_{21} \\& 
    + (b_1-b_2)^2\nu_{13}\nu_{23} + (-Db_2+db_1)(b_1-b_2)(\nu_{11}\nu_{23}-(BD+1)(db_1+1)b_1(BDA_1c_0c_1 \\&
    + BDC_1c_0c_1+A_1c_0c_1+C_1c_0c_1)) + (b_1-b_2)(D-d)(\nu_{13}\nu_{22}-(BD+1)(db_1+1)b_1(-Db_2+db_1)\\&\ (BC_0C_1c_0+BC_0c_0c_1)) + (b_1-b_2)(D-d)(\nu_{12}\nu_{23}-(BD+1)(db_1+1)b_1(BdC_0C_1a_1b_1 \\& 
    + BdC_0C_1b_1c_1 + BC_0C_1a_1+BC_0C_1c_1))\Big] \\
\end{align*}
\end{small}
\vspace{-0.8cm}
\begin{small}
\begin{align*}
\ \ \ \ \ \ \ \ \ \ \ \ \ \ \ \ \ =&\ (db_1+1)^2(BD+1)^2\Big[(-Db_2+db_1)^2\nu_{11}\nu_{21}+(-Db_2+db_1)(D-d)\nu_{11}\nu_{22} + (D-d)\\
    &\ (-Db_2+db_1)\nu_{12}\nu_{21} +(D-d)^2 \nu_{12}\nu_{22} +
    (b_1-b_2)(-Db_2+db_1)\nu_{13}\nu_{21} \\& 
    + (b_1-b_2)^2\nu_{13}\nu_{23} + (-Db_2 + db_1)(b_1-b_2)a_1b_1c_0(BD+1)^2(db_1+1)(A_1+C_1) + (b_1-b_2)(D-d)\\&\ 
    BC_0b_1(db_1+1)(BD+1)(A_1c_0(-Db_2 + db_1) + a_1c_0(-Db_2 + db_1) + dA_1C_1b_1 + dA_1a_1b_1 + dA_1b_1c_1 + dC_1^2b_1 \\& 
    + 2dC_1a_1b_1 + 2dC_1b_1c_1 + da_1^2b_1 + 2da_1b_1c_1 + db_1c_1^2 + A_1C_1 + A_1a_1 + A_1c_1 + C_1^2 + 2C_1a_1 + 2C_1c_1 + a_1^2 \\& + 2a_1c_1 + c_1^2) + (b_1-b_2)(D-d)BA_1C_0b_1(db_1+1)^2(a_1+c_1)(BD+1)\Big]\\
    =&\ (db_1+1)^2(BD+1)^2\Big[(-Db_2+db_1)^2\nu_{11}\nu_{21}+(-Db_2+db_1)(D-d)\nu_{11}\nu_{22} + (D-d)\\
    &\ (-Db_2+db_1)\nu_{12}\nu_{21} +(D-d)^2 \nu_{12}\nu_{22} +
    (b_1-b_2)(-Db_2+db_1)\nu_{13}\nu_{21} \\& 
    + (b_1-b_2)^2\nu_{13}\nu_{23} + (-Db_2 + db_1)(b_1-b_2)a_1b_1c_0(BD+1)^2(db_1+1)(A_1+C_1) + (b_1-b_2)(D-d)\\&\ 
    BC_0b_1(db_1+1)(BD+1)\Big(c_0(A_1+a_1)(-Db_2+db_1) + (db_1+1)(C_1+a_1+c_1)(A_1+C_1+a_1+c_1)\Big) \\& 
    + (b_1-b_2)(D-d)BA_1C_0b_1(db_1+1)^2(a_1+c_1)(BD+1)\Big] > 0 \ \ \text{if and only if}\ \ b_{1}>b_{2}\ \ \text{and}\ \ d<D<\frac{b_1}{b_2}d.
\end{align*}
\end{small}

\noindent Next, we want to show $\mu_1 \mu_2 \mu_3 - \mu_0\mu_3^2 - \mu_1^2\mu_4 >0$.
We notice $\mu_1 \mu_2 \mu_3 - \mu_0\mu_3^2 = (\mu_1 \mu_2 - \mu_0\mu_3)\mu_3 > 0$ since we proved $\mu_1 \mu_2 - \mu_0\mu_3 > 0$. We can show $\mu_1 \mu_2 \mu_3 - \mu_0\mu_3^2 - \mu_1^2\mu_4 >0$ by showing $ (\mu_1 \mu_2 - \mu_0\mu_3)\mu_3 > \mu_1^2\mu_4$.
That is to show:
\vspace{-0.3cm}

\begin{small}
\begin{align*}
    &(db_1+1)^3(BD+1)^3b_1 \Big( (-Db_2 + db_1) \nu_{31} + (D-d) \nu_{32}\Big) \\&
    \Big((-Db_2+db_1)^2\nu_{11}\nu_{21}+(-Db_2+db_1)(D-d)\nu_{11}\nu_{22} + (D-d)\\&
    (-Db_2+db_1)\nu_{12}\nu_{21} +(D-d)^2 \nu_{12}\nu_{22} +
    (b_1-b_2)(-Db_2+db_1)\nu_{13}\nu_{21} \\& + (b_1-b_2)^2\nu_{13}\nu_{23} + (-Db_2 + db_1)(b_1-b_2)a_1b_1c_0(BD+1)^2(db_1+1)(A_1+C_1) + (b_1-b_2)(D-d) \\& 
    BC_0b_1 (db_1+1)(BD+1)(A_1c_0(-Db_2 + db_1) + a_1c_0(-Db_2 + db_1) + dA_1C_1b_1 + dA_1a_1b_1 + dA_1b_1c_1 + dC_1^2b_1 \\& + 2dC_1a_1b_1 
    + 2dC_1b_1c_1 + da_1^2b_1 + 2da_1b_1c_1 + db_1c_1^2 + A_1C_1 + A_1a_1 + A_1c_1 + C_1^2 + 2C_1a_1 + 2C_1c_1 + a_1^2 + 2a_1c_1 + c_1^2) \\& + (b_1-b_2)(D-d)BA_1C_0b_1(db_1+1)^2(a_1+c_1)(BD+1)\Big) >\\& 
    (d b_{1}+1)^3 (B D +1)^3 b_1\Big((-Db_2+db_1)\nu_{11} + (D-d)\nu_{12} + (b_1-b_2)\nu_{13}\Big)^2 B C_{0} c_{0} (D-d)(-D b_{2}+d b_{1}) C_{1} c_{1}. 
\end{align*}
\end{small}

\noindent Here we denote: 
\begin{small}
\begin{align*}
    \zeta =&\ A_1c_0(-Db_2 + db_1) + a_1c_0(-Db_2 + db_1) + dA_1C_1b_1 + dA_1a_1b_1 + dA_1b_1c_1 + dC_1^2b_1 + 2dC_1a_1b_1 + 2dC_1b_1c_1 \\& 
    + da_1^2b_1 + 2da_1b_1c_1 + db_1c_1^2 + A_1C_1 + A_1a_1 + A_1c_1 + C_1^2 + 2C_1a_1 + 2C_1c_1 + a_1^2 + 2a_1c_1 + c_1^2
\end{align*}
\end{small}

\noindent so the above inequality is equivalent to:
\begin{small}
\begin{align*}
    &(db_1+1)^3(BD+1)^3b_1\Big[(-Db_2+db_1)^3\nu_{11}\nu_{21}\nu_{31} + (-Db_2+db_1)^2(D-d)(\nu_{11}\nu_{22}\nu_{31} \\& 
    + \nu_{12}\nu_{21}\nu_{31})+(-Db_2+db_1)(D-d)^2\nu_{12}\nu_{22}\nu_{31}+ (-Db_2+db_1)^2(b_1-b_2)\nu_{13}\\&\nu_{21}\nu_{31}+(-Db_2+db_1)(b_1-b_2)^2\nu_{13}\nu_{23}\nu_{31} + (-Db_2+db_1)^2(b_1-b_2)a_1b_1c_0(BD+1)^2(db_1+1)\\&(A_1+C_1)\nu_{31}+(-Db_2+db_1)(D-d)(b_1-b_2)BC_0b_1(db_1+1)(BD+1)\nu_{31}\zeta+
    (-Db_2+db_1)(D-d)(b_1-b_2)\\&BA_1C_0b_1(db_1+1)^2(a_1+c_1)(BD+1)\nu_{31}+(-Db_2+db_1)^2(D-d)\nu_{11}\nu_{21}\nu_{32} +  
    (-Db_2+db_1)(D-d)^2\\&(\nu_{11}\nu_{22}\nu_{32}+\nu_{12}\nu_{21}\nu_{32}) + (D-d)^3\nu_{12}\nu_{22}\nu_{32} + (-Db_2+db_1)(b_1-b_2)\\&(D-d)\nu_{13}\nu_{21}\nu_{32} + (D-d)(b_1-b_2)^2\nu_{13}\nu_{23}\nu_{32} +  (-Db_2+db_1)(b_1-b_2)(D-d)a_1b_1c_0\\&(BD+1)^2(db_1+1)(A_1+C_1)\nu_{32}+(b_1-b_2)(D-d)^2BC_0b_1(db_1+1)(BD+1)\nu_{32}\zeta + (b_1-b_2)(D-d)^2\\&BA_1C_0b_1(db_1+1)^2(a_1+c_1)(BD+1)\nu_{32}\Big] >\\&  (db_1+1)^3(BD+1)^3b_1\Big[(-Db_2+db_1)^3(D-d)BC_0C_1c_0c_1\nu_{11}^2 + 2(-Db_2+db_1)^2(D-d)^2BC_0C_1c_0c_1\nu_{11}\nu_{12} \\& + 2(-Db_2+db_1)^2(D-d)(b_1-b_2)BC_0C_1c_0c_1\nu_{11}\nu_{13} + (-Db_2+db_1)(D-d)^3BC_0C_1c_0c_1\nu_{12}^2 + 2(-Db_2+db_1)\\&(D-d)^2(b_1-b_2)BC_0C_1c_0c_1\nu_{12}\nu_{13} + (-Db_2+db_1)(D-d)(b_1-b_2)^2BC_0C_1c_0c_1\nu_{13}^2\Big].
\end{align*}
\end{small}

\noindent By comparing the terms between the LHS and RHS, we can easily show the following five inequalities.
\vspace{0.3cm}

$\blacktriangleright$ First, we can show:
\begin{small}
\begin{align*}
           (-Db_2+db_1)^3\nu_{11}\nu_{21}\nu_{31} - (-Db_2+db_1)^3(D-d)BC_0C_1c_0c_1\nu_{11}^2 > 0
\end{align*}
\end{small}

\noindent if and only if
\begin{small}
\begin{align*}
& (-D b_2 +db_1 )^{3} b_1
^{2} c_0^{3} (B D +1)^{2} \Big(A_1^{2} B D c_1 +(D-d)(A_1 B C_0 C_1 +A_1 B C_0 c_1 +  B C_0 C_1^{2}  +B C_0 C_1 c_1  +B C_0 c_1^{2}) +2 A_1 B C_1 D c_1 \\& 
+ A_1 B D c_1^{2}+ B C_1^{2} D c_1
+B C_1 D c_1^{2}+A_1^{2} c_1 +2 A_1 C_1 c_1 +A_1 c_1^{2}+C_1^{2} c_1 +C_1 c_1^{2}\Big) > 0
\end{align*}
\end{small}

\noindent which is true if we assume $d<D<(b_1/b_2)d$.
\vspace{0.3cm}

$\blacktriangleright$ Second, using the same strategy, we can show:
\begin{small}
\begin{align*}
        (-Db_2+db_1)^2(D-d)\nu_{12}\nu_{21}\nu_{31} - 2(-Db_2+db_1)^2(D-d)^2BC_0C_1c_0c_1\nu_{11}\nu_{12} > 0
\end{align*}
\end{small}

\noindent if and only if
\begin{small}
\begin{align*}
    &(D b_2 -b_1 d )^{2}(D -d ) B C_0 b_1^{2} (b_1 d +1) c_0^{2} (B D +1) \Big(A_1^{2} B D c_1 +(D-d)(A_1 B C_0
C_1 +A_1 B C_0 c_1 +B C_0 C_1^{2} +B C_0  c_1^{2}) \\& 
+ 2 A_1 B C_1 D c_1 +A_1 B D c_1^{2}
 +B C_1^{2} D c_1 +B C_1 D c_1^{2}+A_1^{2} c_1 +2
A_1 C_1 c_1 +A_1 c_1^{2}+
C_1^{2} c_1 +C_1 c_1^{2}\Big) > 0
\end{align*}
\end{small}

\noindent which is true if we assume $D>d$.
\vspace{0.3cm}

$\blacktriangleright$ Third, similarly, we can show:
\begin{small}
\begin{align*}
        (-Db_2+db_1)^2(b_1-b_2)\nu_{13}\nu_{21}\nu_{31} - 2(-Db_2+db_1)^2(b_1-b_2)(D-d)BC_0C_1c_0c_1\nu_{11}\nu_{13} > 0
\end{align*}
\end{small}

\noindent if and only if
\begin{small}
\begin{align*}
&(D b_2 -b_1 d )^{2} 
(b_1 -b_2 ) (A_1+C_1 +
a_1+c_1 ) (b_1 d +1) 
(B D +1)^{2} b_1 c_0^{2}
\Big(A_1^{2} B D c_1 +(D-d)(A_1B C_0
C_1 + A_1B C_0c_1+B C_0 C_1^{2} \\& + B C_0 c_1^{2} ) +
 2 A_1B C_1 D
c_1 +A_1B D c_1^{2} 
+B C_1^{2} D c_1 +B C_1 D c_1^{2}+A_1^{2} c_1 +2
A_1C_1 c_1 +A_1c_1^{2}+
C_1^{2} c_1 +C_1 c_1^{2}\Big) >0
\end{align*}
\end{small}

\noindent which is true if we assume $b_1>b_2$ and $D>d$.
\vspace{0.3cm}

$\blacktriangleright$ We can show the fourth inequality as written below holds since:
\begin{small}
\begin{align*}
        (D-d)^3\nu_{12}\nu_{22}\nu_{32} - (-Db_2+db_1)(D-d)^3BC_0C_1c_0c_1\nu_{12}^2>0
\end{align*}
\end{small}

\noindent if and only if
\begin{small}
\begin{align*}
    &(D -d )^{3} B^{3} C_0^{3} b_1^
{2} (b_1 d +1)^{2} C_1 \Big(a_1c_0(-Db_2+db_1) + d C_1
a_1 b_1 + d C_1
b_1 c_1 + d a_1^{2} b_1
+ 2d a_1 b_1 c_1 + d b_1 c_1^{2} \\& + C_1 a_1 + C_1 c_1 +
a_1^{2} + 2 a_1 c_1 + c_1^{2}\Big) > 0
\end{align*}
\end{small}

\noindent which is true if we assume $d<D<(b_1/b_2)d$.
\vspace{0.3cm}

$\blacktriangleright$ Last but not least, we can show the fifth inequality as written below also holds since:
\begin{small}
\begin{align*}
        (-Db_2+db_1)(b_1-b_2)^2\nu_{13}\nu_{23}\nu_{31} - (-Db_2+db_1)(b_1-b_2)^2(D-d)BC_0C_1c_0c_1\nu_{13}^2 > 0
\end{align*}
\end{small}

\noindent if and only if
\begin{small}
\begin{align*}
&(-D b_2 +b_1 d ) (b_1 -b_2 )^{2} (A_1 +C_1 +a_1 +c_1 ) (b_1 d +1)^{2} (B
D +1)^{2} c_0 \Big(A_1^{2} B D a_1 c_1 +A_1^{2} B D c_1^{2}+(D-d)(A_1 B C_0 C_1a_1 \\& + A_1 B C_0 a_1 c_1 +A_1 B
C_0 c_1^{2} + B
C_0 C_1^{2} a_1)+2 A_1 B C_1 D a_1 c_1 +2 A_1 B C_1 D c_1^{2} +B C_1^{2} D
a_1 c_1 +B C_1^{2} D c_1^{2}+A_1^{2} a_1 c_1 +A_1^{2}
c_1^{2} \\& + 2 A_1 C_1 a_1 c_1 +2
A_1 C_1 c_1^{2}+C_1^{2} a_1
c_1 +C_1^{2} c_1^{2}\Big) >0
\end{align*}
\end{small}

\noindent which is true if we assume $b_1 > b_2$ and $d<D<(b_1/b_2)d$.
\vspace{0.3cm}

$\blacktriangleright$ Now, we only have one term left on the RHS and by combining it with one term from LHS,  we can show that:
\begin{small}
\begin{align*}
&(-Db_2+db_1)(D-d)(b_1-b_2)BC_0b_1(db_1+1)(BD+1)\nu_{31}\zeta\ \\& 
- 2(-Db_2+db_1)(b_1-b_2)(D-d)^2BC_0C_1c_0c_1\nu_{12}\nu_{13} > 0.
\end{align*}
\end{small}
We factor $\nu_{12},\ \nu_{13},\ \nu_{31}$ and $\zeta$:
\begin{small}
\begin{align*}
 \nu_{12} =&\ BdC_0b_1^2 + BC_0b_1
 = BC_0b_1(db_1+1) \\
    \nu_{13} =&\ BDdA_1b_1 + BDdC_1b_1 +  BDda_1b_1 +BDdb_1c_1 + BDA_1 + BDC_1 + BDa_1 + BDc_1 + dA_1b_1 +dC_1b_1 \\& + da_1b_1 + db_1c_1 + A_1 + C_1 + a_1 + c_1 = (BD+1)(db_1+1)(A_1+C_1+a_1+c_1) \\
    \nu_{31} =&\ BDA_1c_0c_1 + (D-d)(BC_0C_1c_0 + BC_0c_0c_1)  + BDC_1c_0c_1 + A_1c_0c_1 + C_1c_0c_1 \\
    =&\ (D -d)BC_0c_0(C_1+c_1)+c_0c_1(A_1+C_1)(BD+1) \\
    \zeta =&\ A_1c_0(-Db_2 + db_1) + a_1c_0(-Db_2 + db_1) + dA_1C_1b_1 + dA_1a_1b_1 + dA_1b_1c_1 + dC_1^2b_1 + 2dC_1a_1b_1 + 2dC_1b_1c_1 \\& + 
    da_1^2b_1 + 2da_1b_1c_1 + db_1c_1^2 + A_1C_1 + A_1a_1 + A_1c_1 + C_1^2 + 2C_1a_1 + 2C_1c_1 + a_1^2 + 2a_1c_1 + c_1^2 \\
    =&\ c_0(A_1+a_1)(-Db_2+db_1) + (db_1+1)(C_1+a_1+c_1)(A_1+C_1+a_1+c_1)
\end{align*}
\end{small}

\noindent and we rewrite the above inequality:
\begin{small}
\begin{align*}
    &(-Db_2+db_1)(D-d)(b_1-b_2)BC_0b_1(db_1+1)(BD+1)\Big((D-d)BC_0c_0(C_1+c_1)+c_0c_1(A_1+C_1)(BD+1)\Big)\\&\Big(c_0(A_1+a_1)(-Db_2+db_1) + (db_1+1)(C_1+a_1+c_1)(A_1+C_1+a_1+c_1)\Big)\ \\& 
    - 2(-Db_2+db_1)(b_1-b_2)(D-d)^2BC_0C_1c_0c_1\Big(BC_0b_1(db_1+1)\Big)\Big((BD+1)(db_1+1)(A_1+C_1+a_1+c_1)\Big) > 0
\end{align*}
\end{small}

\noindent if and only if
\begin{small}
\begin{align*}
    &(-Db_2+db_1)(D-d)(b_1-b_2)BC_0b_1(db_1+1)(BD+1)\Big((D-d)BC_0c_0^2(C_1+c_1)(A_1+a_1)(-Db_2+db_1) \\& 
    + (D-d)BC_0c_0(C_1+c_1)(db_1+1)(C_1+a_1+c_1)(A_1+C_1+a_1+c_1) + c_0^2c_1(A_1+C_1)(BD+1)(A_1+a_1)(-Db_2+db_1) \\&
    + c_0c_1(A_1+C_1)(BD+1)(db_1+1)(C_1+a_1+c_1)(A_1+C_1+a_1+c_1)\Big) - 2(-Db_2+db_1)(b_1-b_2)(D-d)^2BC_0C_1c_0c_1\\ 
    &\Big(BC_0b_1(db_1+1)\Big)\Big((BD+1)(db_1+1)(A_1+C_1+a_1+c_1)\Big) > 0  
\end{align*}
\end{small}

\noindent if and only if
\begin{small}
\begin{align*}
    &(-Db_2+db_1)(D-d)(b_1-b_2)BC_0b_1(db_1+1)(BD+1)\Big((D-d)BC_0c_0(C_1+c_1)(db_1+1)(C_1+a_1+c_1)(A_1+C_1+a_1+c_1)\Big) \\&
    - 2(-Db_2+db_1)(b_1-b_2)(D-d)^2BC_0C_1c_0c_1\Big(BC_0b_1(db_1+1)\Big)\Big((BD+1)(db_1+1)(A_1+C_1+a_1+c_1)\Big) \\& 
    + (-Db_{2}+db_1)(D-d)(b_1-b_2)BC_0b_1(db_1+1)(BD+1)\Big((D-d)BC_0c_0^2(C_1+c_1)(A_1+a_1)(-Db_2+db_1) \\& 
    + c_0^2c_1(A_1+C_1)(BD+1)(A_1+a_1)(-Db_2+db_1) + c_0c_1(A_1+C_1)(BD+1)(db_1+1)(C_1+a_1+c_1)(A_1+C_1+a_1+c_1)\Big) > 0.
\end{align*}
\end{small}

\noindent Next, we will demonstrate that this inequality holds. To do this, we will show:
\begin{small}
\begin{align*}
    &(-Db_2+db_1)(D-d)(b_1-b_2)BC_0b_1(db_1+1)(BD+1)\Big((D-d)BC_0c_0(C_1+c_1)(db_1+1)(C_1+a_1+c_1)(A_1+C_1+a_1+c_1) \Big) \\& 
    - 2(-Db_2+db_1)(b_1-b_2)(D-d)^2BC_0C_1c_0c_1\Big(BC_0b_1(db_1+1)\Big)\Big((BD+1)(db_1+1)(A_1+C_1+a_1+c_1)\Big) > 0
\end{align*}
\end{small}

\noindent if and only if
\begin{small}
\begin{align*}
      &(-Db_2+db_1)(D-d)^2(b_1-b_2)B^2C_0^2b_1(db_1+1)^2(BD+1)c_0(A_1+C_1+a_1+c_1)(C_{1}^2+C_{1}a_{1}+a_{1}c_{1}+c_{1}^2) > 0
\end{align*}
\end{small}

\noindent which is true if we assume $b_{1}>b_{2}$ and $d<D<(b_{1}/b_{2})d$. Therefore, we have:
\begin{small}
\begin{align*}
    &(-Db_2+db_1)(D-d)^2(b_1-b_2)B^2C_0^2b_1c_0(db_1+1)^2(BD+1)(A_1+C_1+a_1+c_1)(C_{1}^2+C_{1}a_{1}+a_{1}c_{1}+c_{1}^2) \\& 
    + (-Db_{2}+db_1)(D-d)(b_1-b_2)BC_0b_1(db_1+1)(BD+1)\Big((D-d)BC_0c_0^2(C_1+c_1)(A_1+a_1)(-Db_2+db_1) \\& 
    + c_0^2c_1(A_1+C_1)(BD+1)(A_1+a_1)(-Db_2+db_1) + c_0c_1(A_1+C_1)(BD+1)(db_1+1)(C_{1}+a_{1}+c_{1})(A_1+C_1+a_1+c_1)\Big) > 0
\end{align*}
\end{small}

\noindent which is true if we assume $\ b_1 > b_2\ $ and $\ d<D<(b_1/b_2)d\ $.
\vspace{0.3cm}

\noindent $\blacktriangledown$ Thus, the second condition can be expressed as:
\vspace{-0.4cm}

\begin{small}
\begin{flalign*}
&\mu_1\mu_2\mu_3 - \mu_1^2\mu_4 - \mu_0\mu_3^2 = (db_1+1)^3(BD+1)^3b_1\Big[(-D b_2 +db_1 )^{3} b_1^{2} c_0^{3} (B D +1)^{2} \Big(A_1^{2} B D c_1 +(D-d)(A_1 B C_0 C_1 +A_1 B C_0 c_1 \\
&+ B C_0 C_1^{2}  +B C_0 C_1 c_1 + B C_0 c_1^{2}) + 2 A_1 B C_1 D c_1 + A_1 B D c_1^{2}+ B C_1^{2} D c_1 + B C_1 D c_1^{2}+A_1^{2} c_1 +2 A_1 C_1 c_1 +A_1 c_1^{2}+C_1^{2} c_1 +C_1 c_1^{2}\Big) \\
&+ (-Db_2+db_1)^2(D-d)\nu_{11}\nu_{22}\nu_{31} + (D b_2 -b_1 d )^{2}(D -d ) B C_0 b_1^{2} (b_1 d +1) c_0^{2} (B D +1) \Big(A_1^{2} B D c_1 +(D-d)\\
&(A_1 B C_0 C_1 + A_1 B C_0 c_1 + B C_0 C_1^{2} +B C_0  c_1^{2}) + 2 A_1 B C_1 D
c_1 +A_1 B D c_1^{2}
 +B C_1^{2} D c_1 +B C_1 D c_1^{2}+A_1^{2} c_1 +2
A_1 C_1 c_1 +A_1 c_1^{2}+
C_1^{2} c_1 \\
&+ C_1 c_1^{2}\Big) + (-Db_2+db_1)(D-d)^2\nu_{12}\nu_{22}\nu_{31} + (D b_2 -b_1 d )^{2} 
(b_1 -b_2 ) (A_1+C_1 + a_1+c_1 ) (b_1 d +1) (B D +1)^{2} \\
&b_1 c_0^{2}\Big(A_1^{2} B D c_1 +(D-d)(A_1B C_0C_1 + A_1B C_0c_1 + B C_0 C_1^{2}+B C_0 c_1^{2} ) + 2 A_1B C_1 D
c_1 +A_1B D c_1^{2} 
+B C_1^{2} D c_1 +B C_1 D c_1^{2}+A_1^{2} c_1 \\
& + 2A_1C_1 c_1 +A_1c_1^{2}+
C_1^{2} c_1 +C_1 c_1^{2}\Big) + (D -d )^{3} B^{3} C_0^{3} b_1^{2}(b_1 d +1)^{2} C_1 \Big(a_1c_0(-Db_2+db_1) + d C_1
a_1 b_1 +d C_1
b_1 c_1 + d a_1^{2} b_1\\ 
&+ 2d a_1 b_1 c_1 + d b_1 
c_1^{2}+C_1 a_1 +C_1 c_1 +
a_1^{2}+2 a_1 c_1 +c_1^{2}\Big) + (-D b_2 +b_1 d ) (b_1 -b_2 )^{2} (A_1 +C_1 +a_1 +c_1 ) (b_1 d +1)^{2}(B D +1)^{2}\\
&c_0 \Big(A_1^{2} B D a_1 c_1 +A_1^{2} B D c_1^{2}+(D-d)(A_1 B C_0 C_1a_1 + A_1 B C_0 a_1 c_1 + A_1 B C_0 c_1^{2} + BC_0 C_1^{2} a_1)+2 A_1 B C_1 D a_1 c_1 + 2 A_1 B C_1 D c_1^{2} \\
&+B C_1^{2} D a_1 c_1 +B C_1^{2} D c_1^{2}+A_1^{2} a_1 c_1 +A_1^{2}
c_1^{2}+2 A_1 C_1 a_1 c_1 +2A_1 C_1 c_1^{2} + C_1^{2} a_1c_1 +C_1^{2} c_1^{2}\Big) + (-Db_2+db_1)^2(b_1-b_2)a_1b_1c_0\\
&(BD+1)^2(db_1+1)(A_1+C_1)\nu_{31} + (-Db_2+db_1)(D-d)^2(b_1-b_2)B^2C_0^2b_1c_0(db_1+1)^2(BD+1)(A_1+C_1+a_1+c_1)\\
&(C_1^2+C_1a_1+a_1c_1+c_1^2) + (-Db_{2}+db_1)(D-d)(b_1-b_2)BC_0b_1(db_1+1)(BD+1)\Big((D-d)BC_0c_0^2(C_1+c_1)(A_1+a_1)\\
&(-Db_2+db_1) + c_0^2c_1(A_1+C_1)(BD+1)(A_1+a_1)(-Db_2+db_1) +  c_0c_1(A_1+C_1)(BD+1)(db_1+1)(C_1+a_1+c_1)(A_1+C_1\\
&+a_1+c_1)\Big) + (-Db_2+db_1)(D-d)(b_1-b_2) BA_1C_0b_1(db_1+1)^2(a_1+c_1)(BD+1)\nu_{31} + (-Db_2+db_1)^2(D-d)\nu_{11}\nu_{21}\nu_{32}\\
&+ (-Db_2+db_1)(D-d)^2(\nu_{11}\nu_{22}\nu_{32} + \nu_{12}\nu_{21}\nu_{32}) + (-Db_2+db_1)(b_1-b_2)(D-d)\nu_{13}\nu_{21}\nu_{32} + (D-d)(b_1-b_2)^2\nu_{13}\nu_{23}\nu_{32}\\
&+ (-Db_2+db_1)(D-d)(b_1-b_2)a_1b_1c_0(BD+1)^2(db_1+1)(A_1+C_1)\nu_{32} + (b_1-b_2)(D-d)^2BC_0b_1(db_1+1)(BD+1)\nu_{32}\\
&\Big(c_0(A_1+a_1)(-Db_2+db_1) + (db_1+1)(C_1+a_1+c_1)(A_1+C_1+a_1+c_1)\Big) + (b_1-b_2)(D-d)^2BA_1C_0b_1(db_1+1)^2(a_1+c_1)\\
&(BD+1)\nu_{32}\Big] > 0\ \ \text{if and only if}\ \ b_1>b_2\ \ \text{and}\ \ d<D<\frac{b_1}{b_2}d. &&
\end{flalign*}
\end{small}

\subsection*{Considering the initial mathematical model (2)}
\vspace{0.3cm}

It's important to note that the same conclusions hold (based on the Routh-Hurwitz criterion for the 4th-order polynomials), i.e., the same inequalities hold, even when we work with the initial mathematical model $(2)$. We have:
\vspace{-0.3cm}

\begin{small}
\begin{flalign*}
\mu_0 =&\ (BR+1)^2(b_1-b_2)(rb_1+1)^2 \\[6pt]
\mu_1 =&\ (r b_{1}+1) (B R +1) \Big[\ (-Rb_2+rb_1)\Big(BRb_1(k_5+k_6+k_{13}) +b_1(k_5+k_6+k_{13})\Big) + (R-r)\Big(Br(k_{21}+k_{22}+k_{29})b_1^2 + B(k_{21}+k_{22} \\& +k_{29}) b_1\Big) + (b_1-b_2)\Big(BRrk_{18}b_1 + BRrk_{19}b_1 + BRrk_{2}b_1 + BRrb_1k_{3} + BRk_{18} + BRk_{19} + BRk_{2} + BRk_{3} + rk_{18}b_1 + rk_{19}b_1 \\& + rk_{2}b_1 + rb_1k_{3} + k_{18} + k_{19} + k_{2} + k_{3}\Big)\Big] 
  = (r b_{1}+1) (B R +1) \Big[\ (-Rb_2+rb_1)\overline{\nu}_{11} + (R-r)\overline{\nu}_{12} + (b_1-b_2)\overline{\nu}_{13}\Big]\  &&
  \end{flalign*}
\end{small}
\vspace{-0.5cm}
\begin{small}
\begin{flalign*}
\mu_2 =&\ (r b_{1}+1) (B R +1)\Big[\ (-Rb_2+rb_1)\Big(BRk_{18}b_1(k_5+k_6+k_{13}) + BRk_{19}b_1(k_5+k_6+k_{13}) + BRb_1(k_5+k_6+k_{13})k_{3} + k_{18}b_1(k_5+k_6 \\& +k_{13}) + k_{19}b_1(k_{5}+k_{6}+k_{13}) + b_1(k_5+k_6+k_{13})k_{3}\Big) +  (R-r)\Big(B(k_{21}+k_{22}+k_{29})b_1(k_{5}+k_{6}+k_{13}) (-Rb_2+rb_1) + Br(k_{21}+k_{22} \\& +k_{29})k_{19}b_1^2 + Br(k_{21}+k_{22}+k_{29})k_{2}b_1^2
+ Br(k_{21}+k_{22}+k_{29})b_1^2k_{3} + B(k_{21}+k_{22}+k_{29})k_{19}b_1 + B(k_{21}+k_{22}+k_{29})k_{2}b_1 + B(k_{21} \\& +k_{22}+k_{29})b_1k_{3}\Big) + (b_1-b_2)\Big(BRrk_{18}k_{2}b_1 + BRrk_{18}b_1k_{3} + BRrk_{19}k_{2}b_1 + BRrk_{19}b_1k_{3} + BRk_{18}k_{2} + BRk_{18}k_{3} + BRk_{19}k_{2} \\& + BRk_{19}k_{3} +  rk_{18}k_{2}b_1 + rk_{18}b_1k_{3} + rk_{19}k_{2}b_1 + 
rk_{19}b_1k_{3} + k_{18}k_{2} + k_{18}k_{3} + k_{19}k_{2} + k_{19}k_{3}\Big)\Big]\ \\
=&\ (r b_{1}+1) (B R +1)\Big[\ (-Rb_2+rb_1)\overline{\nu}_{21} + (R-r)\overline{\nu}_{22} + (b_1-b_2)\overline{\nu}_{23}\Big]\ &&
  \end{flalign*}
\end{small}
\vspace{-0.5cm}
\begin{small}
\begin{flalign*}
\mu_3 =&\ b_{1} (r b_{1}+1) (B R +1)\Big[\ (-Rb_2 + rb_1)\Big(BRk_{18}(k_5+k_6+k_{13})k_{3} + (R-r)(B(k_{21}+k_{22}+k_{29})k_{19}(k_{5}+k_{6}+k_{13}) \\
&+ B(k_{21}+k_{22}+k_{29})(k_5+k_6+k_{13})k_{3})  + BRk_{19}(k_5+k_6+k_{13})k_{3} + k_{18}(k_5+k_6+k_{13})k_{3} + k_{19}(k_5+k_6+k_{13})k_{3}\Big) \\
&+ (R-r)\Big( Br(k_{21}+k_{22}+k_{29})k_{19}k_{2}b_1 + Br(k_{21}+k_{22}+k_{29})k_{19}b_1k_{3} + B(k_{21}+k_{22}+k_{29})k_{19}k_{2} + B(k_{21}+k_{22}+k_{29})k_{19}k_{3}\Big)\Big]\ \\ 
 =&\ b_{1} (r b_{1}+1) (B R +1)\Big[\ (-Rb_2 + rb_1)\overline{\nu}_{31} + (R-r)\overline{\nu}_{32}\Big] \\[6pt] 
\mu_4 =&\ b_{1} B (k_{21}+k_{22}+k_{29}) (k_{5}+k_{6}+k_{13}) (r b_{1}+1) (R - r) (-R b_{2}+r b_{1}) (B R + 1) k_{19} k_{3}. && 
   \end{flalign*}
\end{small}

\noindent Thus, we know that $\ \mu_{0},\ \mu_{1},\ \mu_{2},\ \mu_{3},\ \mu_{4} > 0\ $ if and only if $\ b_{1} > b_{2}\ $ and $\ r < R < (b_{1}/b_{2})r$.
\vspace{0.2cm}

\noindent $\blacktriangledown$ The first condition is:
\vspace{-0.4cm}

\begin{small}
\begin{align*}
&\mu_1\mu_2 - \mu_0\mu_3 =\ (rb_1+1)^2(BR+1)^2\Big[\ (-Rb_2+rb_1)^2\overline{\nu}_{11}\overline{\nu}_{21}+(-Rb_2+rb_1)(R-r)\overline{\nu}_{11}\overline{\nu}_{22} + (R-r)(-Rb_2+rb_1)\overline{\nu}_{12}\overline{\nu}_{21} +(R-r)^2\\ &\overline{\nu}_{12}\overline{\nu}_{22} + (b_1-b_2)(-Rb_2+rb_1)\overline{\nu}_{13}\overline{\nu}_{21} + (b_1-b_2)^2\overline{\nu}_{13}\overline{\nu}_{23} + (-Rb_2 + rb_1)(b_1-b_2)k_{2}b_1(k_5+k_6+k_{13})(BR+1)^2(rb_1+1)(k_{18}+k_{19})\\ 
& + (b_1-b_2)(R-r)B(k_{21}+k_{22}+k_{29})b_1(rb_1+1)(BR+1)\Big((k_5+k_6+k_{13})(k_{18}+k_{2})(-Rb_2+rb_1) + (rb_1+1)(k_{19}+k_{2}+k_{3})\\ 
&(k_{18}+k_{19}+k_{2}+k_{3})\Big) + (b_1-b_2)(R-r)Bk_{18}(k_{21}+k_{22}+k_{29})b_1(rb_1+1)^2(k_{2}+k_{3})(BR+1)\Big]\ > 0 \\& \text{if and only if}\ \ b_1>b_2\ \ \text{and}\ \ r<R<\frac{b_1}{b_2}r. &&
\end{align*}
\end{small}

\noindent $\blacktriangledown$ The second condition is:
\begin{small}
\begin{align*}
&\mu_1\mu_2\mu_3 - \mu_1^2\mu_4 - \mu_0\mu_3^2 = (rb_1+1)^3(BR+1)^3b_1\Big[(-R b_2 +rb_1 )^{3} b_1^{2} (k_5+k_6+k_{13})^{3} (B R +1)^{2} \Big(k_{18}^{2} B R k_{3} +(R-r)(k_{18} B (k_{21}+k_{22} \\& +k_{29}) k_{19} + k_{18} B (k_{21} + k_{22}+k_{29}) k_{3} + B (k_{21}+k_{22}+k_{29}) k_{19}^{2}  + B (k_{21}+k_{22}+k_{29}) k_{19} k_{3} + B (k_{21}+k_{22}+k_{29}) k_{3}^{2}) + 2 k_{18} B k_{19} R k_{3} \\& + k_{18} B R k_{3}^{2} + B k_{19}^{2} R k_{3} + B k_{19} R k_{3}^{2} + k_{18}^{2} k_{3} + 2 k_{18} k_{19} k_{3} + k_{18} k_{3}^{2} + k_{19}^{2} k_{3} + k_{19} k_{3}^{2}\Big) + (-Rb_2+rb_1)^2(R-r)\overline{\nu}_{11}\overline{\nu}_{22}\overline{\nu}_{31} + (R b_2 - b_1 r )^{2} \\& (R - r) B (k_{21}+k_{22}+k_{29}) b_1^{2} (b_1 r +1) (k_5+k_6+k_{13})^{2} (B R +1) \Big(k_{18}^{2} B R k_{3} + (R-r)(k_{18} B (k_{21}+k_{22}+k_{29})
k_{19} + k_{18} B (k_{21}+k_{22} + k_{29})\\& k_{3} + B (k_{21}+k_{22}+k_{29}) k_{19}^{2} + B (k_{21}+k_{22}+k_{29})  k_{3}^{2}) + 2 k_{18} B k_{19} R
k_{3} + k_{18} B R k_{3}^{2}
 + B k_{19}^{2} R k_{3} + B k_{19} R k_{3}^{2} + k_{18}^{2} k_{3} + 2 k_{18} k_{19} k_{3} + k_{18} k_{3}^{2}\\& + k_{19}^{2} k_{3} + k_{19} k_{3}^{2}\Big) + (-Rb_2+rb_1)(R-r)^2\overline{\nu}_{12}\overline{\nu}_{22}\overline{\nu}_{31} + (R b_2 -b_1 r )^{2} 
(b_1 - b_2) (k_{18}+ k_{19} + k_{2} + k_{3}) (b_1 r +1) (B R +1)^{2} b_1 (k_5+k_6+k_{13})^{2}\\& \Big(k_{18}^{2} B R k_3 + (R-r)(k_{18}B (k_{21}+k_{22}+k_{29})k_{19} + k_{18}B (k_{21}+k_{22}+k_{29})k_{3} + B (k_{21}+k_{22} +k_{29}) k_{19}^{2} + B (k_{21}+k_{22}+k_{29}) k_{3}^{2} ) + 2 k_{18}B \\& k_{19} R
k_{3} + k_{18}B R k_{3}^{2} 
+B k_{19}^{2} R k_{3} + B k_{19} R k_{3}^{2} + k_{18}^{2} k_{3} + 2k_{18}k_{19} k_{3} + k_{18}k_{3}^{2} +
k_{19}^{2} k_{3} + k_{19} k_{3}^{2}\Big) + (R - r)^{3} B^{3} (k_{21}+k_{22}+k_{29})^{3} b_1^{2}(b_1 r +1)^{2}\\& k_{19} \Big(k_{2}(k_5+k_6+k_{13})(-Rb_2+rb_1) + r k_{19}
k_{2} b_1 + r k_{19}
b_1 k_{3} + r k_{2}^{2} b_1 + 2r k_{2} b_1 k_{3} + r b_1 
k_{3}^{2} + k_{19} k_{2} + k_{19} k_{3} +
k_{2}^{2} + 2 k_{2} k_{3} + k_{3}^{2}\Big) + (-R b_2\\& +b_1 r ) (b_1 - b_2 )^{2} (k_{18} + k_{19} + k_{2} + k_{3}) (b_1 r +1)^{2}(B R +1)^{2}(k_5+k_6+k_{13}) \Big(k_{18}^{2} B R k_{2} k_{3} + k_{18}^{2} B R k_{3}^{2}+(R-r)(k_{18} B (k_{21}+k_{22}+k_{29})\\& k_{19}k_{2} + k_{18} B (k_{21}+k_{22}+k_{29}) k_{2} k_{3} + k_{18} B (k_{21}+k_{22}+k_{29}) k_{3}^{2} + B(k_{21}+k_{22}+k_{29}) k_{19}^{2} k_{2}) + 2 k_{18} B k_{19} R k_{2} k_{3} + 2 k_{18} B k_{19} R k_{3}^{2} + B k_{19}^{2} R\\& k_{2} k_{3} + B k_{19}^{2} R k_{3}^{2} + k_{18}^{2} k_{2} k_{3} + k_{18}^{2}
k_{3}^{2} + 2 k_{18} k_{19} k_{2} k_{3} + 2k_{18} k_{19} k_{3}^{2} + k_{19}^{2} k_{2}k_{3} + k_{19}^{2} k_{3}^{2}\Big)+ (-Rb_2+rb_1)^2(b_1-b_2)k_{2}b_1(k_5+k_6+k_{13})\\&(BR+1)^2(rb_1+1)(k_{18}+k_{19})\overline{\nu}_{31} + (-Rb_2+rb_1)(R-r)^2(b_1-b_2)B^2(k_{21}+k_{22}+k_{29})^2 b_1(k_5+k_6+k_{13})(rb_1+1)^2(BR+1)(k_{18}+k_{19}\\& +k_{2}+k_{3})(k_{19}^2+k_{19}k_{2}+k_{2}k_{3}+k_{3}^2) + (-Rb_{2}+rb_1)(R-r)(b_1-b_2)B(k_{21} + k_{22} + k_{29}) b_1(rb_1+1)(BR+1)\Big((R-r)B(k_{21}+k_{22}+k_{29})\\&(k_5+k_6+k_{13})^2(k_{19}+k_{3})(k_{18}+k_{2})(-Rb_2+rb_1) + (k_5+k_6+k_{13})^2k_{3}(k_{18}+k_{19}) (BR+1)(k_{18}+k_{2})(-Rb_2+rb_1) + (k_5+k_6+k_{13})k_{3}\\& (k_{18}+k_{19})(BR+1)(rb_1+1)(k_{19}+k_{2}+k_{3})(k_{18}+k_{19}+k_{2}+k_{3})\Big) + (-Rb_2+rb_1) (R-r)(b_1-b_2) Bk_{18}(k_{21}+k_{22}+k_{29})b_1(rb_1+1)^2\\& (k_{2}+k_{3})(BR+1)\overline{\nu}_{31} + (-Rb_2+rb_1)^2(R-r)\overline{\nu}_{11}\overline{\nu}_{21}\overline{\nu}_{32} + (-Rb_2+rb_1)(R-r)^2(\overline{\nu}_{11}\overline{\nu}_{22}\overline{\nu}_{32} + \overline{\nu}_{12}\overline{\nu}_{21}\overline{\nu}_{32}) + (-Rb_2+rb_1)(b_1-b_2)\\& (R-r)\overline{\nu}_{13}\overline{\nu}_{21}\overline{\nu}_{32} + (R-r)(b_1-b_2)^2\overline{\nu}_{13}\overline{\nu}_{23}\overline{\nu}_{32} + (-Rb_2+rb_1) (R-r)(b_1-b_2)k_{2}b_1(k_{5}+k_{6} + k_{13})(BR+1)^2(rb_1+1)(k_{18}+k_{19})\overline{\nu}_{32}\\& + (b_1-b_2)(R-r)^2B(k_{21}+k_{22}+k_{29})b_1(rb_1+1)(BR+1)\overline{\nu}_{32}\Big((k_{5}+k_{6}+k_{13}) (k_{18}+k_{2})(-Rb_2+rb_1) + (rb_1+1)(k_{19}+k_{2}+k_{3})\\& (k_{18}+k_{19}+k_{2}+k_{3})\Big) + (b_1-b_2)(R-r)^2Bk_{18}(k_{21}+k_{22}+k_{29})b_1(rb_1+1)^2 (k_{2}+k_{3})(BR+1)\overline{\nu}_{32}\Big] > 0\\& \text{if and only if}\ \ b_{1} > b_{2}\ \ \text{and}\ \ r<R<\frac{b_{1}}{b_{2}}r. &&
\end{align*}
\end{small}


\begin{thebibliography}{99}

\bibitem{acr} Adimy, M.; Crauste, F.; Ruan, S. A mathematical
study of the hematopoiesis process with applications to chronic myelogenous leukemia. {\em SIAM J. Appl. Math.} \textbf{2005}, 65(4): 1328--1352.

\bibitem{af} Afenya, E. Mathematical models for cancer and
their relevant insights. In \emph{Handbook of Cancer Models With Applications} (W.Y. Tan $\&$ L. Hanin eds). Hackensack, NJ, USA: World Scientific, \textbf{2008}.

\bibitem{aaa} Alenzi, F.Q.; Alenazi, B.Q.; Ahmad, S.Y.; Salem, M.L.; Al-Jabri, A.A.; Wyse, R.K.H. The haemopoietic stem cell: between apoptosis and self renewal. {\em Yale J. Biol. Med}. \textbf{2009}, 82(1): 7–18.

\bibitem{am} Andersen, L.K.; Mackey, M.C. Resonance in periodic chemotherapy: a case study of acute myelogenous leukemia. \emph{J. Theor. Biol.} \textbf{2001}, 209(1): 113--130.

\bibitem{aoh} Arber, D.A.; Orazi, A.; Hasserjian, R.P.; Borowitz, M.J.; Calvo, K.R.; Kvasnicka, H.M.; Wang, S.A.; Bagg, A.; Barbui, T.; Branford, S.; Bueso-Ramos, C.E.; Cortes, J.E.; Dal Cin, P.; DiNardo, C.D.; Dombret, H.; Duncavage, E.J.; Ebert, B.L.; Estey, E.H.; Facchetti, F.; Foucar, K.; Gangat, N.; Gianelli, U.; Godley, L.A.; Gökbuget, N.; Gotlib, J.; Hellström-Lindberg, E.; Hobbs, G.S.; Hoffman, R.; Jabbour, E.J.; Kiladjian, J.J.; Larson, R.A.; Le Beau, M.M.; Loh, M.L.C.; Löwenberg, B.; Macintyre, E.; Malcovati, L.; Mullighan, C.G.; Niemeyer, C.; Odenike, O.M.; Ogawa, S.; Orfao, A.; Papaemmanuil, E.; Passamonti, F.; Porkka, K.; Pui, C.H.; Radich, J.P.; Reiter, A.; Rozman, M.; Rudelius, M.; Savona, M.R.; Schiffer, C.A.; Schmitt-Graeff, A.; Shimamura, A.; Sierra, J.; Stock, W.A.; Stone, R.M.; Tallman, M.S.; Thiele, J.; Tien, H.F.; Tzankov, A.; Vannucchi, A.M.; Vyas, P.; Wei, A.H.; Weinberg, O.K.; Wierzbowska, A.; Cazzola, M.; Döhner, H.; Tefferi, A. International consensus classification of myeloid neoplasms and acute leukemias: integrating morphologic, clinical, and genomic data. {\em Blood} \textbf{2022}, 140(11): 1200–1228.

\bibitem{agd} Arbore, D.R.; Galdean, S.M.; Dima, D.; Rus, I.; Kegyes, D.; Ababei, R.G.; Dragancea, D.; Tomai, R.A.; Trifa, A.P.; Tomuleasa, C. COVID-19 Impact on chronic myeloid leukemia patients. {\em J. Pers. Med.} \textbf{2022}, 12(11): 1886. 

\bibitem{bhm} Badralexi, I.; Halanay, A.-D.; Mghames, R. Stability analysis of equilibria for a model of maintenance therapy in acute lymphoblastic leukemia. {\em Mathematics} \textbf{2022}, 10(3): 313.

\bibitem{bbf} Barile, M.; Busch, K.; Fanti, A.-K.; Greco, A.; Wang, X.; Oguro, H.; Zhang, Q.; Morrison, S.J.; Rodewald, H.-R.; Höfer, T. Hematopoietic stem cells self-renew symmetrically or gradually proceed to differentiation. {\em BioRxiv} \textbf{2020}, 2020-2008.

\bibitem{bpmr} Bianca, C.; Pennisi, M.; Motta, S.; Ragusa,
M.A. Immune system network and cancer vaccine. \emph{AIP Conf. Proc.} \textbf{2011}, 1389(1): 945--948.

\bibitem{bppr} Bianca, C.; Pappalardo, F.; Pennisi, M.;
Ragusa, M.A. Persistence analysis in a Kolmogorov-type model for
cancer-immune system competition. \emph{AIP Conf. Proc.} \textbf{2013}, 1558(1): 1797--1800.

\bibitem{bes} Breccia, M.; Efficace F.; Scalzulli, E.; Ciotti, G.; Maestrini, G.; Colafigli, G.; Martelli, M. Measuring prognosis in chronic myeloid leukemia: what’s new? {\em Expert Rev. Hematol}. \textbf{2021}, 14(6): 577–585.

\bibitem{c} Carr, J.H. Clinical Hematology Atlas-E-Book. {\em Elsevier Health Sciences}, \textbf{2021}.

\bibitem{czc} Cheng, H.; Zheng, Z.; Cheng, T. New paradigms on hematopoietic stem cell differentiation. {\em Protein} {\em Cell} \textbf{2020}, 11(1): 34-44.

\bibitem{cdc} Cisneros, T.; Dillard, D.; Castro, M.; Arredondo-Guerrero, J.; Krams, S.; Esquivel, C.; Martinez, O. The role of natural killer cells in recognition and killing of stem cells and stem cell-derived hepatoblasts. {\em Am. J. Transplant}. \textbf{2017}, 17(Suppl. 3): 115.

\bibitem{cld} Clapp, G.; Levy, D. A review of mathematical models for leukemia and lymphoma. \emph{Drug Discov. Today Dis. Models}.
\textbf{2015}, 16: 1--6.

\bibitem{cl} Coddington, E.A.; Levinson, N. Theory of Ordinary Differential Equations. {\em New Delhi, India: Tata McGraw-Hill}, \textbf{1972}.

\bibitem{cm} Colijn, C.; Mackey, M.C. A mathematical model of
hematopoiesis-I. Periodic chronic myelogenous leukemia. {\em J. Theor. Biol.} \textbf{2005}, 237(2): 117--132.

\bibitem{cp} Cucuianu, A.; Precup, R. A hypothetical-mathematical model of acute myeloid leukaemia pathogenesis.
{\em Comput. Math. Methods Med}. \textbf{2010}, 11(1): 49--65.

\bibitem{dkll} DeConde, R.; Kim, P.S.; Levy, D.; Lee, P.P.
Post-transplantation dynamics of the immune response to chronic myelogenous leukemia. \emph{J. Theor. Biol.} \textbf{2005}, 236(1): 39--59.

\bibitem{dm} Dingli, D.; Michor, F. Successful therapy must eradicate cancer stem cells. {\em Stem Cells} \textbf{2006}, 24(12): 2603--2610.

\bibitem{dtm} Dingli, D.; Traulsen, A.; Michor, F. (A)symmetric stem cell replication and cancer. {\em PLoS Comput. Biol.} \textbf{2007}, 3(4): e83.

\bibitem{ds} Djulbegovic, B.; Svetina, S. Mathematical model of acute myeloblastic leukaemia: an investigation of the relevant kinetic parameters. \emph{Cell Prolif.} \textbf{1985}, 18(3): 307--319.

\bibitem{d} Domen, J. The role of apoptosis in regulating hematopoietic stem cell numbers. {\em Apoptosis} \textbf{2001}, 6: 239--252.

\bibitem{dj} Doumic-Jauffret, M.; Kim, P.S.; Perthame, B. Stability analysis of a simplified yet complete model for chronic myelogenous leukemia. \emph{Bull. Math. Biol.} \textbf{2010}, 72: 1732--1759.

\bibitem{dbcsb} Driessens, G.; Beck, B.; Caauwe, A.; Simons, B.D.; Blanpain, C. Defining the mode of tumour growth by clonal analysis. \emph{Nature} \textbf{2012}, 488(7412): 527--530.

\bibitem{dlgll} Du, J.; Li, Z.; Gong, Y.; Lan, Y.; Liu, B. Integrative cross-species transcriptome analysis reveals earlier occurrence of myelopoiesis in pre-circulation primates compared to mice. {\em Zool. Res.} \textbf{2024}, 45(6): 1276-1286.

\bibitem{dwm} Duffy, K.R.; Wellard, C.J.; Markham, J.F.; Zhou, J.H.; Holmberg, R.; Hawkins, E.D.; Hasbold, J.; Dowling, M.R.; Hodgkin, P.D. Activation-induced B cell fates are selected by intracellular stochastic competition. {\em Science} \textbf{2012}, 335(6066): 338-341.

\bibitem{ekdhs} Enriquez-Navas, P.; Kam, Y.; Das, T.; Hassan, S.; Silva, A.; Foroutan, P.; Ruiz, E.; Martinez, G.; Minton, S.; Gillies, R.; \textit{et al.} Exploiting evolutionary principles to prolong tumor control in preclinical models of breast cancer. {\em Sci. Transl. Med.} \textbf{2016}, 8(327): 327ra24-327ra24.

\bibitem{f} Fokas, A.S.; Keller, J.B.; Clarkson, B.D. Mathematical model of granulocytopoiesis and chronic myelogenous leukemia. \emph{Cancer Res.} \textbf{1991}, 51(8): 2084--2091.

\bibitem{fm} Foley, C.; Mackey, M.C. Dynamic hematological
disease: a review. \emph{J. Math. Biol.} \textbf{2009}, 58: 285--322.

\bibitem{fdc} Foo, J.; Drummond, M.W.; Clarkson, B.; Holyoake, T.; Michor, F. Eradication of chronic myeloid leukemia stem cells: A novel mathematical model predicts no therapeutic benefit of adding G-CSF to Imatinib. {\em PLoS Comput. Biol}. \textbf{2009}, 5(9): e1000503.

\bibitem{fc} Fuchs, E.; Chen, T. A matter of life and death: self-renewal in stem cells. {\em EMBO Rep.} \textbf{2013}, 14(1): 39-48.

\bibitem{grhle} Gerlinger, M.; Rowan, A.; Horswell, S.; Larkin, J.; Endesfelder, D.; Gronroos, E.; Martinez, P.; Matthews, N.;
Stewart, A.; Tarpey, P.; \textit{et al.} Intratumor heterogeneity and branched evolution revealed by multiregion sequencing. {\em N. Engl. J. Med.} \textbf{2012}, 366(10): 883--892.

\bibitem{glp} Gómez-López, S.; Lerner, R.G.; Petritsch, C. Asymmetric cell division of stem and progenitor cells during homeostasis and cancer. {\em Cell. Mol. Life Sci.} \textbf{2014}, 71: 575-597.

\bibitem{get} Grinenko, T.; Eugster, A.; Thielecke, L.; Ramasz, B.; Krüger, A.; Dietz, S.; Glauche, I.; Gerbaulet, A.; von Bonin, M.; Basak, O.; Clevers, H.; Chavakis, T.; Wielockx, B. Hematopoietic stem cells can differentiate into restricted myeloid progenitors before cell division in mice. {\em Nat. Commun.} \textbf{2018}, 9(1): 1898.

\bibitem{h} Hehlmann, R. How I treat CML blast crisis. {\em Blood} \textbf{2012}, 120(4): 737–747. 

\bibitem{hhb} Howard, M.R.; Hamilton, P.; Britton, R. Haematology, {\em London, UK: Churchill Livingstone,} \textbf{2013}, pp. 10–120.

\bibitem{hl} Hu, Y.; Li, S. Survival regulation of leukemia stem cells. {\em Cell. Mol. Life Sci.} \textbf{2016}, 73: 1039-1050.

\bibitem{jbz} Jagannathan-Bogdan, M.; Zon, L.I. Hematopoiesis. {\em Development} \textbf{2013}, 140(12): 2463–2467.

\bibitem{j} Jaiswal, S. Clonal hematopoiesis and nonhematologic disorders. {\em Blood} \textbf{2020}, 136(14): 1606–1614.

\bibitem{jg} Jilkine, A.; Gutenkunst, R.N. Effect of dedifferentiation on time to mutation acquisition in stem cell-driven cancers. \emph{PLoS Comput. Biol.} \textbf{2014}, 10(3): e1003481.

\bibitem{jps} Jones, D.S.; Plank, M.J.; Sleeman, B.D. Differential Equations and Mathematical Biology. {\em London, UK: CRC Press}, \textbf{2010}.

\bibitem{kg} Kaplan, D.; Glass, L. Understanding Nonlinear Dynamics. {\em New York, NY: Springer}, \textbf{1995}.

\bibitem{K} Kim, P.S.; Lee, P.P.; Levy, D. Mini-transplants
for chronic myelogenous leukemia: A Modeling Perspective, {\em Biology and Control Theory: Current Challenges}; Lecture Notes in Control and Information Sciences, vol. 357. Berlin/Heidelberg, Germany: Springer, \textbf{2007}, pp. 3--20.

\bibitem{kll} Kim, P.S.; Lee, P.P.; Levy, D. Dynamics and
potential impact of the immune response to chronic myelogenous leukemia. \emph{PLoS Comput. Biol.} \textbf{2008}, 4(6): e1000095.

\bibitem{ks} Klein, A.M.; Simons, B.D. Universal patterns of stem
cell fate in cycling adult tissues. {\em Development} \textbf{2011}, 138(15): 3103-3111.

\bibitem{ko} Komarova, N.L. Mathematical modeling of cyclic treatments of chronic myeloid leukemia. \emph{Math. Biosci. Eng.} \textbf{2011}, 8(2): 289--306.

\bibitem{lv} Lee, M.; Vasioukhin, V. Cell polarity and cancer–cell and tissue polarity as a non-canonical tumor suppressor. {\em J. Cell Sci.} \textbf{2008}, 121(8): 1141-1150.

\bibitem{lms} Lin, Q.; Mao, L.; Shao, L.; Zhu, L.; Han, Q.; Zhu, H.; Jin, J.; You, L. Global, regional, and national burden of chronic myeloid leukemia, 1990–2017: a systematic analysis for the global burden of disease study 2017. {\em Front. Oncol.} \textbf{2020}, 10: 580759.

\bibitem{lksw} Lopez-Garcia, C.; Klein, A.M.; Simons, B.D.;
Winton, D.J. Intestinal stem cell replacement follows a pattern of neutral drift. {\em Science} \textbf{2010}, 330(6005): 822--825.

\bibitem{m1} Mackey, M.C. Unified hypothesis for the origin of aplastic anemia and periodic hematopoiesis. \emph{Blood} \textbf{1978}, 51(5): 941--956.

\bibitem{m} Mackey, M.C. Cell kinetic status of haematopoietic stem cells. {\em Cell Prolif.} \textbf{2001} 34(2): 71-83.

\bibitem{mg} Mackey, M.C.; Glass, L. Oscillation and chaos in physiological control systems, {\em Science} \textbf{1977}, 197(4300): 287–289.

\bibitem{mls} MacLean, A.L.; Lo Celso, C.; Stumpf, M.P.H.
Population dynamics of normal and leukaemia stem cells in the haematopoietic stem cell niche show distinct regimes where leukaemia will be controlled. \emph{J. R. Soc. Interface} \textbf{2013}, 10(81): 20120968.

\bibitem{mfs} MacLean, A.L.; Filippi, S.; Stumpf, M.P.H. The
ecology in the hematopoietic stem cell niche determines the clinical outcome in chronic myeloid leukemia. \emph{Proc. Natl. Acad. Sci. USA} \textbf{2014}, 111(10): 3883--3888.

\bibitem{ms} Marciniak-Czochra, A.; Stiehl, T. Mathematical models of hematopoietic reconstitution after stem cell transplantation. \emph{Model Based Parameter Estimation: Theory and Applications} (Bock, H.G., Carraro, T., Jaeger, W., Koerkel, S., Rannacher, R., Schloeder, J.P., eds.). {Berlin/Heidelberg: Springer}, \textbf{2013}, pp. 191--206.

\bibitem{mi} Michor, F. Mathematical models of cancer stem
cells. \emph{J. Clin. Oncol.} \textbf{2008}, 26(17): 2854--2861.

\bibitem{mhi} Michor, F.; Hughes, T.P.; Iwasa, Y.; Branford, S.; Shah, N.P.; Sawyers, C.L.; Nowak, M.A. Dynamics of chronic myeloid leukaemia. {\em Nature}, \textbf{2005}, 435(7046): 1267--1270.

\bibitem{mgj} Mitchell S.R.; Gopakumar J.; Jaiswal S. Insights into clonal hematopoiesis and its relation to cancer risk. {\em Curr. Opin. Genet. Dev.} \textbf{2021}, 66: 63–69. 

\bibitem{mta} Molina-Peña, R.; Tudon-Martinez, J.C.; Aquines-Gutiérrez, O. A mathematical model of average dynamics in a stem cell hierarchy suggests the combinatorial targeting of cancer stem cells and progenitor cells as a potential strategy against tumor growth. {\em Cancers} \textbf{2020}, 12(9): 2590.

\bibitem{mk} Morrison, S.J.; Kimble, J. Asymmetric and symmetric stem-cell divisions in development and cancer. {\em Nature} \textbf{2006}, 441(7097): 1068–1074.

\bibitem{n} Neiman, B. A Mathematical Model of Chronic Myelogenous Leukemia. {\em Oxford University: Oxford UK}, \textbf{2000}. Available online: \url{https://ora.ox.ac.uk/objects/uuid:4bb8a627-0747-4629-a9fd-a8e542409174}

\bibitem{py} Palis, J.; Yoder, M.C. Yolk-sac hematopoiesis: the first blood cells of mouse and man. {\em Exp. Hematol.} \textbf{2001}, 29(8): 927-936.

\bibitem{plg} Parajdi, L.G. Stability of the equilibria of a dynamic system modeling stem cell transplantation. {\em Ricerche Mat.} \textbf{2020}, 69(2): 579--601.

\bibitem{p} Parajdi, L.G. Analysis of Some Mathematical Models of Cell Dynamics in Hematology, {\em Cluj-Napoca: Casa C\u{a}r\c{t}ii de \c{S}tiin\c{t}\u{a}}, \textbf{2021}.

\bibitem{pp} Parajdi, L.G.; Precup, R. Analysis of a planar differential system arising from hematology. {\em Stud. Univ. Babeș-Bolyai Math.} \textbf{2018}, 63(2): 235–244.

\bibitem{ppbt} Parajdi, L.G.; Precup, R.; Bonci, E.A.; Tomuleasa, C. A mathematical model of the transition from normal hematopoiesis to the chronic and accelerated-acute stages in myeloid leukemia. {\em Mathematics} \textbf{2020}, 8(3): 376.

\bibitem{pc} Parekh, C.; Crooks, G.M. Critical differences in hematopoiesis and lymphoid development between humans and mice. {\em J. Clin. Immunol.} \textbf{2013}, 33(4): 711-715.

\bibitem{pf} Pinho, S.; Frenette, P.S. Haematopoietic stem cell activity and interactions with the niche. {\em Nat. Rev. Mol. Cell Biol.} \textbf{2019}, 20(5): 303–320.

\bibitem{pacs} Precup, R.; Arghirescu. S.; Cucuianu, A.;
Șerban, M. Mathematical modeling of cell dynamics after allogeneic bone marrow transplantation. \emph{Int. J. Biomath.} \textbf{2012}, 5(02): 1250026.

\bibitem{pst} Precup, R.; Șerban, M.A.; Trif, D. Asymptotic
stability for a model of cell dynamics after allogeneic bone marrow transplantation. \emph{Nonlinear Dyn. Syst. Theory} \textbf{2013}, 13(1): 79--92.

\bibitem{pdtsp} Precup, R.; Dima, D.; Tomuleasa, C.; Șerban, M.-A.; Parajdi, L.-G. Theoretical models of hematopoietic cell dynamics related to bone marrow transplantation. {\em Frontiers in Stem Cell and Regenerative Medicine Research, Volume 8 (Rahman, A. $\&$ Anjum, S., eds). Bentham
Science Publishers}, \textbf{2018}, pp. 202--241.

\bibitem{rch} Radulescu, I.R.; Candea, D.; Halanay, A. A
study on stability and medical implications for a complex delay model for CML with cell competition and treatment. \emph{J. Theor. Biol.} \textbf{2014}, 363: 30--40.

\bibitem{rr} Ragusa, M.A.; Russo, G. ODEs approaches in
modeling fibrosis: comment on ``towards a unified approach in the modeling of fibrosis: a review with research perspectives'' by Martine ben Amar and Carlo Bianca. \emph{Phys. Life Rev.} \textbf{2016}, 17: 112--113.

\bibitem{rso} Riether, C.; Schürch, C.M.; Ochsenbein, A.F. Regulation of hematopoietic and leukemic stem cells by the
immune system. {\em Cell Death Differ.} \textbf{2015}, 22(2): 187–198.

\bibitem{rs} Riffelmacher, T.; Simon, A.-K. Mechanistic roles of autophagy in hematopoietic differentiation. {\em FEBS J.} \textbf{2017}, 284(7): 1008-1020.

\bibitem{rl1} Rubinow, S.I.; Lebowitz, J.L. A mathematical
model of neutrophil production and control in normal man. \emph{J. Math. Biol.} \textbf{1975}, 1(3): 187--225.

\bibitem{rl2} Rubinow, S.I.; Lebowitz, J.L. A mathematical
model of the acute myeloblastic leukemic state in man. \emph{Biophys. J.} \textbf{1976}, 16(8): 897--910.

\bibitem{ssd} Savino W.; Smaniotto, S.; Dardenne, M. Hematopoiesis, Chapter in: Varela-Nieto, I., Chowen, J.A. (eds) The Growth Hormone/Insulin-like Growth Factor Axis During Development. {\em Adv. Exp. Med. Biol. Vol 567}, \textbf{2005}.

\bibitem{sk} Shahriyari, L.; Komarova, N.L. Symmetric vs. asymmetric stem cell divisions: an adaptation against cancer? \emph{PloS one} \textbf{2013}, 8(10): e76195.

\bibitem{sh} Smith, H. Dynamics of competition. {\em Lecture Notes in Mathematics. Springer}, \textbf{1999}, pp. 191--240.

\bibitem{sfseb} Snippert, H.J.; Van Der Flier, L.G.; Sato, T.;
Van Es, J.H.; Van Den Born, M.; Kroon-Veenboer, C.; Barker, N.; Klein, A.M.; Van Rheenen, J.; Simons, B.D.; Clevers, H. Intestinal crypt homeostasis results from neutral competition between symmetrically dividing Lgr5 stem cells. {\em Cell} \textbf{2010}, 143(1): 134--144.

\bibitem{s} Stiehl, T.; Marciniak-Czochra, A. Mathematical modeling of leukemogenesis and cancer stem cell dynamics. \emph{Math. Mod. Nat. Phenom.} \textbf{2012}, 7(1): 166--202.

\bibitem{sh1} Stiehl, T.; Ho, A.D.; Marciniak-Czochra, A. The
impact of CD34+ cell dose on engraftment after SCTs: personalized estimates based on mathematical modeling. \emph{Bone Marrow Transplant.} \textbf{2014}, 49(1): 30--37.

\bibitem{vi} Vincent, P.C.; Rutzen-Loesevitz, L.; Tibken, B.;
Heinze, B.; Hofer, E.P.; Fliedner, T.M. Relapse in chronic myeloid leukemia after bone marrow transplantation: biomathematical modeling as a new approach to understanding pathogenesis. \emph{Stem Cells} \textbf{1999}, 17(1): 9--17.

\bibitem{vtb} Vivier, E.; Tomasello, E.; Baratin, M.; Walzer, T.; Ugolini, S. Functions of natural killer cells. {\em Nat. Immunol.} \textbf{2008}, 9(5): 503–510.

\bibitem{vgh} Vuelta, E.; García-Tuñón, I.; Hernández-Carabias, P.; Méndez, L.; Sánchez-Martín, M. Future approaches for treating chronic myeloid leukemia: CRISPR therapy. {\em Biology} \textbf{2021}, 10(2): 118.

\bibitem{wlowb} Wilson, A.; Laurenti, E.; Oser, G.; van der
Wath, R.C.; Blanco-Bose, W.; Jaworski, M.; Offner, S.; Dunant, C.F.; Eshkind, L.; Bockamp, E.; Lió, P.; Macdonald, H.R.; Trumpp, A. Hematopoietic stem cells reversibly switch from dormancy to self-renewal during homeostasis and repair. \emph{Cell} \textbf{2008}, 135(6): 1118--1129.

\bibitem{ya} Yang, J.; Axelrod, D.E.; Komarova, N.L. Determining the control networks regulating stem cell lineages in colonic crypts. \emph{J. Theor. Biol.} \textbf{2017}, 429(1): 190--203.

\bibitem{y} Young, N.S.; Gerson, S.L.; High, K.A. Clinical hematology. {\em Mosby.} Philadelphia, PA, USA, \textbf{2006}.

\bibitem{zdm} Zhang, S.; Davidson, D.D.; Cheng, L. Conceptual evolution in cancer biology. {\em Molecular Genetic Pathology. Springer. Totowa, NJ: Humana Press,} \textbf{2008}, pp. 185-208.

\end{thebibliography}
\end{document}